\documentclass[11pt]{amsart}
\pdfoutput=1
\usepackage{latexsym, graphicx, epsfig, amsmath, amsfonts,amssymb}
\usepackage{multirow}
 \usepackage{diagbox}
\usepackage{color}
 \usepackage[percent]{overpic}
\usepackage{algorithm}
\usepackage{algorithmicx}
\usepackage{algpseudocode}
\floatname{algorithm}{Algorithm}
\usepackage{threeparttable}
\usepackage{booktabs}

\usepackage{epstopdf,mathtools,bbm}
\usepackage[left=1.25in,right=1.25in,top=1in]{geometry}

\def\mb{\mathbf}

\newcounter{Rownumber}

\newtheorem{theorem}{Theorem}[section]

\newtheorem{assumption}[theorem]{Assumption}

\title{Failure-informed adaptive sampling for PINNs}

\author{Zhiwei Gao}
\thanks{School of Mathematics, Southeast University, Nanjing 210096, China.}

\author{Liang Yan}
\thanks{School of Mathematics, Southeast University, Nanjing 210096; Nanjing Center for Applied Mathematics,  Nanjing 211135,  China. Email: yanliang@seu.edu.cn. L. Yan is supported by NSF of China (No.11771081)}

\author{Tao Zhou}
\thanks{LSEC, Institute of Computational Mathematics, Academy of Mathematics and Systems
Science, Chinese Academy of Sciences, Beijing 100190, China. Email: tzhou@lsec.cc.ac.cn. T. Zhou is partially supported  by the National Key R$\&$D Program of China(No. 2020YFA0712000), the NSF of China (under grant numbers 11822111, 11688101and 11731006),  the Strategic Priority Research Program of Chinese Academy of Sciences (No.  XDA25000404) and youth innovation promotion association (CAS)}

\begin{document}
\graphicspath{figures/}
\maketitle

\begin{abstract}
Physics-informed neural networks (PINNs) have emerged as an effective technique for solving PDEs in a wide range of domains. It is noticed, however, the performance of PINNs can vary dramatically with different sampling procedures. For instance, a fixed set of (prior chosen) training points may fail to capture the effective solution region (especially for problems with singularities). To overcome this issue, we present in this work an adaptive strategy, termed the failure-informed PINNs (FI-PINNs), which is inspired by the viewpoint of reliability analysis. The key idea is to define an effective failure probability based on the residual, and then, with the aim of placing more samples in the failure region, the FI-PINNs employs a failure-informed enrichment technique to adaptively add new collocation points to the training set, such that the numerical accuracy is dramatically improved. In short, similar as adaptive finite element methods, the proposed FI-PINNs adopts the failure probability as the posterior error indicator to generate new training points. We prove rigorous error bounds of FI-PINNs and illustrate its performance through several problems.
\end{abstract}

%\begin{keyword}
%Physics-Informed Neural Networks,  failure probability, adaptive sampling, deep learning.
%\end{keyword}

\pagestyle{myheadings}
\thispagestyle{plain}
\markboth{Z. Gao, L. Yan and T. Zhou}
{FI-PINNS}

\section{Introduction}
Partial differential equations (PDEs) are important tools for modeling many real-world phenomena.  Traditional numeric solvers such as finite difference method and finite element method have been widely used to solve PDEs for many decades.
However,  for high-dimensional problems, the computational cost will be prohibitively expensive when applying those methods due to the dramatic increase in the number of grids and thus the curse of dimensionality is inevitable.  Driven by their well-documented performance in fields of machine learning, particularly deep learning is being increasing used in scientific computing. As a result, the area of scientific machine learning has emerged, see e.g., \cite{E_CICP,Edeep18,raissi2019physics,sirignano2018dgm,zang2020weak}.

Physics-Informed neural networks (PINNs) \cite{lu2021deepxde,raissi2019physics} is one of the popular machine learning methods for solving PDEs with deep neural networks. PINNs adopt automatic differentiation to solve PDEs by penalizing the PDE in the loss function at a random set of points in the domain of interest. PINNs have been successfully applied to simulate a variety of forward and inverse problems for PDEs, see \cite{wang2021understanding,karniadakis2021physics,lu2021deepxde,GWZ_2022,GWZ_2022_JCP} and references therein. Being a mesh less approach, a key issue in PINNs is how to choose the training points. Obviously, for PDEs with complex solution structures, a fixed set of (prior chosen) training points may fail to capture the effective solution region \cite{krishnapriyan2021characterizing,wang2021understanding,wang2022and}. Another limitation is that when implementing PINNs to PDEs in an unbounded domain, due to the lack of prior knowledge about the PDE solutions, choosing an effective training set can be challenging. To address these issues, some adaptive sampling strategies have been proposed \cite{daw2022rethinking,gao2021active,lu2021deepxde,subramanian2022adaptive,tang2021deep,FZZ_2022}, see also a recent review \cite{wu2022comprehensive}. Among others, the most common strategy is the residual-based adaptive refinement method (RAR) in \cite{lu2021deepxde}. The RAR method provides an insight to adaptively select samples based on residual errors. However, the adaptive samples in RAR are chosen based on a large set of candidate points (either grid-based or sampled with a prior distribution), making it not effective for high dimensional problems.
Another interesting idea is the so-called deep adaptive sampling (DAS) method \cite{tang2021deep}, in which one first train a generative model based on the residual, and then the generative model can be used to generate new training points. While DAS is certainly a more general approach, however, the main issue for DAS is that the complexity for training a generative model can be comparable for solving the PDEs.

Motivated by the above discussions, we present in this work the failure-informed PINNs (FI-PINNs) that can effectively generate adaptive training points. Our approach is inspired by the viewpoint of reliability analysis. The key idea is to define an effective failure probability based on the residual, and then, with the aim of placing more samples in the failure region, the FI-PINNs employs a failure-informed enrichment technique to adaptively add new collocation points to the training set, in such a way that the numerical accuracy is dramatically improved.  We summarize the main features of FI-PINNs as follows:
\vspace{0.25cm}
\begin{itemize}
\item In each iteration, based on the residual, we define the failure probability as an error indicator. This is similar as in the adaptive finite element methods, where one drive posterior error indicator for adjusting the meshes.
\vspace{0.25cm}
\item In each iteration, we use the simple (truncated) Gaussian model to estimate the failure probability and to generate new training points, and this can be done in a very efficient way (compared to RAR and DAS).
\vspace{0.25cm}
\item We prove rigorous error bound for FI-PINNs in terms of the error tolerance and failure probability tolerance.
\vspace{0.25cm}
\item We test FI-PINNs for various PDE problems, which include PDEs with singular solutions,  PDEs in unbounded domains,  and time-dependent PDEs. It is shown that for all the test problems, FI-PINNs can effectively capture the solution structure.
\vspace{0.25cm}
\item Our adaptive sampling strategy, especially the important sampling idea, may be extended to other formulations such as the deep Ritz method \cite{Edeep18} and weak adversarial networks \cite{zang2020weak}.
 \end{itemize}
\vspace{0.25cm}
The rest of  this paper is organized as follows. In the next section, we briefly introduces the basic knowledge of PINNs. In Section \ref{failure_probability_criteria} we present FI-PINNs and propose the adaptive sampling strategy.  In Section \ref{analysis_of_convergence}, we present the convergence analysis of FI-PINNs. Finally, we provide numerical experiments in Section \ref{Numerical_experiments}, and this is followed by some concluding remarks in Section \ref{sec:con}.

\section{Preliminaries of PINNs}\label{Preliminaries_of_pinns}

We briefly review physics-informed neural networks (PINNs). Let $\Omega \in \mathbb{R}^d$ be a spatial domain, and we denote by $\mb{x}\in \Omega$ the spatial variable. Consider the following partial differential equations:
\begin{equation}
\begin{split} \label{nonlinear-pde}
&\mathcal{A}(\mb{x};u(\mb{x})) = 0, \quad \mb{x}\in \Omega,\\
&\mathcal{B}(\mb{x};u(\mb{x})) = 0, \quad \mb{x}\in \partial  \Omega,\\
% &u(\mb{x}, 0) = g(\mb{x}), \quad \mb{x}\in \Omega
\end{split}
\end{equation}
where $\mathcal{A}$ is a linear or non-linear differential operator, $\mathcal{B}$ is the boundary operator, and $u(\mb{x})$ is the unknown solution.

The basic idea of PINNs is to use a deep neural network (DNN)  $u(\mb{x};\theta)$ with parameters $\theta$ to approximate the unknown solution $u(\mb{x})$.  The PDE solution is then obtained by choosing the best parameters for the following soft constrained optimization problem:
\begin{equation}\label{loss-function}
\min_{\theta \in \Theta}\mathcal{L}(\theta) = \min_{\theta\in\Theta} \mathcal{L}_{c}(\theta)  + \lambda\mathcal{L}_{b}(\theta),
\end{equation}
where  $\Theta$ is the parameter space and $\lambda$ is a penalty factor that balances the PDE loss $\mathcal{L}_{c}(\theta)$ and the boundary loss $\mathcal{L}_{b}(\theta)$. A common choice for $\mathcal{L}(\theta)$ is the $L^2$ loss, i.e., $$\mathcal{L}_{c}(\theta)=\|r(\mb{x};\theta)\|_{2, \Omega}^{2}, \quad \mathcal{L}_{b}(\theta)=\|b(\mb{x};\theta)\|_{2,  \partial\Omega}^{2},$$
where $r(\mb{x};\theta) =\mathcal{A}(\mb{x}; u(\mb{x};\theta))$, and $b(\mb{x};\theta)=\mathcal{B}(\mb{x}; u(\mb{x};\theta))$ measure how well $u(\mb{x}, \theta)$ satisfies the PDE and the boundary operators, respectively. The above $L^2$-norm is defined as $\|u\|_{2, \Omega}^{2} = \int_{\Omega } u(\mb{x})^{2}\omega(\mb{x})d\mb{x}$, where $\omega(\mb{x})$ is a prior distribution and is usually set to be 1 when the problem domain is bounded.
 \begin{figure}[H]
    \centering
    \includegraphics[width = 0.7\textwidth]{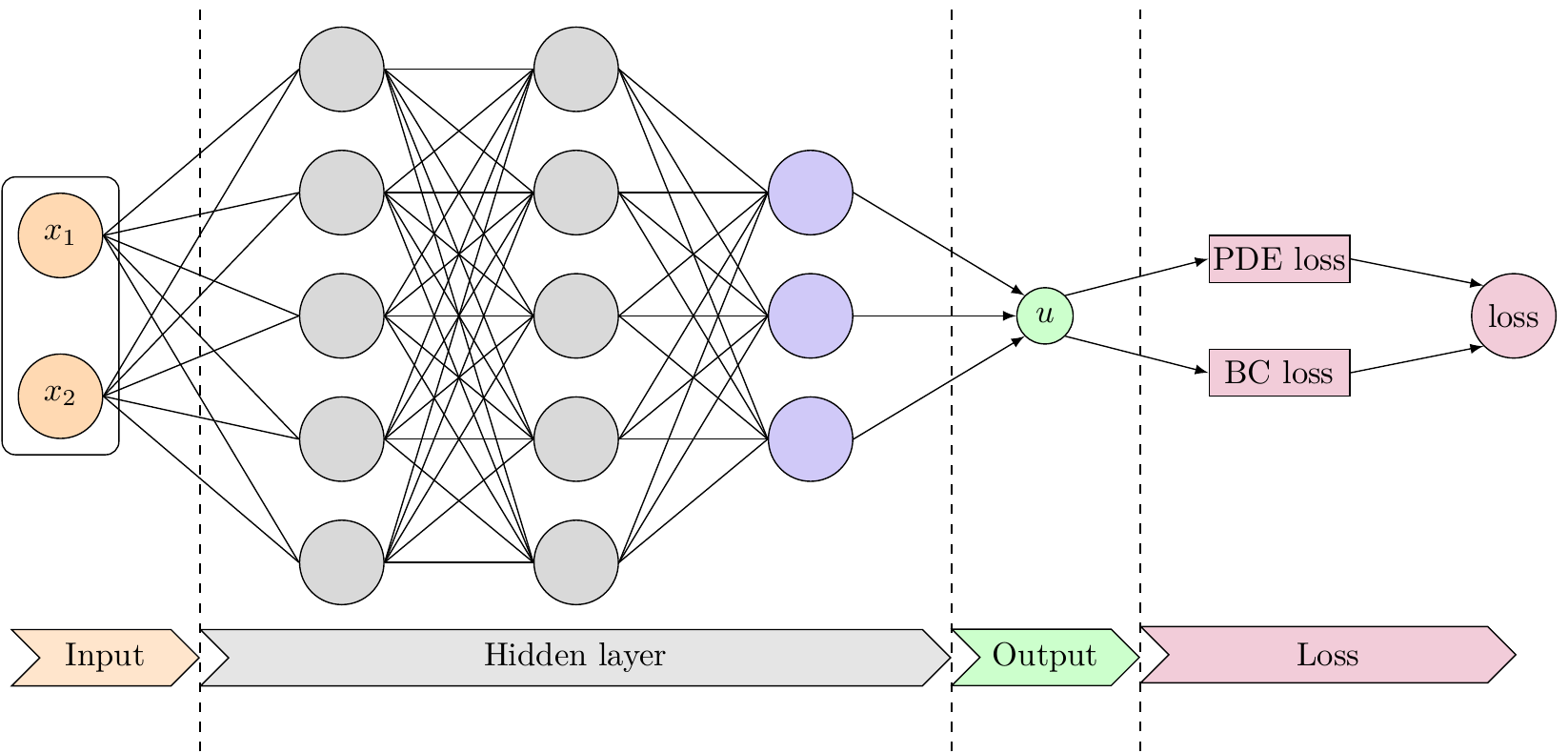}
    \caption{The framework of PINNs}
    \label{pinn_framework}
\end{figure}

We present in Fig.\ref{pinn_framework} the general PINNs framework. In practice, the loss functional $\mathcal{L}(\theta)$ should be discretized
by giving a prior distributed training dataset, which consists of the discrete points (collocation points) $\mathcal{D}_{c} = \{\mb{x}_{i}^{c}\}_{i=1}^{N_{c}}$ and the boundary points $\mathcal{D}_{b} = \{\mb{x}_{i}^{b}\}_{i=1}^{N_{b}}$.
We can then consider the discrete loss function
  \begin{equation}
    \label{optimal_parameters}
\hat{\mathcal{L}}(\theta) = \hat{\mathcal{L}}_{c}(\theta) + \lambda\hat{\mathcal{L}}_{b}(\theta),
  \end{equation}
 where
 \begin{equation}\label{discrete_loss_function}
\hat{\mathcal{L}}_{c}(\theta) = \frac{1}{N_{c}}\sum_{i = 1}^{N_{c}}\left|r(\mb{x}^{c}_{i};\theta)\right|^{2}, \quad
\hat{\mathcal{L}}_{b}(\theta)  = \frac{1}{N_{b}}\sum_{i=1}^{N_{b}}\left|b(\mb{x}^{b}_{i};\theta) \right|^{2}.
\end{equation}
It should be noted that the penalty factor $\lambda$ is often chosen as $\lambda=1$  in most cases, but it can be carefully tuned to achieve better performance \cite{wang2021understanding,wang2022and}. The discrete loss function (\ref{optimal_parameters}) can then be minimized using stochastic gradient-based algorithms \cite{bottou2012stochastic,kingma2014adam}.

Although PINNs have been successfully applied to simulate a variety of forward and inverse problems for PDEs, they can sometimes have difficulty in converging to the true solution, as shown in several studies depicting the \textit{failure modes} of PINNs \cite{krishnapriyan2021characterizing,wang2021understanding,wang2022and}. These potential failure modes are caused by the mechanism of PINNs, which makes optimizing the loss function extremely difficult due to its complicated landscape. To address these challenges, specialized solutions have been proposed independently.

 \paragraph{Modify the loss function} Modifying the structure of the loss function has become a popular tool for training PINNs. This includes (I) applying adaptive weights for each individual point \cite{mcclenny2020self}, (II) embedding gradient information of the residual function in the total loss \cite{yu2022gradient}, and (III) combing the augmented Lagrangian relaxation to balance different loss components \cite{huang2022augmented}.

\paragraph{Construct new learning schemes}  Recent ideas along this line include the extended PINNs (XPINNs) \cite{jagtap2021extended} and sequence to sequence learning\cite{krishnapriyan2021characterizing}.

\paragraph{Adaptive sampling strategies} Adaptive sampling strategies are essential for training PINNs \cite{daw2022rethinking,gao2021active,lu2021deepxde,subramanian2022adaptive,tang2021deep,wu2022comprehensive}. As discussed in the last section, \cite{lu2021deepxde} proposes to update the training dataset by selecting points  from a uniformly distributed set with larger residual values. While \cite{tang2021deep} proposed to use the generative model, e.g., the KRnet, to generate the training points.

Our contribution belongs to the third category, i.e., we focus on developing an efficient adaptive sampling strategy for PINNs.

\section{Failure-informed PINNs (FI-PINNs)}
\label{failure_probability_criteria}
In this section, we introduce the concept of failure-informed PINNs. To this end, We shall first discuss how to establish an adaptive sampling framework from the  view  of reliability analysis, and then we shall present a self-adaptively importance sampling procedure in FI-PINNs.

To begin, we first define the so-called {\it limit-state function (LSF)} $g(\mb{x}):\Omega\rightarrow \mathbb{R}$ $$g(\mb{x}) = \mathcal{Q}(\mb{x}) - \epsilon_{r}.$$ Here, $\epsilon_r$ is a predefined maximum allowed threshold, and $\mathcal{Q}(\mb{x}):\Omega\rightarrow \mathbb{R}$ maps the domain to a quantity of interest (QoI) that characterizes the  system's performance. In PINNs, we can simply choose $\mathcal{Q}(\mb{x})= |r(\mb{x};\theta)|$ with $r(\mb{x};\theta)$ being the residual such that
\begin{equation}\label{Power_function}
g(\mb{x}) = |r(\mb{x};\theta)| - \epsilon_{r}.
\end{equation}
The failure hypersurface defined by $g(\mb{x})=0$ divides the physics domain into two subsets: the safe set $\Omega_{\mathcal{S}} = \{\mb{x}:g(\mb{x}) < 0\}$ and the failure set $\Omega_{\mathcal{F}} = \{\mb{x}:g(\mb{x}) > 0\}.$   To describe the reliability of the PINNs, we define the {\it failure probability} $P_{\mathcal{F}}$ under a prior distribution $\omega(\mb{x})$:
 \begin{equation}\label{failure_probability}
P_{\mathcal{F}} = \int_{\Omega}\omega(\mb{x})\mathbb{I}_{\Omega_{\mathcal{F}}}(\mb{x})d\mb{x}.
\end{equation}
Here $\mathbb{I}_{\Omega_{\mathcal{F}}} :\Omega \rightarrow \{0,1\}$ represents the indicator function, which takes values $\mathbb{I}_{\Omega_{\mathcal{F}}}(\mb{x})=1$ when $\mb{x} \in \Omega_{\mathcal{F}}$, and $\mathbb{I}_{\Omega_{\mathcal{F}}}(\mb{x})=0$ otherwise.

If the failure probability is smaller than a given tolerance $\epsilon_{p}$, we say that PINNs are reliable; otherwise, the failure region indicates that the system is unreliable, and the performance of PINNs should be improved. Generally speaking, as the failure region gets smaller, the failure probability decreases, meaning that the system is becoming more reliable. Motivated by this fact, we may design adaptive strategies to add new training points from the failure region $\Omega_{\mathcal{F}}$ and retrain the PINNs if $P_{\mathcal{F}}>\epsilon_{p}$. The failure probability can thus be used as a stopping indicator during the training of PINNs, similar as the posterior error indicator in adaptive finite element methods.

    \begin{figure}
        \centering
        \includegraphics[width = 0.6\textwidth]{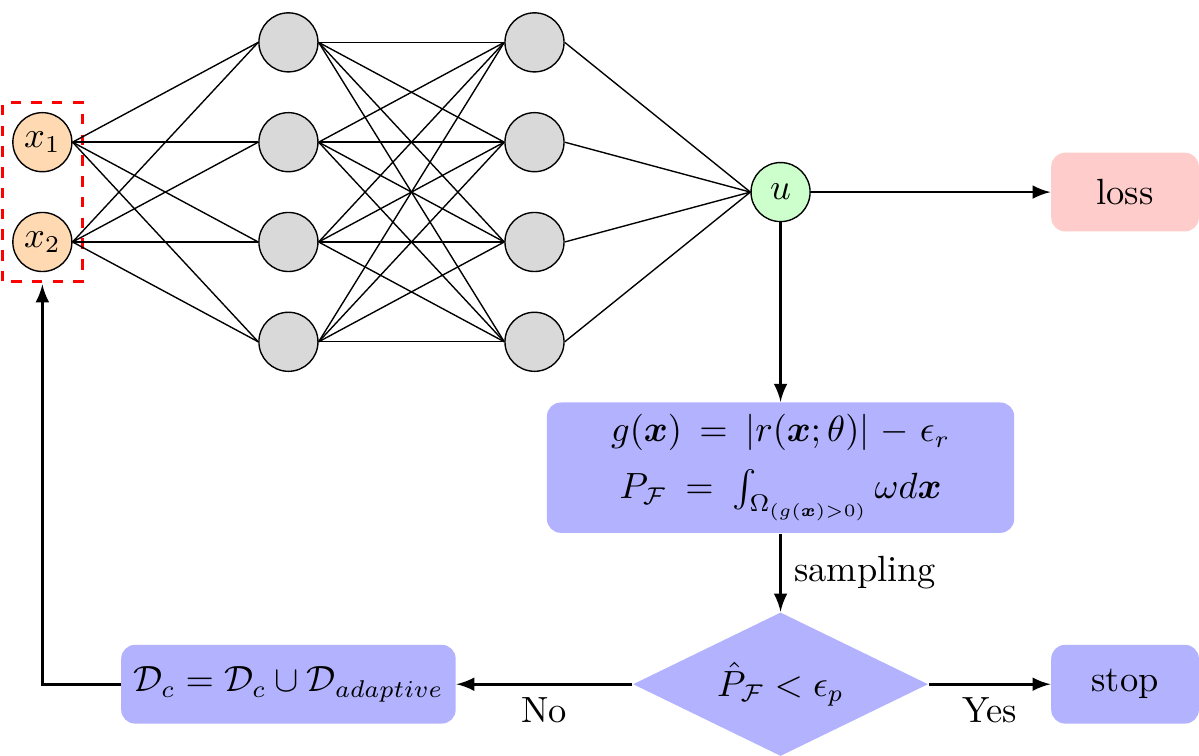}
        \caption{Workflow of FI-PINNs}
        \label{IS_PINN}
    \end{figure}

     \begin{algorithm}[t]
        \caption{Failure informed PINNs (FI-PINNs)}
        \label{Algorithm1}
        \begin{algorithmic}[1]
        \Require  ~ A DNN solution $u(\mb{x};\theta)$, boundary points $\mathcal{D}_{b}$, collocation points $\mathcal{D}_{c}$. Maximum iterations $M$,
        residual tolerance $\epsilon_{r}$ and failure probability tolerance $\epsilon_{p}$.
        \State $s\leftarrow 1$.
        \While  {$s \le M$}
        \State Train $u(\mb{x}; \theta)$ using the training dataset $\mathcal{D}_{b}, \mathcal{D}_{c}$.
        \State Set LSF $g$ using Eq.\eqref{Power_function}.
        \State Estimate the failure probability $\hat{P}_{\mathcal{F}}$ using some sampling technique.
        \If{$\hat{P}_{\mathcal{F}}<\epsilon_{p}$}
        \State \textbf{Break}
        \EndIf
        \State  Generate a new dataset $\mathcal{D}_{adaptive}$ from the failure region $\Omega_{\mathcal{F}}$.
        \State Set $\mathcal{D}_{c} = \mathcal{D}_{c}\cup \mathcal{D}_{adaptive}, s = s+1$.
        \EndWhile
        \end{algorithmic}
    \end{algorithm}

We present in Fig. \ref{IS_PINN} a general workflow of FI-PINNs, which can be used to combine with adaptive sampling strategies. The detailed algorithm of FI-PINNs is summarized in Algorithm \ref{Algorithm1}.  Notice that in Algorithm \ref{Algorithm1}, the key issue is to estimate the failure probability $P_{\mathcal{F}}$ in Eq.\eqref{failure_probability}, as the integral cannot be calculated analytically. Moreover, one needs also to design effective model to sampling the failure region in an efficient way.  These two issues will be addressed in the forthcoming subsections.

\subsection{Monte Carlo and importance sampling} \label{self adaptive_importance_sampling}

To estimate the failure probability, a nature idea is to use the Monte Carlo (MC) methods. In this case, one may generate a set of randomly distributed locations $\mathcal{S} = \{\mb{x}_{1}, \mb{x}_{2},\ldots, \mb{x}_{|\mathcal{S}|}\}$ from the prior $\omega(\mb{x})$ (e.g., the uniform distribution),
then the Monte Carlo estimator is given by $$\hat{P}_{\mathcal{F}}^{MC} = \frac{1}{|\mathcal{S}|}\sum_{\mb{x}\in \mathcal{S}}\mathbb{I}_{\Omega_{\mathcal{F}}}(\mb{x}).$$  Using this estimator, we can propose an adaptive sampling scheme.  If $\hat{P}_{\mathcal{F}}^{MC}>\epsilon_{p}$,  we add the new collocation points $\{\mb{x}_{i}\}$ which satisfy $\{\mb{x}_{i}:\mb{x}_i \in \mathcal{S}, \, g(\mb{x}_{i}) > 0\}$ into the training set and retrain the network. Notice that in each iteration, the number of the new training points $m$ may be different; it is determined by the number of points $\{\mb{x}_{i}\}\in \mathcal{S}$ that fall in the failure region $\Omega_{\mathcal{F}}$. The residual-based adaptive refinement (RAR) method \cite{lu2021deepxde} proposed to select a fixed number of $m$ new points with the largest values of LSF $g$ in $\mathcal{S}.$
%\begin{algorithm}[t]
%    \caption{FI-PINNs using Monte Carlo}
%    \label{Monte Carlo method}
%    \begin{algorithmic}[1]
%        \State Generate a set of randomly distributed locations $\mathcal{S} = \{\mb{x}_{1}, \mb{x}_{2},\ldots, \mb{x}_{|\mathcal{S}|}\}$ from the prior $\omega$ (e.g., uniform distribution).
%    \State Estimate the failure probability $\hat{P}_{\mathcal{F}}^{MC}$:
%    \begin{equation*}
%        \hat{P}_{\mathcal{F}}^{MC} = \frac{1}{|\mathcal{S}|}\sum_{\mb{x}\in \mathcal{S}}\mathbb{I}_{\Omega_{\mathcal{F}}}(\mb{x})
%    \end{equation*}
%    \State Stop if $\hat{P}_{\mathcal{F}}^{MC}<\epsilon_{p}$. Otherwise, add the $(m)$ new points $\mb{x}_{i}\in \mathcal{S}$ satisfying $\{\mb{x}_{i}:g(\mb{x}_{i}) > 0\}$, retrain the network, and go to 1.
%    \end{algorithmic}
%\end{algorithm}

%\begin{algorithm}[t]
%    \caption{FI-PINNs using importance sampling}
%    \label{importance_sampling}
%    \begin{algorithmic}[1]
%        \State Choose a set of randomly distributed samples $\mathcal{S} = \{\mb{x}_{1}, \mb{x}_{2}, \ldots, \mb{x}_{|\mathcal{S}|}\}$ from the proposal distribution $h$.
%        \State Estimate the failure probability $\hat{P}_{\mathcal{F}}^{IS}$ using Eq. \eqref{3.1.3}.
%        \State Stop if $\hat{P}_{\mathcal{F}}^{IS}<\epsilon_{p}$. Otherwise, add the $(m)$ new points $\mb{x}_{i}\in \mathcal{S}$ satisfy $\{\mb{x}_{i}:g(\mb{x}_{i}) > 0\}$, retrain the network, and go to 1.
%    \end{algorithmic}
%\end{algorithm}
Since the failure region may be relatively small compared to the entire problem domain, the above MC sampling strategy can be ineffective in generating effective samples. This is especially true when the PDEs exhibit local behaviors (e.g., in presence of sharp, or very localized features). Consequently, the sample size $|\mathcal{S}|$ is typically in the range of $\mathcal{O}(10^4\sim10^6)$, meaning that the procedure is extremely expensive.

Importance sampling (IS) can be considered to generate effective samples as it can reduce variance by selecting appropriate proposal distributions. In this case, we have
     \begin{equation}
         \label{failure_probabilty_importance_sampling}
         P_{\mathcal{F}} = \int_{\Omega}\mathbb{I}_{\Omega_{\mathcal{F}}}(\mb{x})\frac{\omega(\mb{x})}{h(\mb{x})}h(\mb{x})d\mb{x} =
          \mathbb{E}_{h}\left[\mathbb{I}_{\Omega_{\mathcal{F}}}(\mb{x})R(\mb{x})\right],
     \end{equation}
where $h(\mb{x})$ is the proposal distribution. Here $R(\mb{x})$ represents the weight function $\frac{\omega(\mb{x})}{h(\mb{x})}$ that transfers the proposal distribution $h(\mb{x})$ to the prior distribution $\omega(\mb{x})$.  By choosing a set of randomly distributed samples $\mathcal{S} = \{\mb{x}_{1}, \mb{x}_{2}, \ldots, \mb{x}_{|\mathcal{S}|}\}$ from the proposal distribution $h$, we can approximate the failure probability $P_{\mathcal{F}}$ as
     \begin{equation}
         \label{3.1.3}
\hat{P}_{\mathcal{F}}^{IS} = \frac{1}{|\mathcal{S}|}\sum_{\mb{x}\in \mathcal{S}}\mathbb{I}_{\Omega_{\mathcal{F}}}(\mb{x})R(\mb{x}).
     \end{equation}

%The FI-PINNs using importance sampling is presented in Algorithm \ref{importance_sampling}.
If the support of $h$ contains the intersection of the support of $\omega$ and the failure set, \eqref{3.1.3} gives an unbiased estimator for $P_{\mathcal{F}}$.  Theoretically, there exists an optimal proposal distribution,
 \begin{equation}
    \label{optimal_importance_distribution}
    h_{opt}(\mb{x}) = \frac{\mathbb{I}_{\Omega_{\mathcal{F}}(\mb{x})}\omega(\mb{x})}{P_{\mathcal{F}}} = \frac{\mathbb{I}_{g(\mb{x}) > 0}\omega(\mb{x})}{\int_{\Omega}\mathbb{I}_{g(\mb{x}) > 0}\omega(\mb{x})d\mb{x}},
 \end{equation}
which leads to a zero-variance estimator.  However, $h_{opt}(\mb{x})$ is not available in practice due to the normalizing constant. Note that, the optimal IS density $h_{opt}(\mb{x}) $ can be interpreted as a  ``posterior-failure" density of the collocation points given the occurrence of the failure set  $\Omega_{\mathcal{F}}$, where the indicator function $\mathbb{I}_{\Omega_{\mathcal{F}}(\mb{x})}$ is the likelihood function,  $\omega(\mb{x})$ is the prior density, and $P_{\mathcal{F}}$ is the evidence.

Although it is impossible to evaluate and sample directly from $h_{opt}$, this still serves as a guide for selecting an IS proposal distribution.  We will provide a self-adaptive importance sampling (SAIS) method for approximating the posterior-failure density $h_{opt}$ in the following section.

\subsection{Self-adaptive importance sampling (SAIS)}
\label{sec:sais}

As discussed above, the optimal  density $h_{opt}$ in \eqref{optimal_importance_distribution} is in general not available. Consequently, our strategy is to develop an adaptive procedure that begins with an initial proposal and iteratively updates the intermediate proposal using samples.  For efficiency reason, we simply choose the intermediate proposal distributions $h_{k}$ as truncated Gaussian restricted to $\Omega$, denoted as $\mathcal{N}_{T}(\mu_{k}, \Sigma_{k})$.

To establish the adaptive scheme, we first set the initial proposal distribution as  $h_{1}(\mb{x})=\omega(\mb{x})$. Then, at the $k$ step, we generate $N_{1}$ samples $\{\mb{x}_{i}^{k}\}_{i=1}^{N_{1}}$ from $h_{k}(\mb{x})$ and  sort these samples according to their LSF values in a descending order to obtain the candidate points set $\mathcal{D}_k:=\{\widetilde{\mb{x}}_{i}^{k}\}_{i=1}^{N_{1}}$.  Let $N_{p} = \lfloor p_{0}N_{1}\rfloor$ denotes the minimum number of samples used to  approximate the optimal distribution ${h}_{opt}$, here $0<p_0<1$ is a fixed parameter.  Let $N_{\eta}$ represents the number of points where $\mathcal{D}_k$ falls in the failure region. We use this number as an indicator: if $N_{\eta} <N_{p}$, this means that the intermediate proposal distribution $h_k$ needs to be refined. Otherwise, if $N_{\eta} $ is bigger than $N_{p}$, it means that the number of points is acceptable and we just use those points to approximate the optimal distribution $\hat{h}_{opt}$.

To update $h_k$ into $h_{k+1},$ we consider the truncated Gaussian model, and use the first $N_{p}$ samples to estimate the mean vector and the covariance matrix of $h_{k+1}$  as follow:
\begin{equation}
        \label{3.1.9}
        \begin{split}
        &\mu_{k+1} = \frac{1}{N_{p}}\sum_{i = 1}^{N_{p}}\widetilde{\mb{x}}_{i}^{k},\\
        &\Sigma_{k+1} = \frac{1}{N_{p} - 1}\sum_{i=1}^{N_{p}}(\widetilde{\mb{x}}_{i}^{k} - \mu_{k+1}) \otimes (\widetilde{\mb{x}}_{i}^{k} - \mu_{k+1}).
        \end{split}
\end{equation}
Notice that this procedure can be done in a very efficient way.

Notice also that once the iterative scheme stopped, the mean vector  $ \mu_{opt}$ and the  covariance matrix  $\Sigma_{opt} $  are approximated using the first $N_{p}$ samples  from the last iteration:
    \begin{equation}
        \label{optimal_sampling_center}
        \begin{split}
        \mu_{opt} &= \frac{\sum_{i=1}^{N_{p}}\widetilde{\mb{x}}_{i}\omega(\widetilde{\mb{x}}_{i})}{\sum_{i=1}^{N_{p}}\omega(\widetilde{\mb{x}}_{i})}\\
        \Sigma_{opt} &= \frac{1}{N_{p} - 1}\sum_{i=1}^{N_{p}}(\widetilde{\mb{x}}_{i} - \mu_{opt})\otimes (\widetilde{\mb{x}}_{i} - \mu_{opt}).
        \end{split}
    \end{equation}
Thus, the approximated optimal distribution $\hat{h}_{opt}(\mb{x})$ is a truncated Gaussian with mean vector $\mu_{opt}$ and covariance matrix $\Sigma_{opt}$. By generating $N_{2}$ samples from $\hat{h}_{opt}$, the failure probability can be approximated as
    \begin{equation}
        \label{failure_probability_estimation}
        \hat{P}_{\mathcal{F}}^{SAIS} = \frac{1}{N_{2}}\sum_{i=1}^{N_{2}}\frac{\omega(\mb{x}_{i})}{\hat{h}_{opt}(\mb{x}_{i})}\mathbb{I}_{\Omega_{\mathcal{F}}}(\mb{x}_{i}).
    \end{equation}

\begin{algorithm}[t]
    %    \setstretch{1.35}
 \caption{Self-adaptive importance sampling (SAIS)}
\label{Algorithm2}
\begin{algorithmic}[1]
 \Require   ~ Number of samples $N_{1}$ and $N_{2}$, the parameter $p_{0}$, the LSF $g$, the prior distribution $\omega$ and the maximum updated number $M$.
        \State $k \leftarrow 1$, set $h_1 =\omega$
         \While {$k\leq M$}
         \State Generate $N_{1}$ samples $\{\mb{x}_{i}\}_{i=1}^{N_{1}}$ from $h_k$.
          \State Sort the samples according to the corresponding LSF  in a descending order to obtain $\widetilde{\mb{x}}_{1}, \ldots, \widetilde{\mb{x}}_{N_{1}}$.
          \State Let $N_{\eta} = \max_{1\leq i\leq N_{1}}\{i|g(\widetilde{\mb{x}}_{i})>0$\} and $N_{p} = \lfloor p_{0}N_{1} \rfloor$.
           \If{$N_{\eta} < N_{p}$}
           \State  Compute $\mu_{k+1}$ and  $ \Sigma_{k+1} $  using Eq.\eqref{3.1.9}.
              \State Set $h_{k+1} = \mathcal{N}_{T}(\mu_{k+1}, \Sigma_{k+1})$.
           \State  $k  \leftarrow  k+1$.
            \Else
            \State \textbf{Break}
                \EndIf
               \EndWhile
               \vspace{0.15cm}
               \State  Compute  $\mu_{opt}$ and $\Sigma_{opt}$  using   \eqref{optimal_sampling_center}.
            % \vspace{0.05cm}
               \State Generate $N_{2}$ samples $\mathcal{S}=\{\mb{x}_{1}, \cdots, \mb{x}_{N_{2}}\}$ from $\mathcal{N}_{T}(\mu_{opt}, \Sigma_{opt})$.
               \State Calculate $\hat{P}_{\mathcal{F}}^{SAIS}$ using Eq.\eqref{failure_probability_estimation}.\\
               \Return  $\hat{P}_{\mathcal{F}}^{SAIS}, \mathcal{D}_{adaptive} = \{\mb{x}_{i}| g_{k}(\mb{x}_{i})>0, \, \mb{x}_{i} \in \mathcal{S}\}$.
           \end{algorithmic}
       \end{algorithm}

We present in Algorithm \ref{Algorithm2} the summary of the self-adaptive importance sampling (SAIS) procedure. In our experiments, we observe that by selecting $p_0=0.1,$ SAIS can self terminate rapidly with high numerical accuracy.

\paragraph{Remark 1:}
The above strategy can be easily applied to time-dependent problems and problems in unbounded domains. For instance, when $\Omega$ is unbounded, the intermediate proposal distribution can be simply chosen as Gaussian, which has the same form with the prior distribution $\omega$ (see our numerical experiments in Section 5).
Moreover, in this work, we simply consider a truncated Gaussian distribution to approximate the ``posterior failure" density $h_{opt}.$  One can easily generalize the truncated Gaussian model into more complex models such as the Gaussian mixture model;  see Appendix \ref{mGau}. The optimal choice of the ``posterior failure" density $h_{opt}$ will be part of our future works.

\paragraph{Remark 2:}
It is important to note that our FI-PINNs framework supports a wide range of performance functions. Despite the fact that we only pay attention to the residual function, it is easy to generalize to other concepts, such as its gradient function or causality-based weighted residual\cite{wang2022respecting}; see also Appendix \ref{ACeq}.

\section{Convergence analysis}\label{analysis_of_convergence}
This section is devoted to the convergence analysis  of FI-PINNs. We shall need the following assumptions.
 \vspace{0.15cm}
\begin{assumption}\label{assumption1}
Let $\mathcal{A}$ be a linear operator that maps $\mathcal{X}\rightarrow \mathcal{X}$ in problem \eqref{nonlinear-pde}, where $\mathcal{X}\subset \Omega$ is a Hilbert space. We assume that the operator $\mathcal{A}$ and the boundary operator $\mathcal{B}$ satisfy
        \begin{equation}
            C_{1}\|v\|_{2, \Omega} \leq \|\mathcal{A}v\|_{2, \Omega} + \|\mathcal{B}v\|_{2, \partial\Omega} \leq C_{2}\|v\|_{2, \Omega},\quad \forall v\in \mathcal{X}
        \end{equation}
        where the positive constants $C_1$ and $C_2$ are independent of $v$.
\end{assumption}
 \vspace{0.15cm}

\begin{assumption}
        \label{Assumption2}
 Assume the neural network $u(\mb{x};\theta)$ can be sufficiently trained  to ensure that  the  residual function is bounded, which is $
        M:= \max_{\mb{x}\in \Omega} |r(\mb{x};\theta^{*})| < \infty,
$
and that for any $\epsilon_{b}$,
     \begin{equation}
            \big\|\mathcal{B}(u(\mb{x}) - u(\mb{x};\theta^{*}))\big\|_{2, \partial\Omega} \leq \epsilon_{b}.
        \end{equation}
     \end{assumption}

\begin{assumption}
    \label{Assumption3}
    Suppose the residual function $r(\mb{x};\theta^{*})$ induced by the network $u(\mb{x};\theta^{*})$ is bounded, that is,
    \begin{equation}
        M:= \max_{\mb{x}\in \Omega} |r(\mb{x};\theta^{*})| < \infty
    \end{equation}
\end{assumption}

Assumption 4.2 simply supposes that we have successfully approximated the boundary data. We are now ready to present the following theorem.

\begin{theorem} \label{theorem1}
Assume the problem domain $\Omega$ is bounded, and let $u(\mb{x};\theta^{*})$ be an FI-PINNs solution of \eqref{nonlinear-pde}.  Suppose Assumptions \eqref{assumption1}, \eqref{Assumption2} and \eqref{Assumption3} are satisfied. Then the following error estimate holds
\begin{equation}
\big\|u(\mb{x}) - u(\mb{x};\theta^{*})\big\|_{2,\Omega} \leq \small{\sqrt{2}C_{1}^{-1}\big(S_{\Omega}(M^{2}\epsilon_{p} + \epsilon_{r}^{2}) + \epsilon_{b}^{2}\big)^{\frac{1}{2}}},
\end{equation}
where $$M = \max_{\mb{x}\in \Omega}|r(\mb{x};\theta^{*})|,$$ $S_{\Omega}$ is the area of $\Omega$ and $\epsilon_{r},\epsilon_{p}$ are the pre-given tolerances.
        \end{theorem}
  \vspace{0.05cm}
\begin{proof}
Let $v  = u(\mb{x}) - u(\mb{x};\theta^{*})$.
Using Assumption \ref{assumption1}, it can be shown that
 \begin{equation*}
            \begin{split}
             &\|u(\mb{x}) - u(\mb{x};\theta^{*})\|_{2, \Omega} \\[3mm]
             &\leq C_{1}^{-1}\left(\|\mathcal{A}(u(\mb{x}) - u(\mb{x};\theta^{*}))\|_{2,\Omega} + \|\mathcal{B}(u(\mb{x}) - u(\mb{x};\theta^{*}))\|_{2,\partial\Omega}\right)\\[3mm]
            & = C_{1}^{-1}\left(\|r(\mb{x};\theta^{*})\|_{2,\Omega} + \|\mathcal{B}(u(\mb{x}) - u(\mb{x};\theta^{*}))\|_{2,\partial\Omega}\right)\\
             & \leq \sqrt{2}C_{1}^{-1}\left(\|r(\mb{x};\theta^{*})\|_{2,\Omega}^{2} + \|\mathcal{B}(u(\mb{x}) - u(\mb{x};\theta^{*}))\|_{2,\partial\Omega}^{2}\right)^{\frac{1}{2}}.
            \end{split}
 \end{equation*}

From Assumption \ref{Assumption2}, we have
         \begin{equation}
            \label{estimate}
            \begin{split}
            \|u(\mb{x}) - u(\mb{x};\theta^{*})\|_{2, \Omega}
            \leq \sqrt{2}C_{1}^{-1}\left(\|r(\mb{x};\theta^{*})\|_{2,\Omega}^{2} +\epsilon_{b}^{2}\right)^{\frac{1}{2}}.
            \end{split}
         \end{equation}

For the bounded domain $\Omega$, we  have
         \begin{equation*}
            \begin{split}
            \|r(\mb{x};\theta^{*})\|_{2,\Omega}^2 &= \int_{\Omega}r(\mb{x};\theta^{*})^{2}d\mb{x}\\
            & = \int_{\Omega_{\mathcal{F}}}r(\mb{x};\theta^{*})^{2}d\mb{x} + \int_{\Omega_{\mathcal{S}}}r(\mb{x};\theta^{*})^{2}d\mb{x}.
            \end{split}
         \end{equation*}

Using the definition of the failure probability $P_{\mathcal{F}}$ on bounded domain, it can be shown that
         \begin{equation}
P_{\mathcal{F}} = \int_{\Omega}\mathbb{I}_{\Omega_{\mathcal{F}}}(\mb{x})d\mb{x}=  \int_{\Omega_{\mathcal{F}}}d\mb{x} = \frac{S_{\Omega_{\mathcal{F}}}}{S_{\Omega}} ,
         \end{equation}
where $S_{\Omega},S_{\Omega_{\mathcal{F}}}$ denote the areas of $\Omega,\Omega_{\mathcal{F}}$ respectively. If $P_{\mathcal{F}} < \epsilon_{p}$, we have
         \begin{eqnarray}
             S_{\Omega_{\mathcal{F}}} < S_{\Omega}\epsilon_{p}.
         \end{eqnarray}

Thus,
\begin{equation}
    \label{the_first_part}
    \int_{\Omega_{\mathcal{F}}}r(\mb{x};\theta^{*})^{2}d\mb{x} \leq \max_{\mb{x}\in \Omega_{\mathcal{F}}} |r(\mb{x};\theta^{*})|^{2} S_{\Omega_{\mathcal{F}}} \leq M^{2}S_{\Omega}\epsilon_{p},
\end{equation}
where $M = \max_{\mb{x}\in \Omega}|r(\mb{x};\theta^{*})|$.

Note that in the safe region $\Omega_{\mathcal{S}}$, the PDE residual $|r(\mb{x};\theta)|$ satisfy  $|r(\mb{x};\theta^{*})|< \epsilon_{r}$, which implies
\begin{equation}
    \label{the_second_part}
    \int_{\Omega_{\mathcal{S}}}r(\mb{x};\theta^{*})^{2}d\mb{x} \leq \epsilon_{r}^{2} S_{\Omega_{\mathcal{S}}} \leq S_{\Omega}\epsilon_{r}^{2}.
\end{equation}
Combing Eq.\eqref{the_first_part} and Eq.\eqref{the_second_part}, we can obtain
\begin{equation}
    \label{first_part_bound}
    \|r(\mb{x};\theta^{*})\|_{2,\Omega}^{2} \leq S_{\Omega}(M^{2}\epsilon_{p} + \epsilon_{r}^{2}).
\end{equation}
The desired result follows by substituting Eq.\eqref{first_part_bound} into Eq.\eqref{estimate}. \qquad \qquad \qquad \qquad
     \end{proof}

\section{Numerical experiments} \label{Numerical_experiments}

In this section, we present several numerical experiments to verify the effectiveness of FI-FINNs.

\subsection{Experiment Setup} \label{experiment_setup}

To better present the numerical results, we shall perform the following three types of approaches:
\begin{itemize}
 \item The conventional PINNs \cite{raissi2019physics}, or baseline PINNs approach.  This method update the training dataset by using the uniform sampling strategy.
\item The residual-based adaptive refinement method (RAR) \cite{lu2021deepxde}. This method update the training set by selecting collocation points with the first $(m)$ largest  residual values, in a large set of candidate points.
\item The FI-FINNs with SAIS presented in Section \eqref{self adaptive_importance_sampling}.
\end{itemize}

   \begin{figure}[t]
        \centering
\includegraphics[width = 0.5\linewidth]{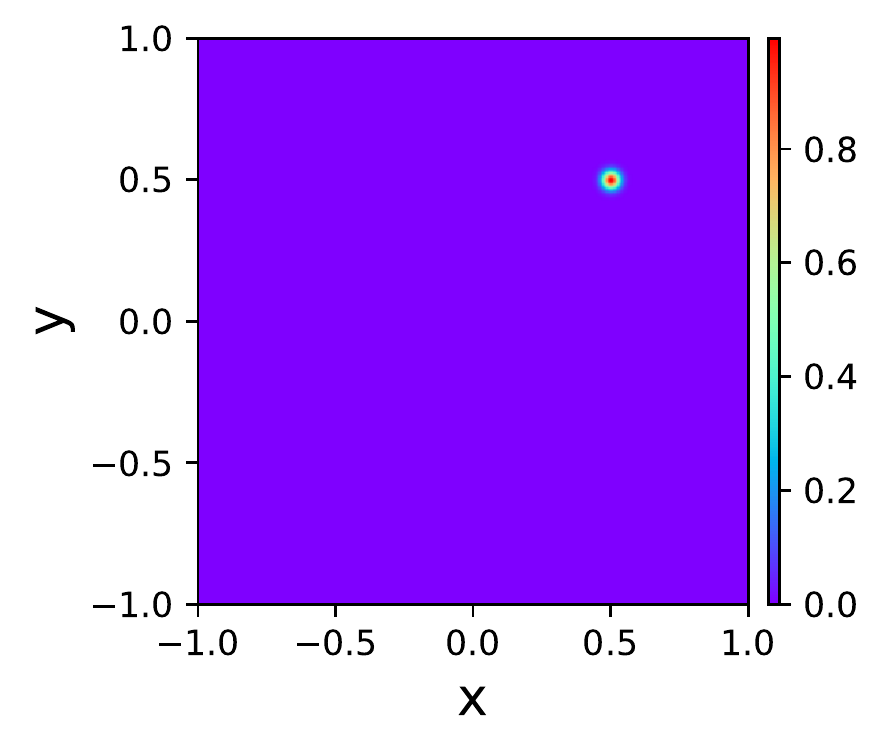}
    \caption{Exact solution for the two-dimensional peak test problem.}
    \label{one_peak_exact_solution}
\end{figure}

      \begin{figure}[t]
        \begin{center}
        \begin{overpic}[width=0.45\textwidth, clip=true,tics=10]{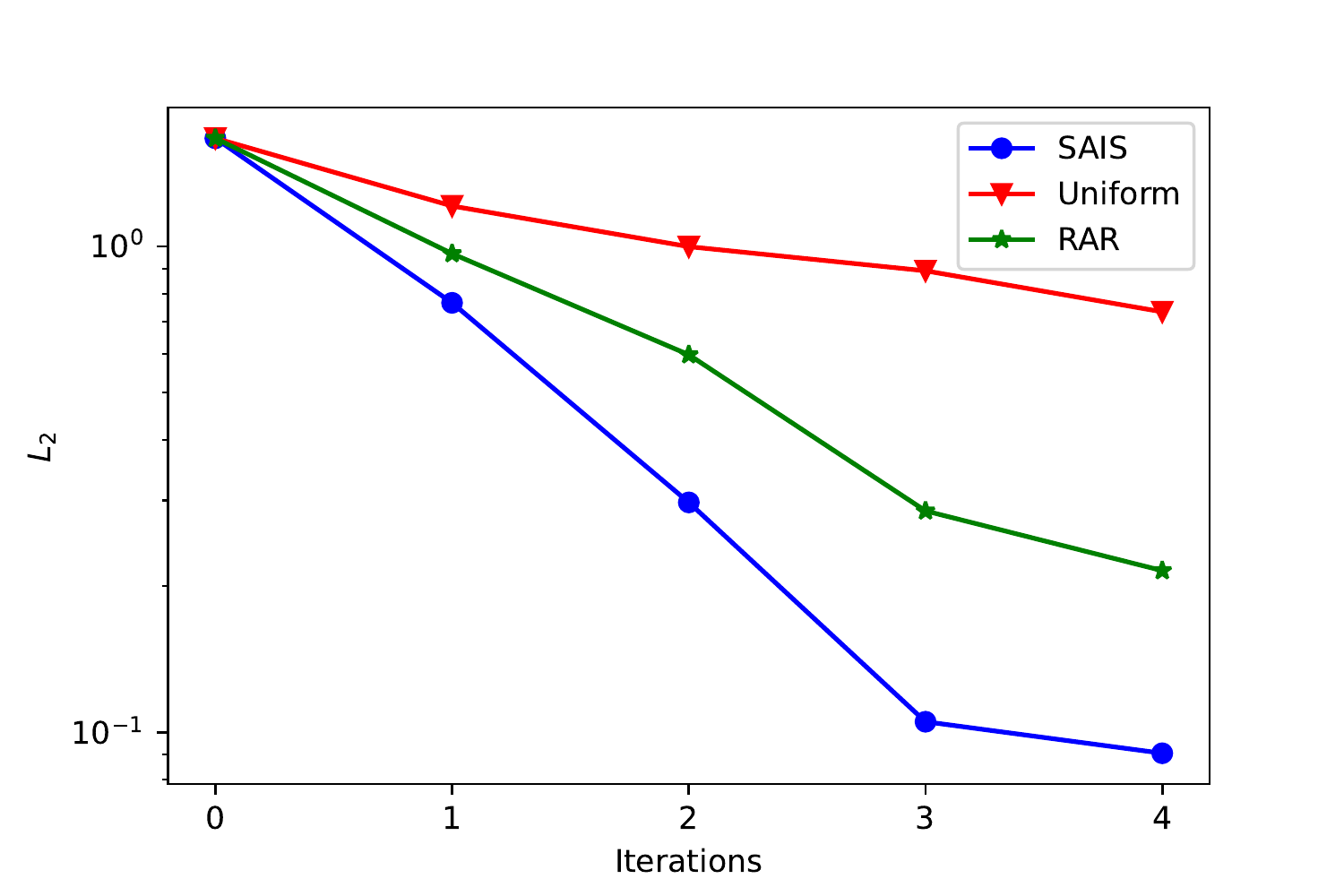}
        % \put (30,30) {\scriptsize {\bf Samples distribution by  ARA}}
        \end{overpic}
        \begin{overpic}[width=0.45\textwidth, clip=true,tics=10]{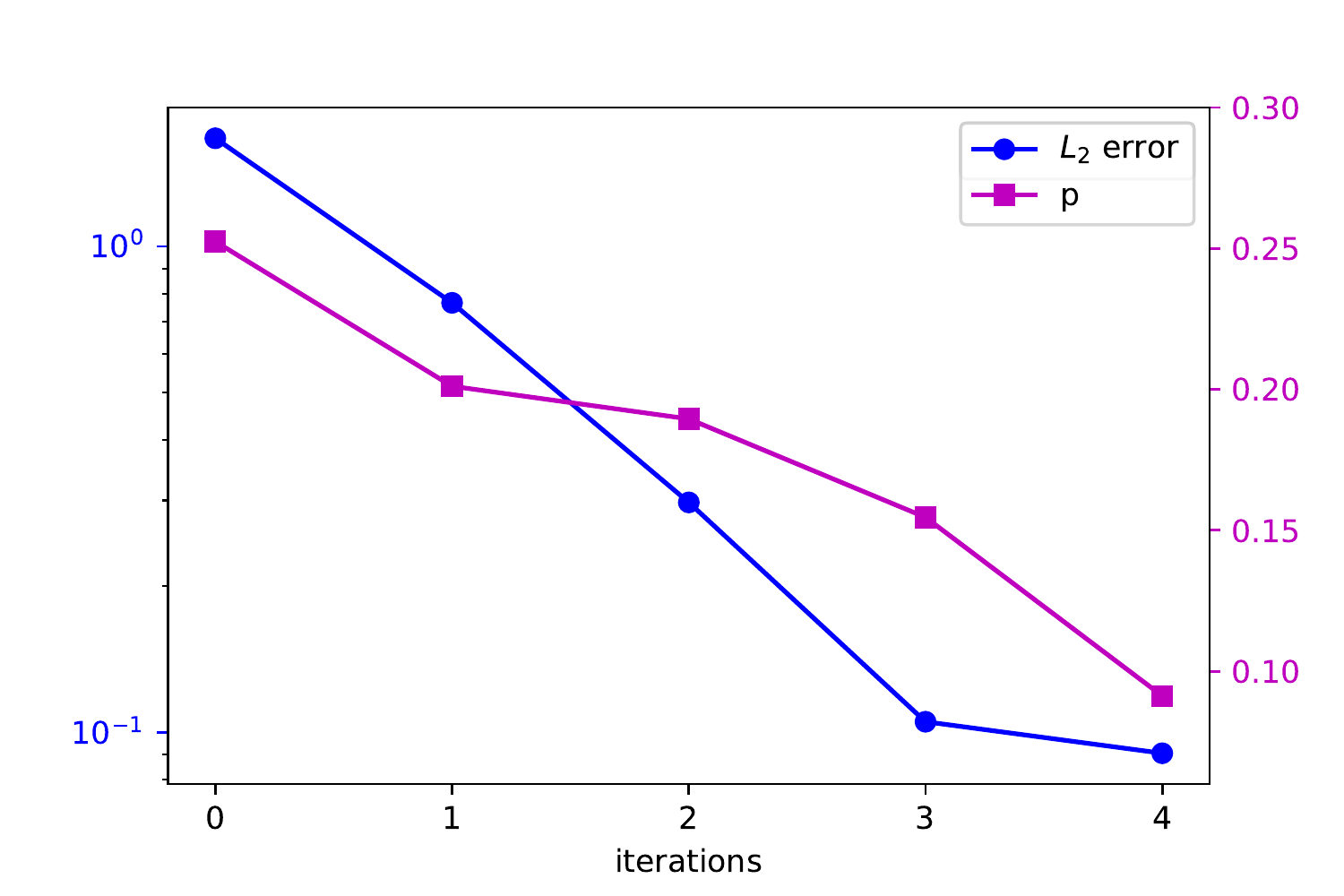}
        % \put (30,30) {\scriptsize {\bf Samples distribution by  ARA}}
        \end{overpic}
        \end{center}
        \caption{Relative $L_{2}$ error (left) and the  estimated failure probability (right) over updates.}
        \label{one_peak_error}
        \end{figure}

In all our numerical results, we will refer to ``\textit{Uniform}" as the traditional PINNs, ``\textit{RAR}" as the residual-based adaptive refinement approach, and ``\textit{SAIS}" as our FI-PINNs approach. To make a fair comparison, we shall first train the network completely using a small uniformly distributed dataset, and then we will train the three distinct models that are descended from the original model using datasets with the same size but produced using three different strategies.

Unless otherwise specified, we shall use the following parameters:  $N_{1} = 300, p_{0} = 0.1, N_2 = 1000$ in \textit{SAIS}.  The network used to implement the experiments is  a fully connected neural network with 7 hidden layers with 20 neurons in each layer, and the activation function is chosen as the \texttt{tanh} function. In the beginning, we always use 2000 collocation points and 200 boundary points to train the network. The \texttt{Adam} optimizer is adopted to optimize the loss function, with a learning rate 0.0001 and 10000 iterations.  The maximum adaptive iteration is set to $M=10.$ In order to test the validity of the method, we use the following  relative $L_{2}$ error:
   \begin{equation}
       \label{L_2 error}
       err_{L_{2}} = \frac{\sqrt{\sum_{i=1}^{N}\left|\hat{u}(\mb{x}_{i}) - u(\mb{x}_{i})\right|^{2}}}{\sqrt{\sum_{i=1}^{N}\left|u(\mb{x}_{i})\right|^{2}}},
   \end{equation}
   where $N$ represents the total number of test points chosen randomly in the domain, and $\hat{u}(\mb{x}_{i})$ and  $u(\mb{x}_{i})$ represent the predicted and the exact solution values, respectively.

 \subsection{Two-dimensional Poisson equation}
 \label{sec:52}

  \begin{figure}[t]
    \begin{center}
    \begin{overpic}[width=0.7\textwidth,trim= 35 0 45 15, clip=true,tics=10]{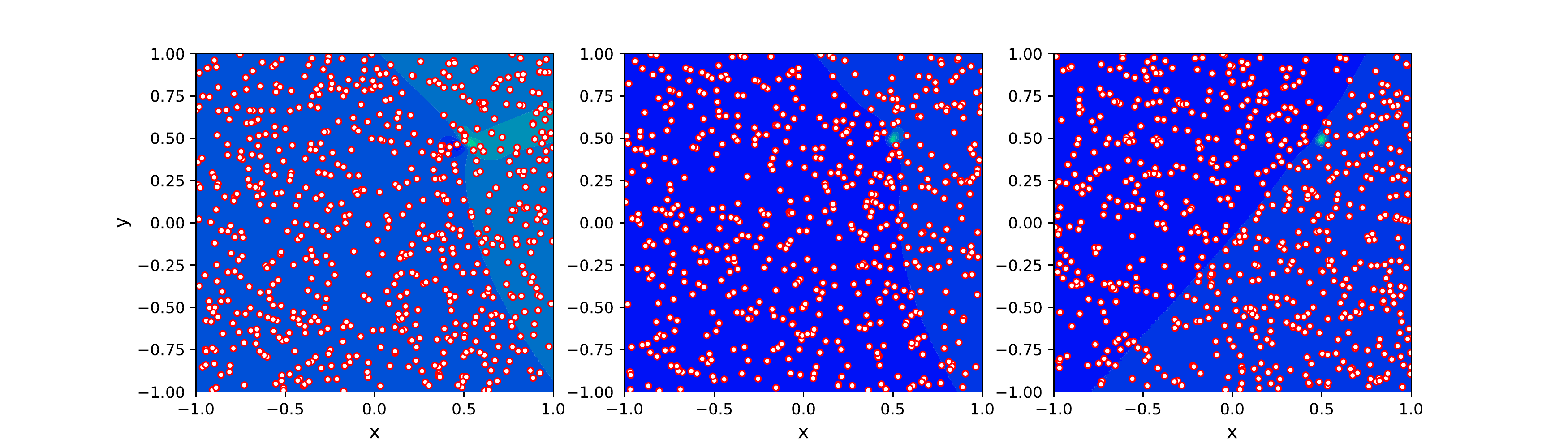}
    \put (38,30) {\scriptsize {\bf  Samples  by  \textit{Uniform}}}
    \end{overpic}
    \end{center}
    \vspace{0.3cm}
    \begin{center}
        \begin{overpic}[width=0.7\textwidth,trim= 35 0 45 15, clip=true,tics=10]{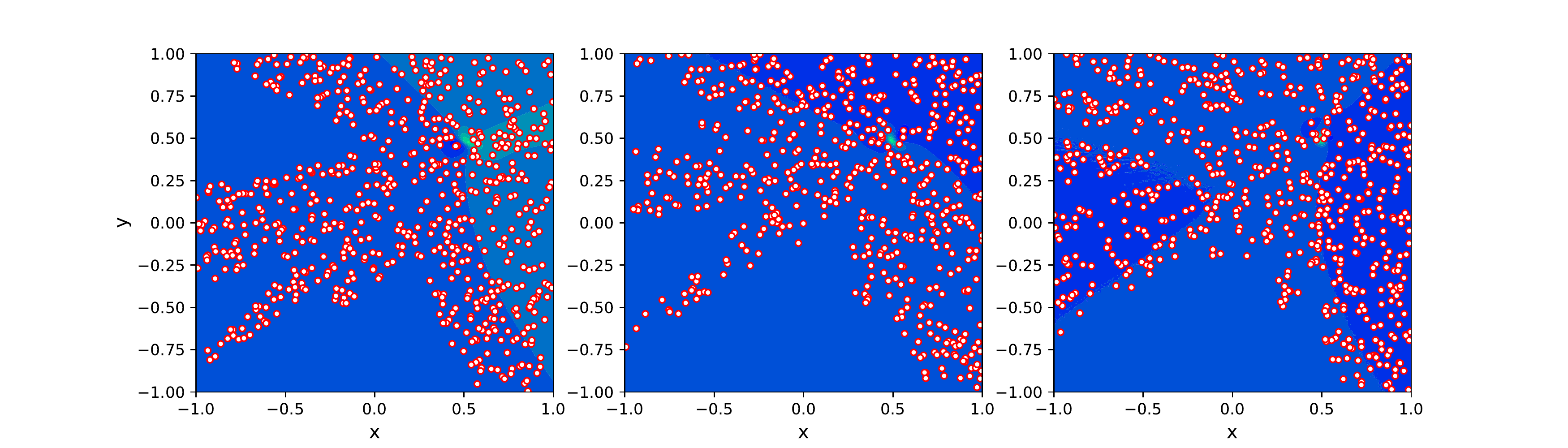}
             \put (38,30) {\scriptsize {\bf Samples  by  \textit{RAR}}}
      \end{overpic}
    \end{center}
    \vspace{0.3cm}
    \begin{center}
      \begin{overpic}[width=0.7\textwidth,trim= 35 0 45 15, clip=true,tics=10]{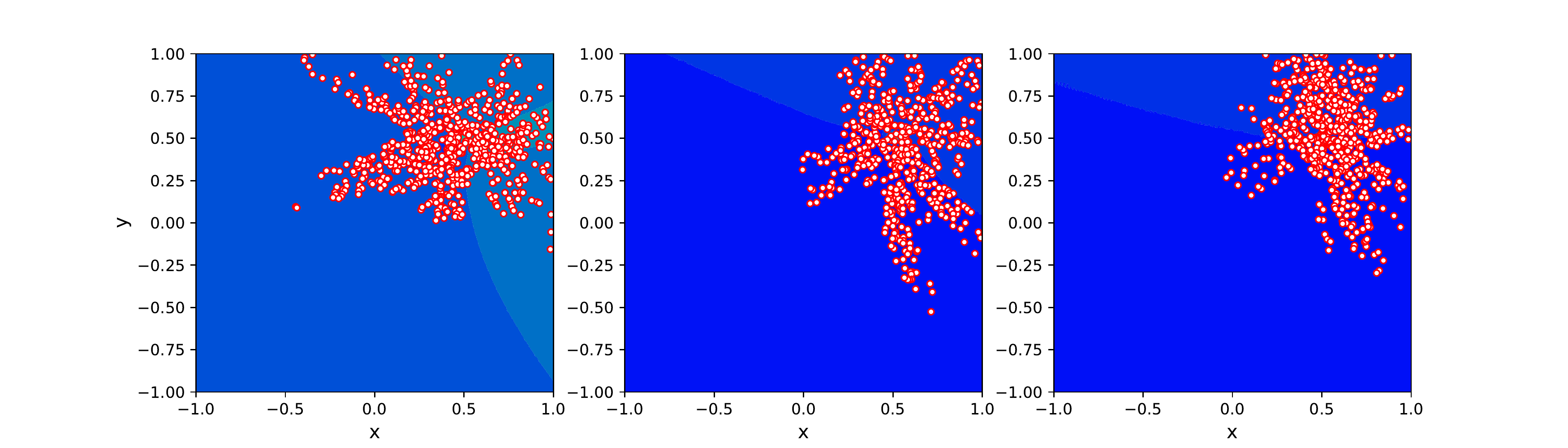}
           \put (38,30) {\scriptsize {\bf Samples  by \textit{SAIS}}}
      \end{overpic}
      \end{center}
      \vspace{-0.2cm}
    \caption{The distribution of collocation points obtained by  three sampling methods for the first three updates.}
    \label{one_peak_samples}
      \end{figure}

We first consider the following two-dimensional Poisson equation \cite{tang2021deep}:
   \begin{equation}
       \label{one_peak_equation}
       \begin{split}
       -\Delta u(x, y) &= f(x, y) \quad \mbox{in} \, \Omega \\
       u(x, y) &= g(x, y) \quad \mbox{on} \, \partial \Omega,
       \end{split}
   \end{equation}
   where  $\Omega$ is $[-1,1]^2$ and we specify the true solution as
   \begin{equation}
    \label{one_peak_true_solution}
    u(x, y) = \exp\left(-1000\left[(x - 0.5)^{2} + (y - 0.5)^{2}\right]\right),
       \end{equation}
which has a peak at $(0.5, 0.5)$ and decreases rapidly away from $(0.5, 0.5)$ (see Fig.\ref{one_peak_exact_solution}).

We first compare  the numerical results obtained by three different sampling strategies.  To construct our FI-PINNs framework, we set  $\epsilon_{p} = 0.1$ and $\epsilon_{r}=0.1$ in this example.   The relative $L_{2}$ errors are shown in the left plot of Fig. \ref{one_peak_error}.  The error achieved by \textit{SAIS} method decreases much faster than the \textit{Uniform} and \textit{RAR} strategies, and it is clearly seen that the mean $L_{2}$ error obtained by \textit{SAIS} is much smaller, which is $9.04\times 10^{-2}$, compared to $7.36\times 10^{-1}$ and $2.15\times 10^{-1}$ obtained by \textit{Uniform} and \textit{RAR} sampling strategies respectively. In the right plot of Fig. \ref{one_peak_error}, we show that the estimated failure probability has the same trend with the relative $L_{2}$ error. After four iterations, the training can be stopped when the estimated failure probability is smaller than the tolerance, meaning that our stop criteria is  effectiveness.

\begin{figure}[t]
    \begin{center}
             \begin{overpic}[width=0.9\textwidth,trim= 35 0 45 15, clip=true,tics=10]{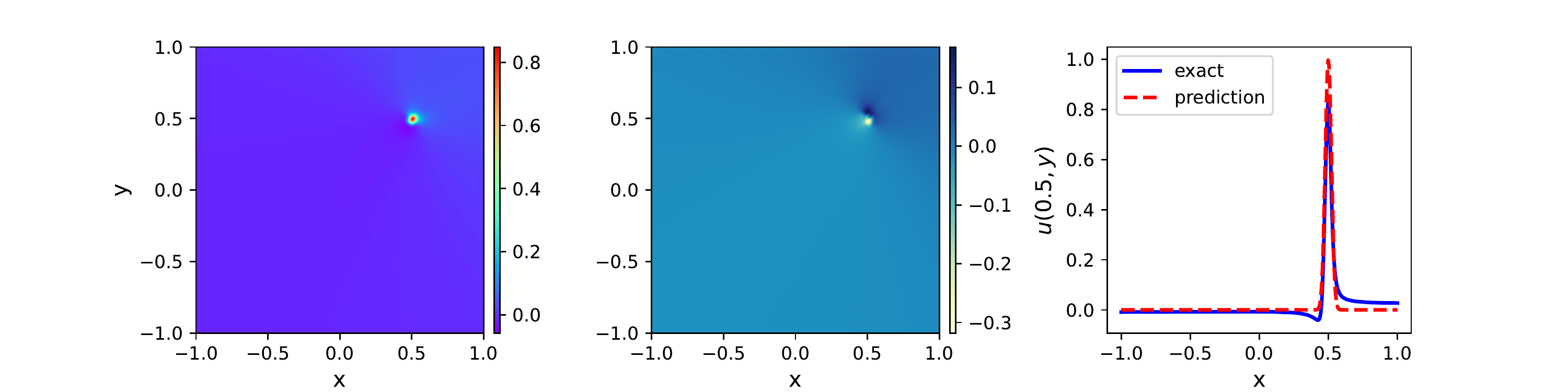}
             \put (30,26) {\scriptsize {\bf Numerical results  obtained by  \textit{Uniform}}}
            \end{overpic}
        \end{center}
        \vspace{0.3cm}
      \begin{center}
        \begin{overpic}[width=0.9\textwidth,trim= 35 0 45 15, clip=true,tics=10]{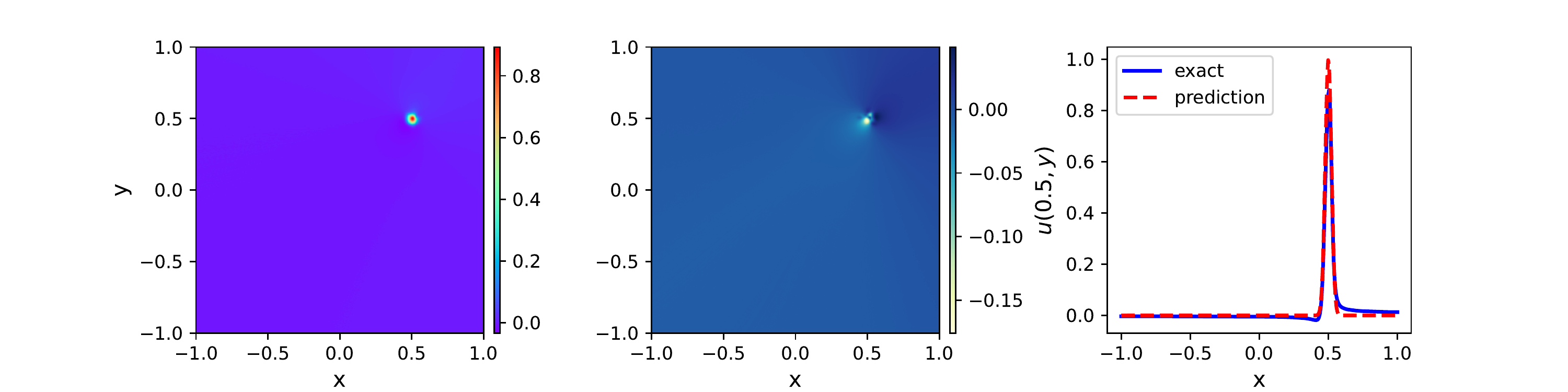}
    %      \caption{Solution errors obtained by ARA}
             \put (30,26) {\scriptsize {\bf Numerical results  obtained by \textit{RAR}}}
      \end{overpic}
    \end{center}
    \vspace{0.3cm}
        \begin{center}
      \begin{overpic}[width=0.9\textwidth,trim= 35 0 45 15, clip=true,tics=10]{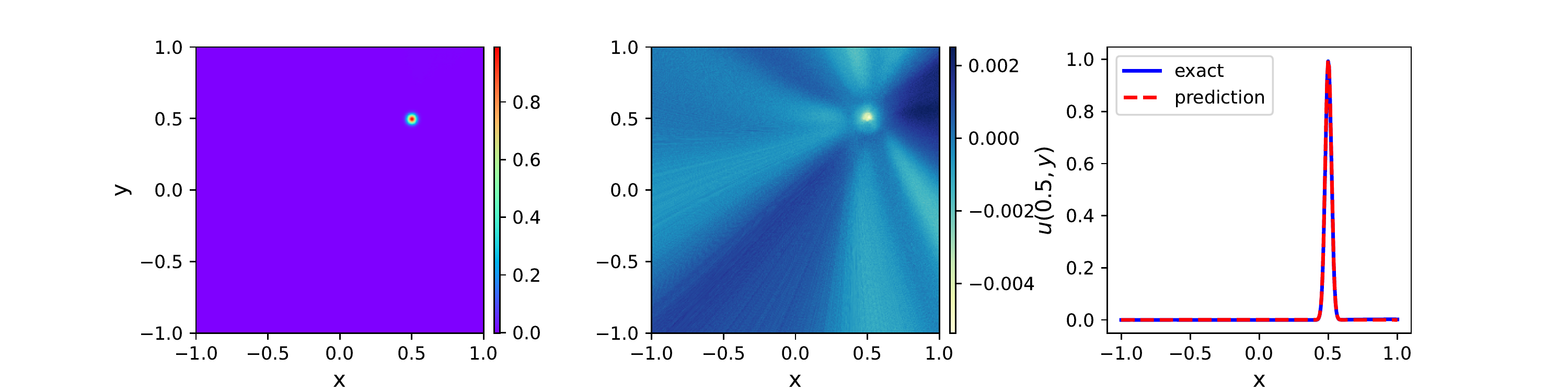}
           \put (30,26) {\scriptsize {\bf Numerical results  obtained  by \textit{SAIS}}}
      \end{overpic}
      \end{center}
           \vspace{-.2cm}
    \caption{The predicted solution(Left), absolute error (Middle) and the predicted curves at $x = 0.5$(Right)  obtained by the three different sampling strategies.}
    \label{one_peak_solution}
      \end{figure}

Fig.\ref{one_peak_samples} displays the distributions of the updated collocation points from the first three iterations using various sampling strategies. It is clear seen that the samples generated by \textit{SAIS} are much more concentrated around the peak $(0.5,0.5)$. The associated  predicted values,  absolute error and the predicted solution curve at $x = 0.5$ obtained by the three sampling strategies are shown in Fig.\ref{one_peak_solution}. It is evident that the absolute error obtained by \textit{SAIS} is much smaller, and there is hardy noticeable difference between the exact values and the predicted values at $x = 0.5$, indicating that the \textit{SAIS} sampling strategy outperforms the other two sampling strategies.

\begin{figure}
    \begin{center}
        \begin{overpic}[width = 0.45\textwidth, clip = true, tics = 10]{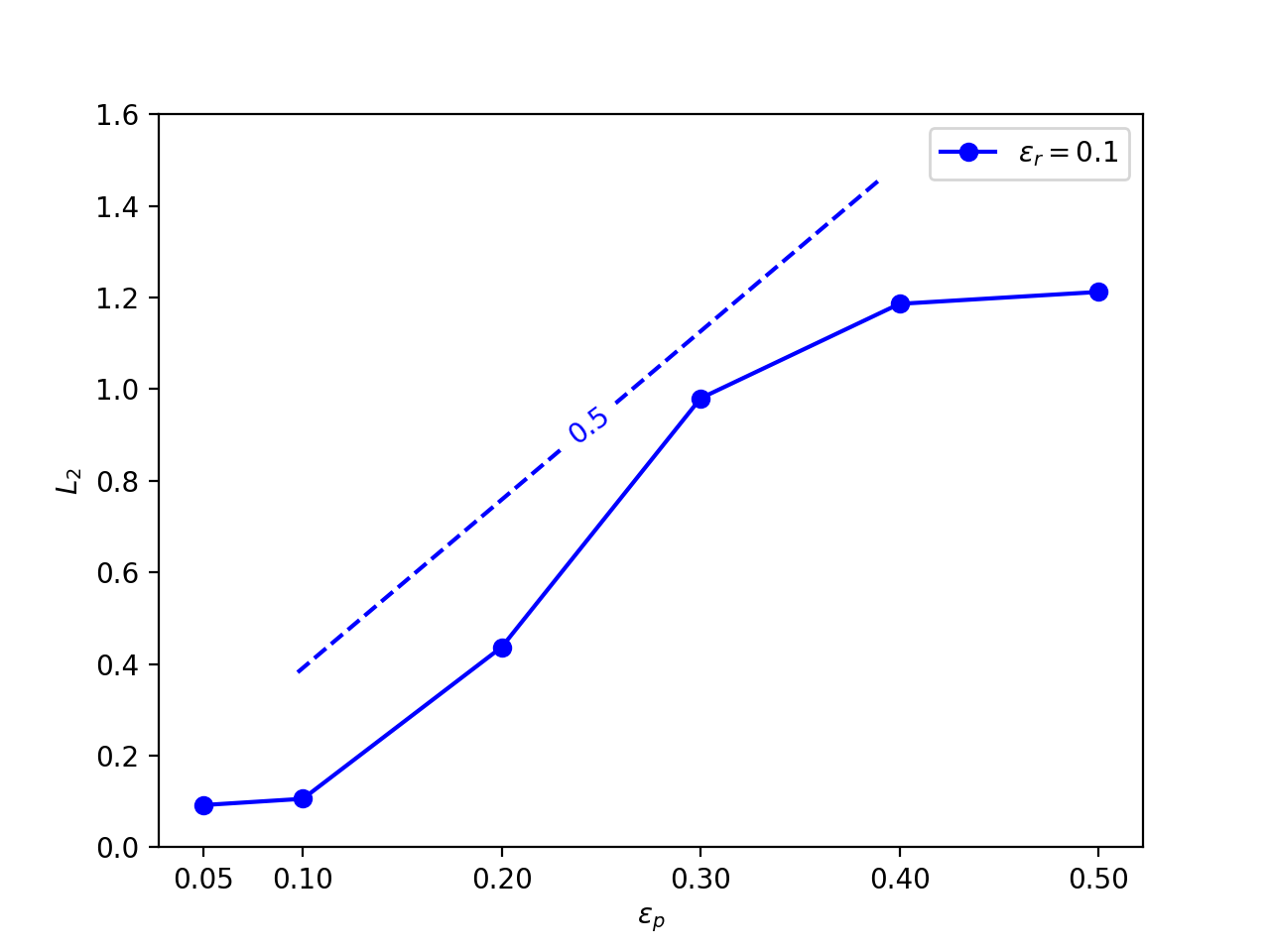}
        \end{overpic}
            \begin{overpic}[width = 0.45\textwidth, , clip = true, tics = 10]{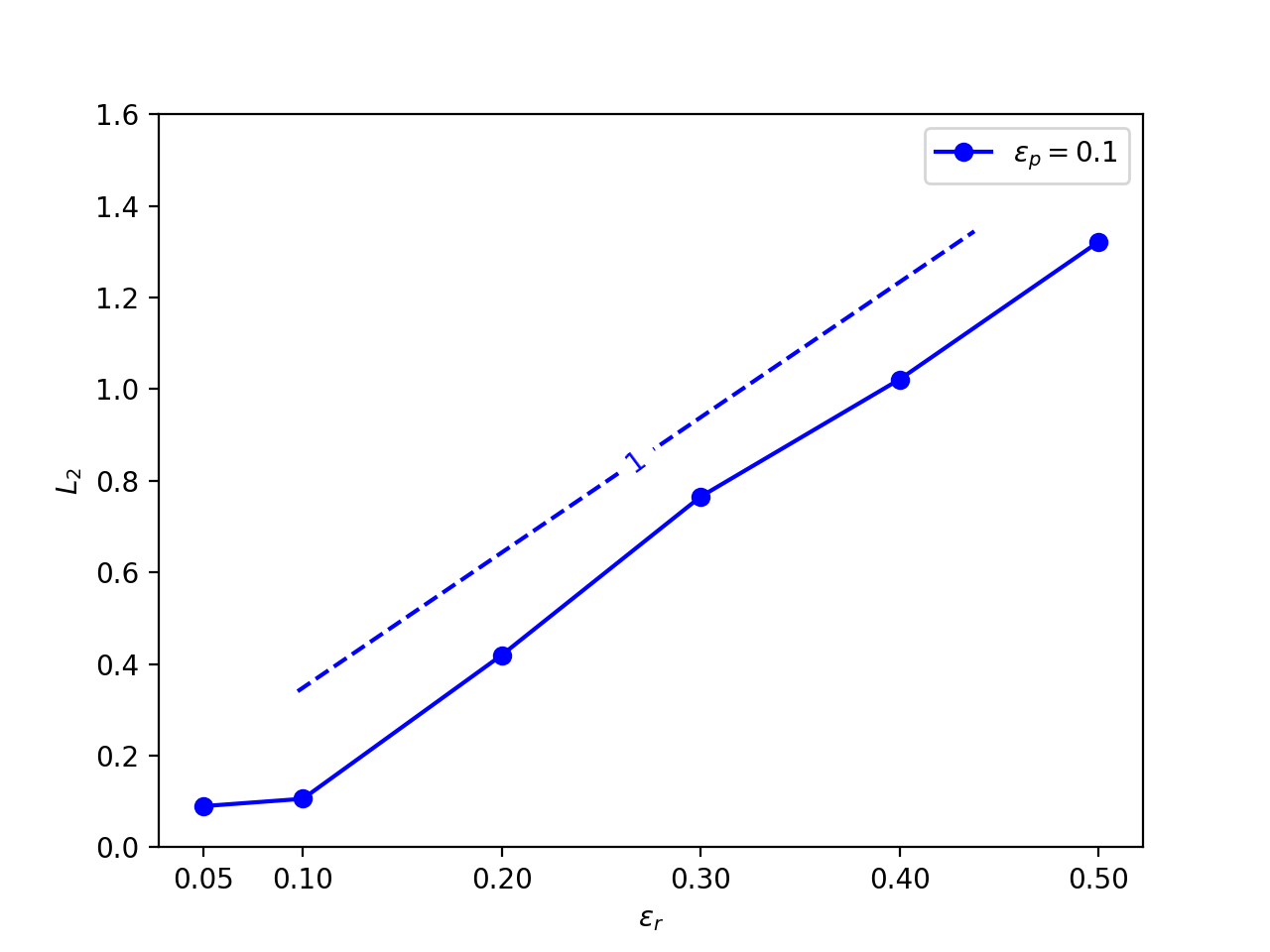}
            \end{overpic}
    \end{center}
    \caption{Error profiles when solving a 2D Poisson equation with one peak. Left: the prediction errors with respect to the $\epsilon_{p}$ with $\epsilon_{r} = 0.1$. Right: the prediction errors with respect to the  $\epsilon_{r}$ with $\epsilon_{p} = 0.1$.}
    \label{residual_failure_one_peak}
\end{figure}

We now examine the convergence properties of FI-PINNs. To show the convergence properties of FI-PINNs, we fix one tolerance parameter and train the network to obtain the prediction error when the other tolerance changes. The results are presented in Fig.\ref{residual_failure_one_peak}. It is clearly seen that the convergence rate with respect to  $\epsilon_{p}$ and $\epsilon_{r}$ are 1/2 and 1, respectively. This agrees well with Theorem \ref{theorem1}.

\subsection{Burgers' equation}
        We consider the following Burgers' equation:
        \begin{equation}
            \label{burgers_equation}
            \begin{split}
               &u_{t} + uu_{x} - \frac{0.01}{\pi}u_{xx} = 0,\\
               &u(x,0) = -\sin(\pi x),\\
               &u(t,-1) = u(t,1) = 0,
            \end{split}
        \end{equation}
where $(x, t)\in [-1,1]\times[0,1]$. Notice that there exists a stiff mutation at $x = 0$ in the solution, which is used here to test the viability of the FI-PINNs and the \textit{SAIS} method for time-dependent nonlinear problems.
\begin{figure}[t]
            \begin{center}
        \begin{overpic}[width=0.45\textwidth,trim= 0 0 0 0, clip=true,tics=10]{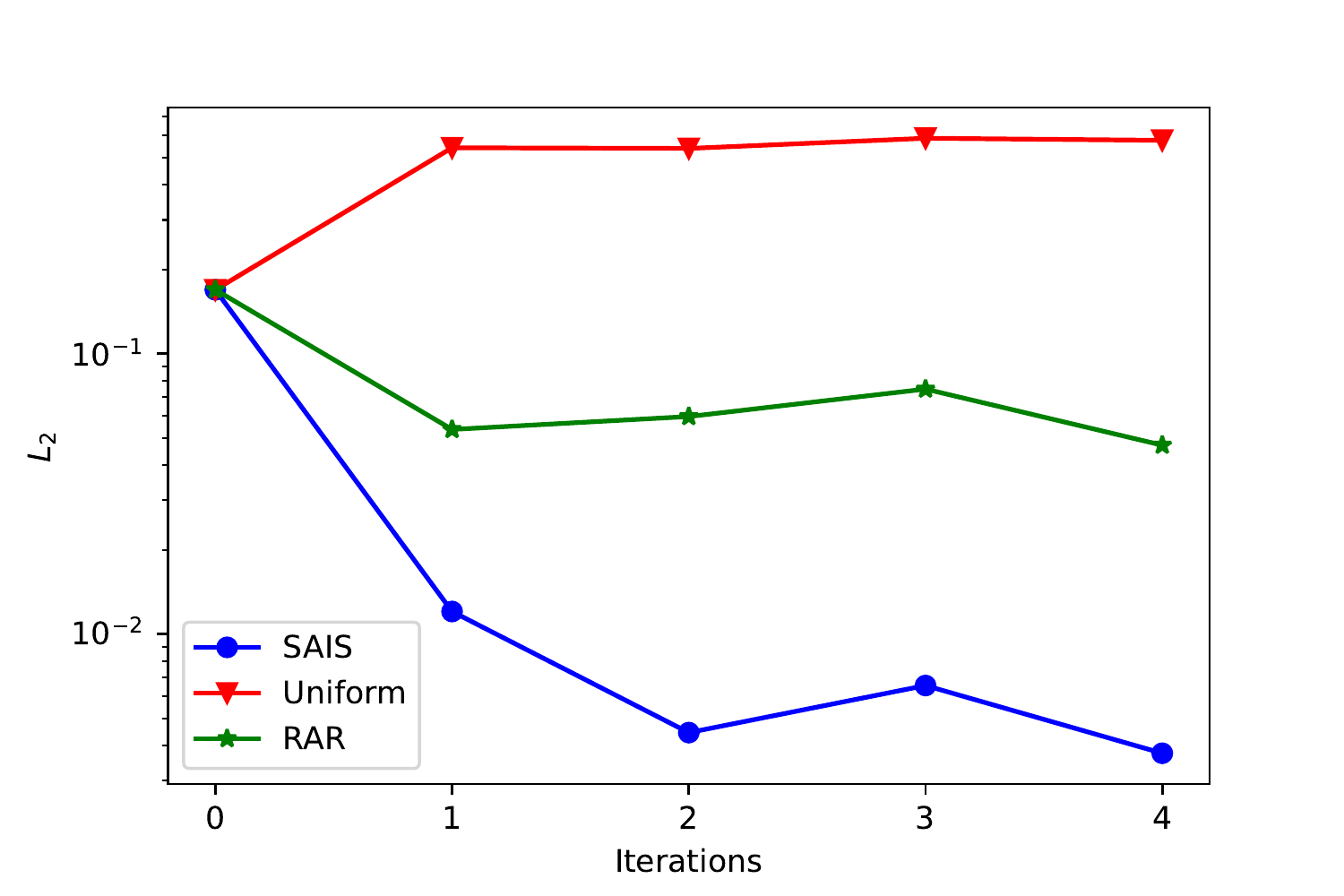}
        \end{overpic}
        \begin{overpic}[width=0.45\textwidth,trim= 0 0 0 0, clip=true,tics=10]{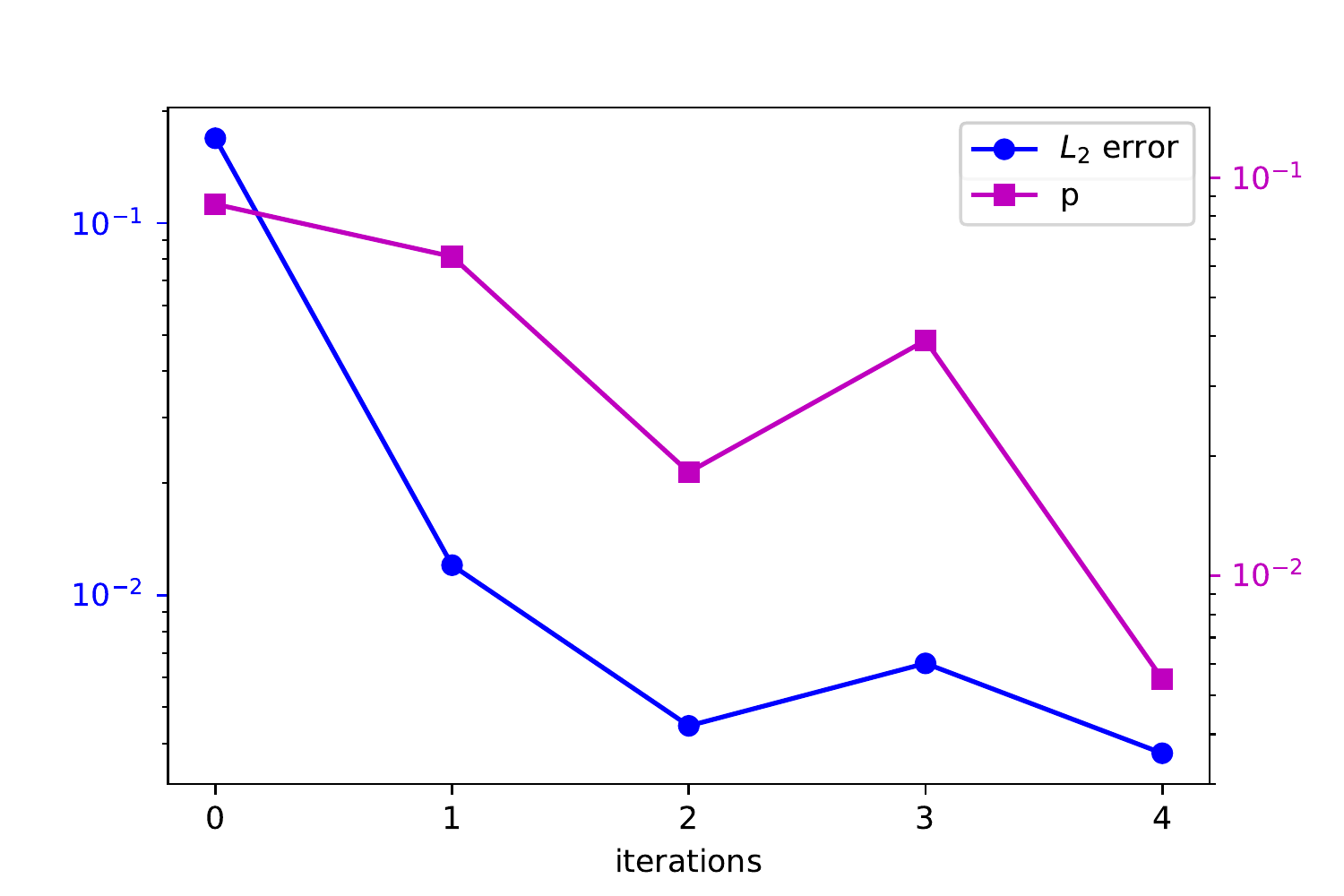}
        \end{overpic}
        \end{center}
        \caption{Relative $L_{2}$ error (left) and the  estimated failure probability (right) over updates.}
        \label{burgers_failure}
\end{figure}

\begin{figure}[t]
            \begin{center}

            \begin{overpic}[width=\textwidth,trim= 35 10 45 15, clip=true,tics=10]{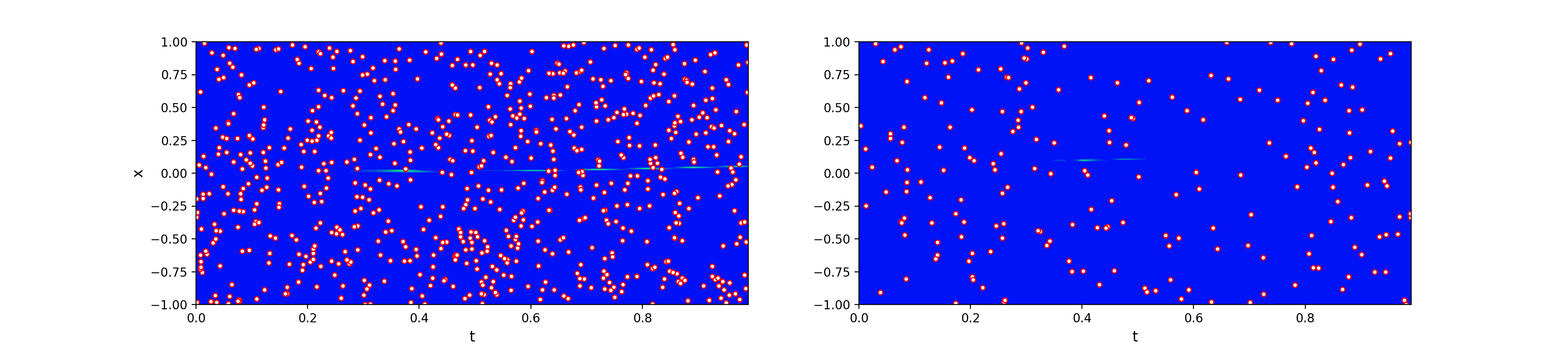}
            \put (42,22) {\scriptsize {\bf  Samples  by \textit{Uniform}}}
            \end{overpic}
            % \vspace{-.2cm}
            \end{center}
            \vspace{0.3cm}
            \begin{center}
            \begin{overpic}[width=\textwidth,trim= 35 10 45 15, clip=true,tics=10]{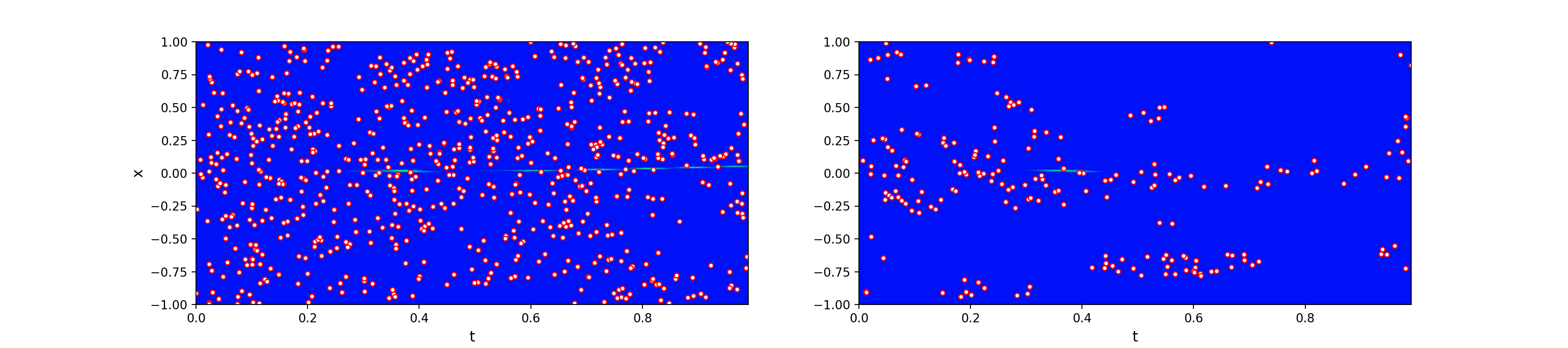}
            \put (42,22) {\scriptsize {\bf Samples  by \textit{RAR}}}
            \end{overpic}
            \end{center}
            \vspace{.3cm}
            \begin{center}
            \begin{overpic}[width=\textwidth,trim= 35 10 45 15, clip=true,tics=10]{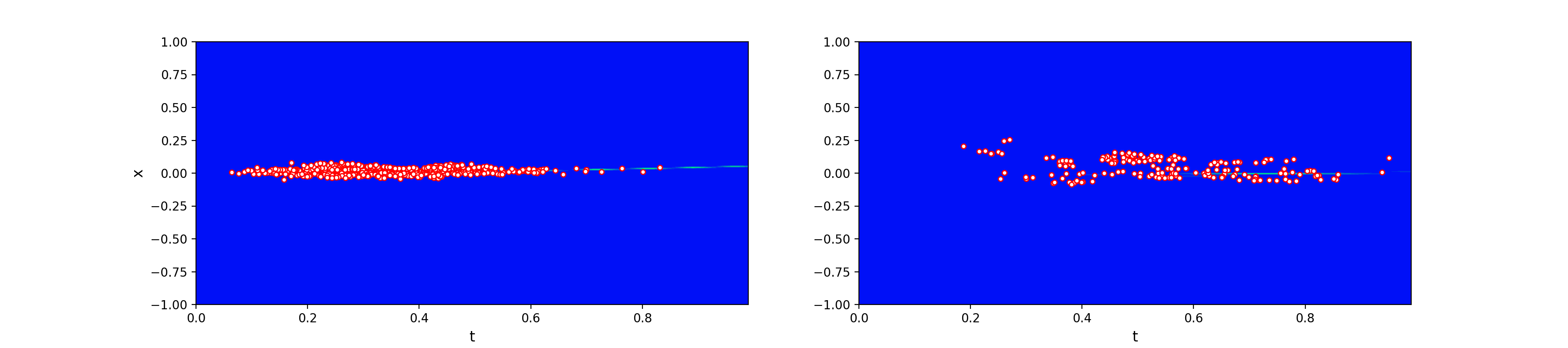}
            \put (42,22) {\scriptsize {\bf Samples  by \textit{SAIS}}}
            \end{overpic}
            \end{center}
            \vspace{-0.2cm}
            \caption{The distribution of collocation points obtained by  three sampling methods for the first two updates.}
            \label{burgers_samples}
\end{figure}

In this example,  the residual tolerance  $\epsilon_{r}$ and failure probability tolerance $\epsilon_{p}$ are both set to be 0.01.  To implement the experiment, we first train the network using  the \texttt{LBFGS} optimizer with a maximum of 50000 iterations. If the network needs to be updated (with new training points),  we first retrain it using the \texttt{Adam} optimizer  for a total of 10000 iterations, and then we utilize the \texttt{LBFGS} optimizer to fine tune the network with maximum 50000 iterations.

Figure \ref{burgers_failure} shows the prediction errors for three strategies, as well as the estimated failure probability with respect to iteration number. Again, it is clearly seen that the estimated failure probability and the prediction errors using \textit{SAIS}  have a similar pattern. The estimated failure probability is below the  tolerance $\epsilon_{p}$ after four updates.  Additionally,  the prediction errors of \textit{SAIS} are significantly lower ($5.59\times 10^{-3}$).

\begin{figure}[t]
        \begin{center}
                 \begin{overpic}[width=\textwidth,trim= 35 0 45 15, clip=true,tics=10]{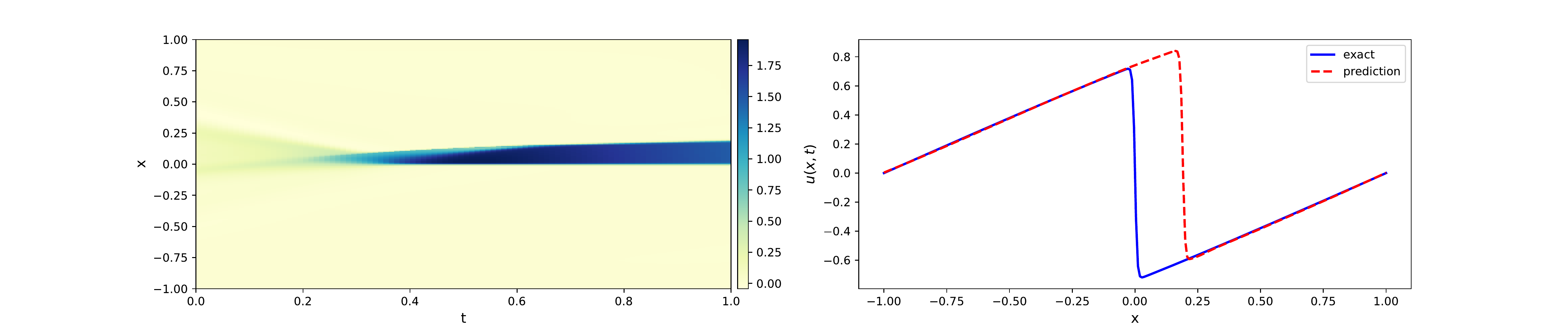}
                 \put (35,22) {\scriptsize {\bf  Numerical results obtained by  \textit{Uniform}}}
                \end{overpic}
            \end{center}
               \vspace{0.3cm}
           \begin{center}
            \begin{overpic}[width=\textwidth,trim= 35 0 45 15, clip=true,tics=10]{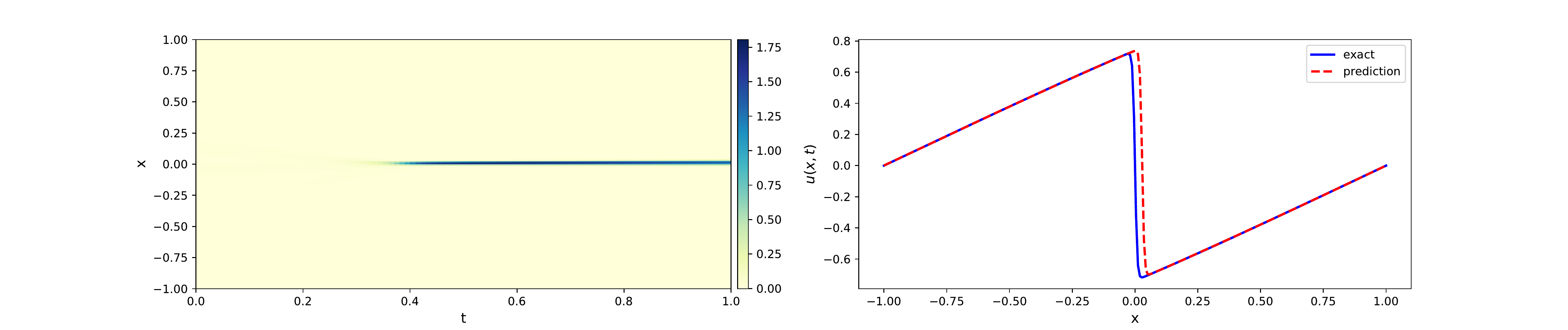}
                 \put (35,22) {\scriptsize {\bf Numerical results obtained by  \textit{RAR}}}
          \end{overpic}
        \end{center}
                  \vspace{.3cm}
        \begin{center}
          \begin{overpic}[width=\textwidth,trim= 35 0 45 15, clip=true,tics=10]{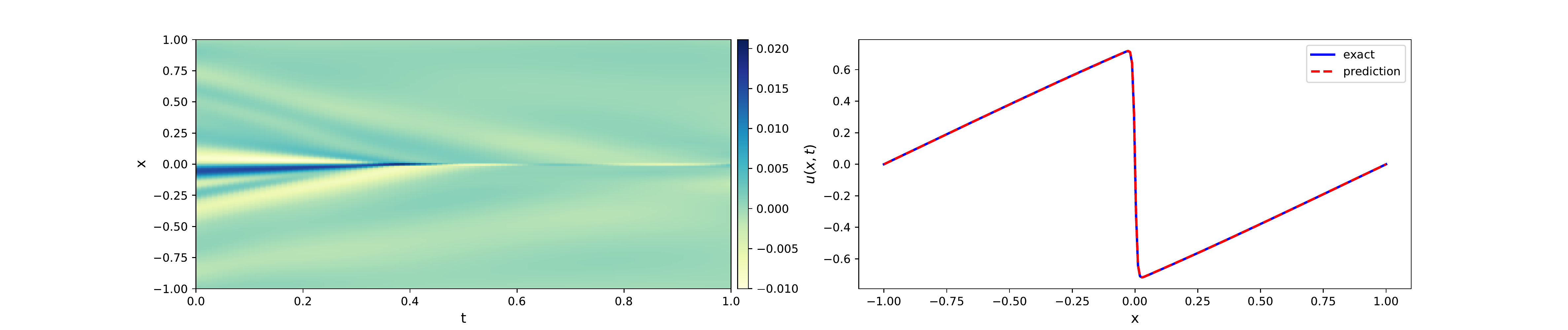}
               \put (35,22) {\scriptsize {\bf Numerical results obtained by \textit{SAIS}}}
          \end{overpic}
          \end{center}
        \caption{The predicted error (Left) and the predicted curves at $t = 1$(Right)  obtained by the three different sampling strategies.}
        \label{burgers_solution}
          \end{figure}

\begin{figure}
        \begin{center}
            \begin{overpic}[width = 0.45\textwidth, clip = true, tics = 10]{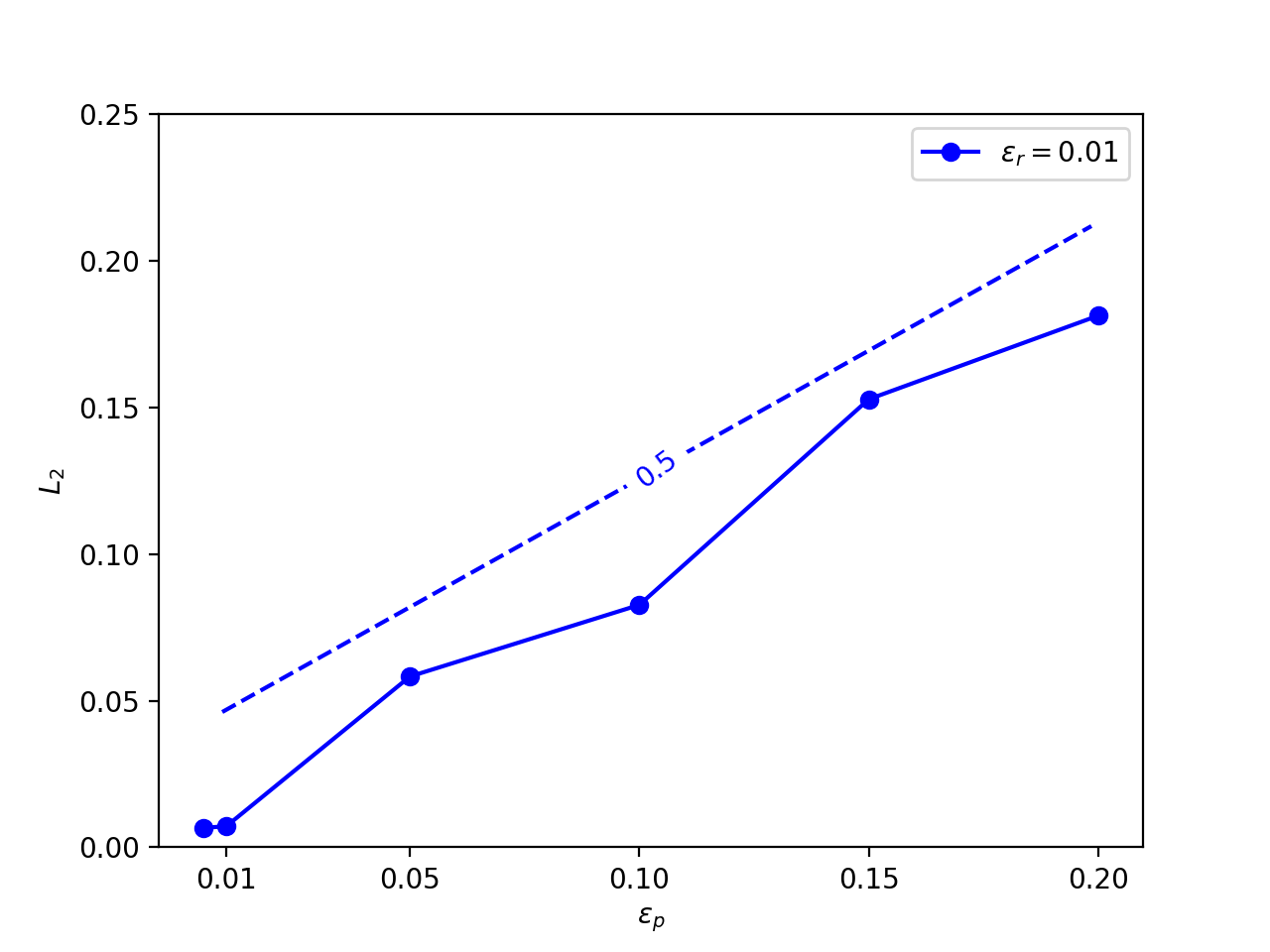}
            \end{overpic}
                \begin{overpic}[width = 0.45\textwidth, clip = true, tics = 10]{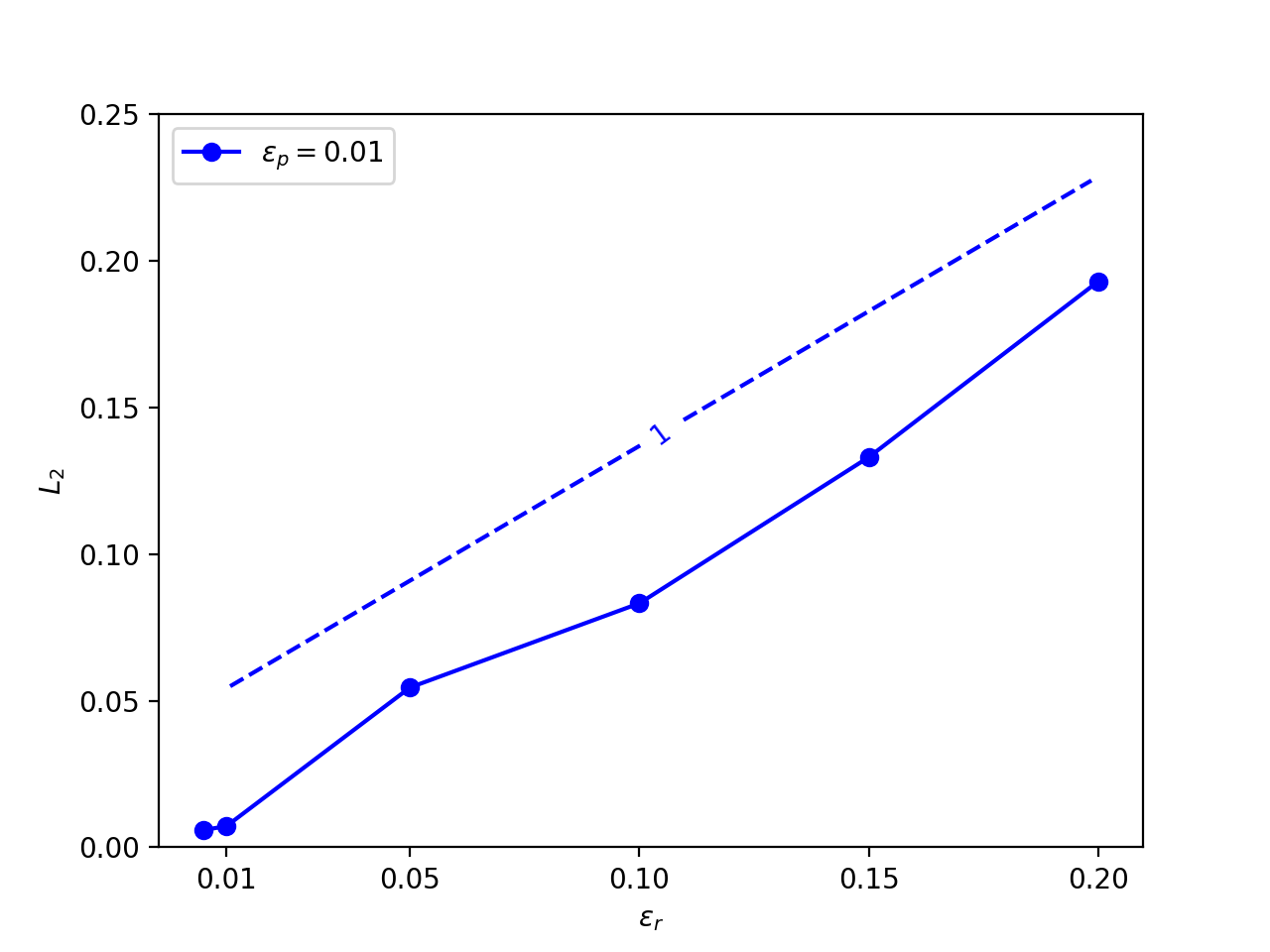}
                \end{overpic}
        \end{center}
        \caption{Errors profiles when solving  Burgers' equation. Left: the prediction errors with respect to the $\epsilon_{p}$ with $\epsilon_{r} = 0.01$. Right: the prediction errors with respect to the  $\epsilon_{r}$ with $\epsilon_{p} = 0.01$.}
        \label{burgers_residual_failure}
    \end{figure}

New training points generated in the first two iterations are presented in Fig. \ref{burgers_samples}. We observe that the samples produced by the \textit{SAIS} strategy are concentrated more in the regions where the residual error is high. In Fig.\ref{burgers_solution}, we show the final absolute prediction error as well as the predicted solution curve at the time $t=1$. It is obvious that the predicted values produced by the \textit{SAIS} match well with the exact solution. Although we have only established the theoretical results for the linear model. We indeed observe in Fig. \ref{burgers_residual_failure} the convergence rates also hold for this nonlinear example.

\subsection{The high-dimensional Poisson equation}

Consider the the following $d$-dimensional elliptic equation
\begin{equation*}
    \label{high dimensional}
    -\Delta u(\mb{x}) = f(\mb{x}),\quad \mb{x}\in [-1,1]^{d},
\end{equation*}
with an exact solution $$u(\mb{x}) = e^{-10\|\mb{x}\|_{2}^{2}}.$$
We set $d=9$, and use 20000 collocation points and 900 boundary points to train the initial network. The hyperparameter for \textit{SAIS}  is set at $N_{2} = 5000$.  Additionally, we define  the residual and failure probability tolerances as  $\epsilon_{r} = 0.01, \epsilon_{p} = 0.005$ respectively.

\begin{figure}[t]
    \begin{center}
    \begin{overpic}[width=0.45\textwidth, clip=true,tics=10]{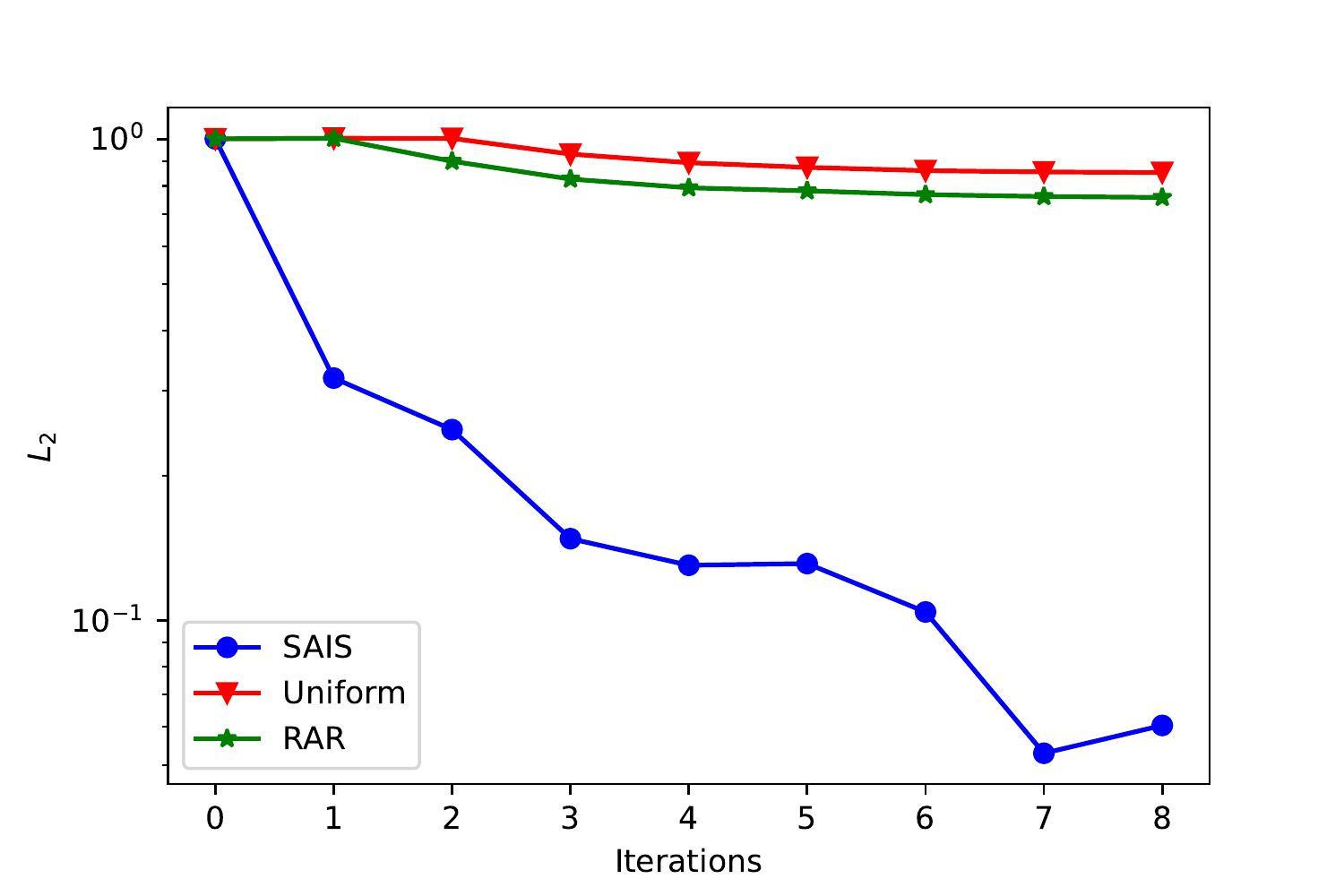}
    \end{overpic}
    \begin{overpic}[width=0.45\textwidth, clip=true,tics=10]{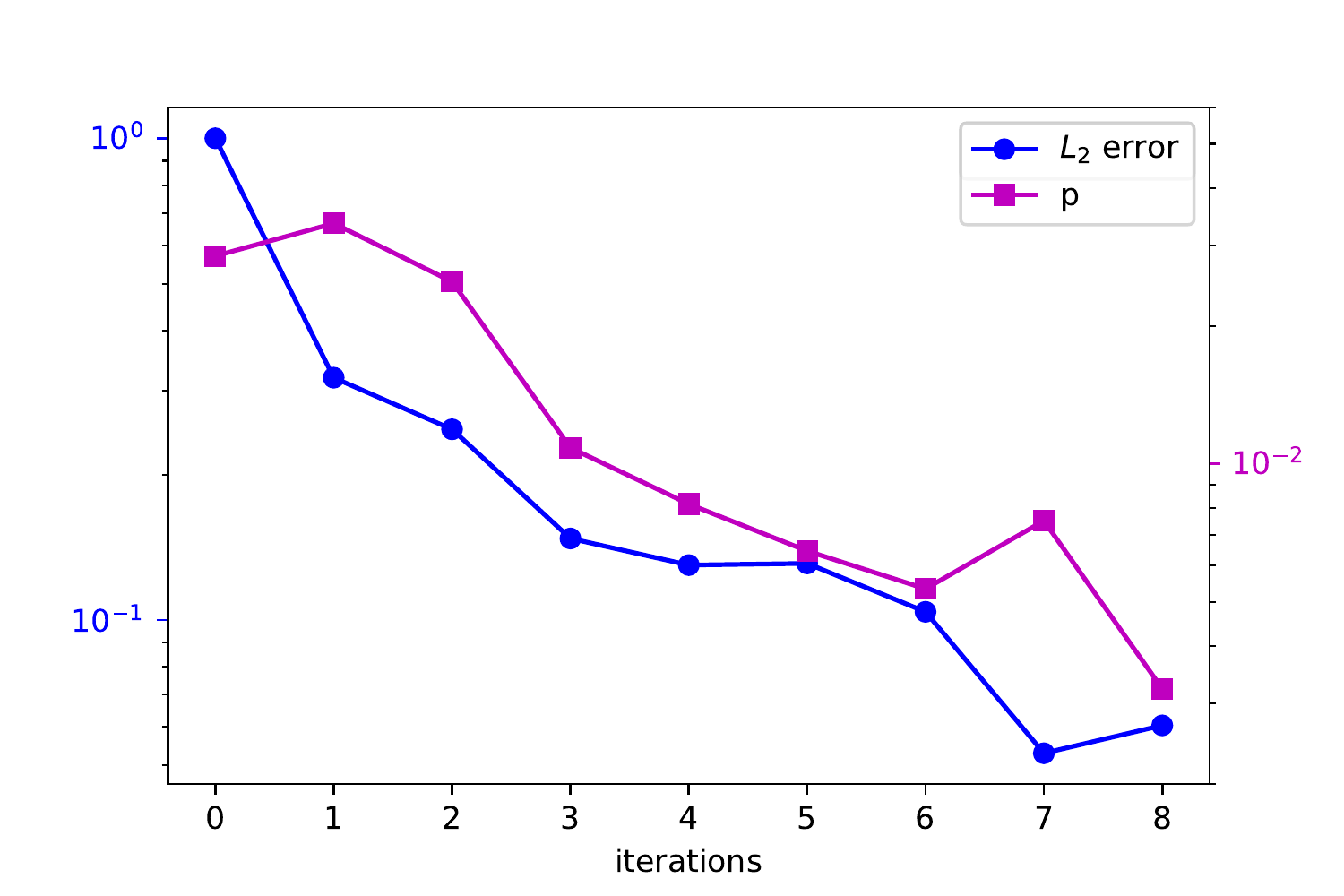}
    \end{overpic}
    \end{center}
    \caption{Relative $L_{2}$ error (left) and the  estimated failure probability (Right);  $d = 9$.}
    \label{high_dimensioanl_error}
    \end{figure}

The left plot of Fig. \ref{high_dimensioanl_error} displays the relative errors achieved  using three different sampling strategies. It is shown that the error decreases smoothly when employing \textit{SAIS}, whereas the \textit{Uniform} and \textit{RAR} seems not converge for this example.  This demonstrates that \textit{SAIS} can outperform  other techniques in terms of efficiency. The right plot of Fig.\ref{high_dimensioanl_error}  further shows  that the estimated failure probability has a similar trend with the mean $L_{2}$ error, confirming the efficiency of the stop criteria of FI-PINNs.

\begin{figure}
    \begin{center}
             \begin{overpic}[width=0.7\textwidth,trim= 35 0 45 15, clip=true,tics=10]{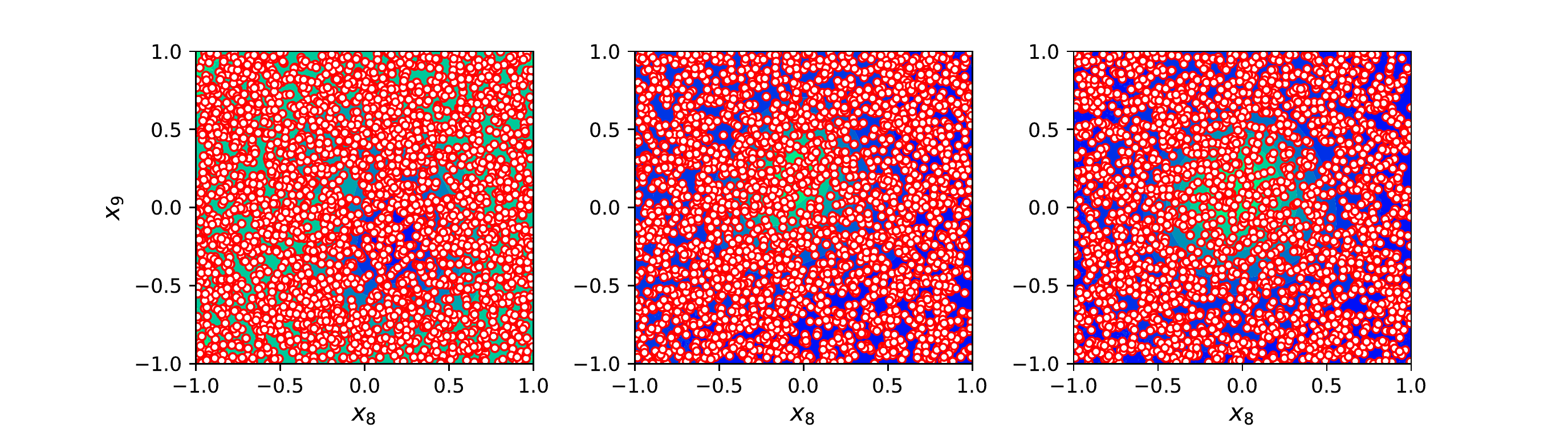}
             \put (34,29) {\scriptsize {\bf Samples  by \textit{Uniform}}}
            \end{overpic}
        \end{center}

        \vspace{0.3cm}
      \begin{center}
      \begin{overpic}[width=0.7\textwidth,trim= 35 0 45 15, clip=true,tics=10]{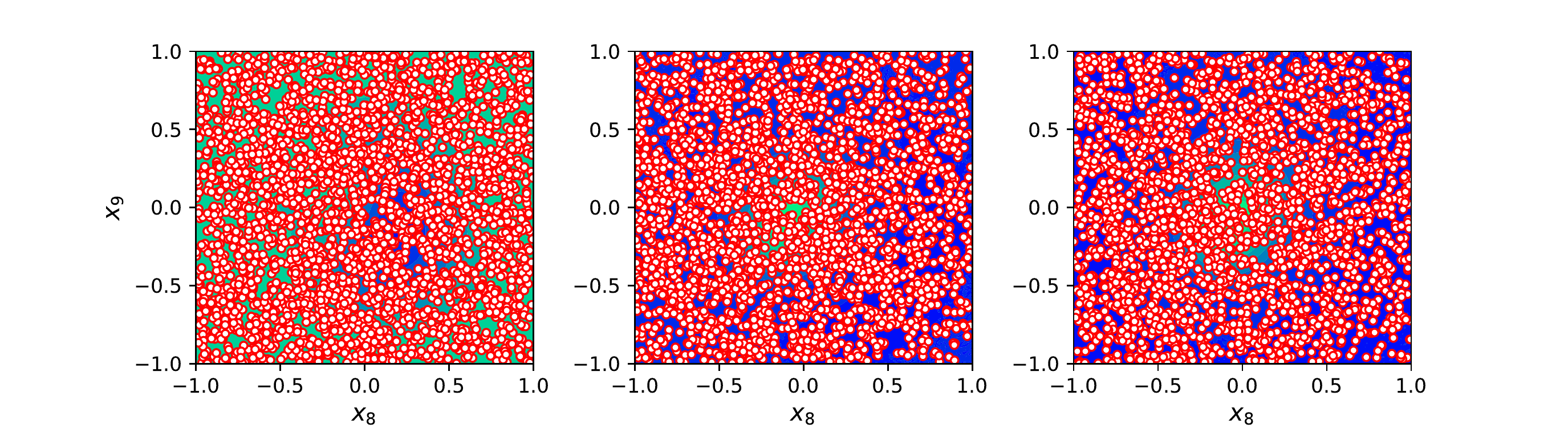}
        \put (34,29) {\scriptsize {\bf Samples  by \textit{RAR}}}
       \end{overpic}
    \end{center}
    \vspace{0.3cm}
          \begin{center}
          \begin{overpic}[width=0.7\textwidth,trim= 35 0 45 15, clip=true,tics=10]{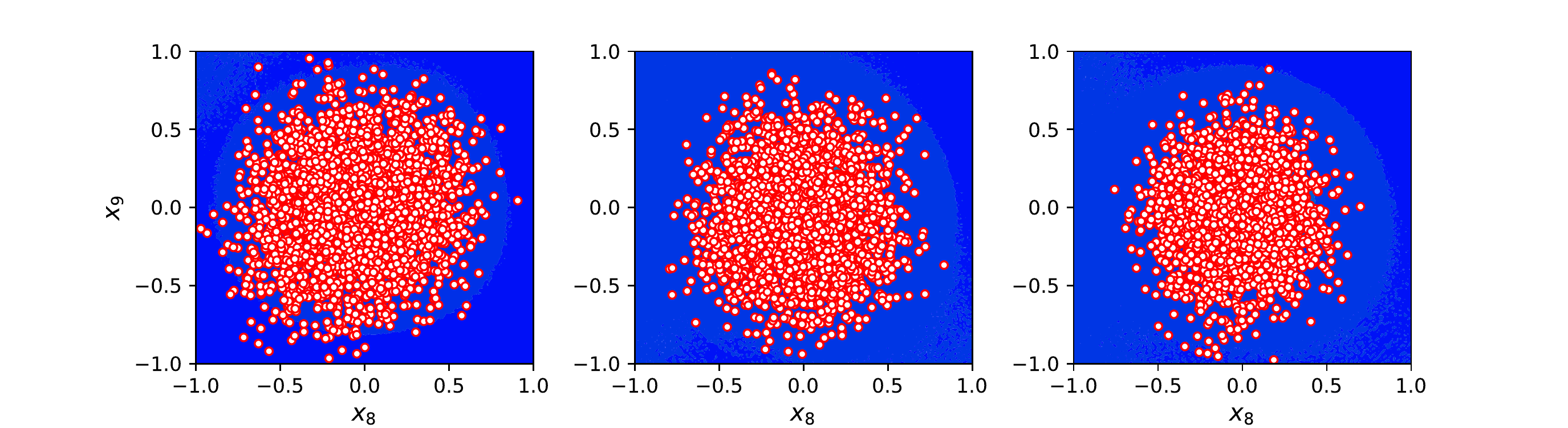}
            \put (34,29) {\scriptsize {\bf Samples  by \textit{SAIS}}}
           \end{overpic}
      \end{center}
 \caption{The $2nd, 4th, 6th$(from left to right) samples distribution of $x_{8}, x_{9}$  obtained by the three sampling strategies respectively; $d = 9$ }
    \label{high_dimensioanl_samples}
      \end{figure}
The new generated samples for the $2nd, 4th, 6th$ iterations in the space of $x_{8}, x_{9}$ are shown in Fig. \ref{high_dimensioanl_samples}. Again, it is clearly seen that the samples generated by \textit{SAIS} are more concentrated around $(0,0)$, which significantly increases the effective sample size to train the network and results in a smaller $L_{2}$ error; see also Fig.\ref{high_dimensioanl_error}.

 \subsection{Two-dimensional unbounded problem}

 Consider the following two-dimensional Possion problem
 \begin{equation}
     \label{unbounded_2d_problem}
     \begin{split}
     -\Delta u(x, y) = f(x, y), \quad (x, y)\in \mathbb{R}^2\setminus \Omega \\
     u(x,y) = s(x, y),\quad (x, y) \in \partial \Omega,
     \end{split}
 \end{equation}
 where the true solution is
 \begin{equation}
     \label{unbounded_2d_true_solution}
     u(x, y) = e^{-(x - 4)^{2} - (y - 4)^{2}}.
 \end{equation}
The boundary of $\Omega$ is defined as
 \[\partial\Omega  = \left(\cos(t) - \frac{\cos(5t)\cos(t)}{4}, \sin(t) - \frac{\cos(5t)\sin(t)}{4}\right), \quad  t\in [0, 2\pi].\]
In this unbounded example, we set the initial  proposal distributions $h_{1}$ to be a two-dimensional Gaussian distribution with mean vector $\mu_{1} = [0,0]$, covariance matrix $\Sigma_{1} = 3I_2$.  Additionally, we set $\epsilon_{r}$ and  $\epsilon_{p}$ in FI-PINNs to be 0.01 and 0.0001, respectively.
We set the hyperparameter of the \textit{SAIS} algorithm to $N_{2} = 1000$.

 \begin{figure}[t]
    \centering
    \includegraphics[width = 1\textwidth]{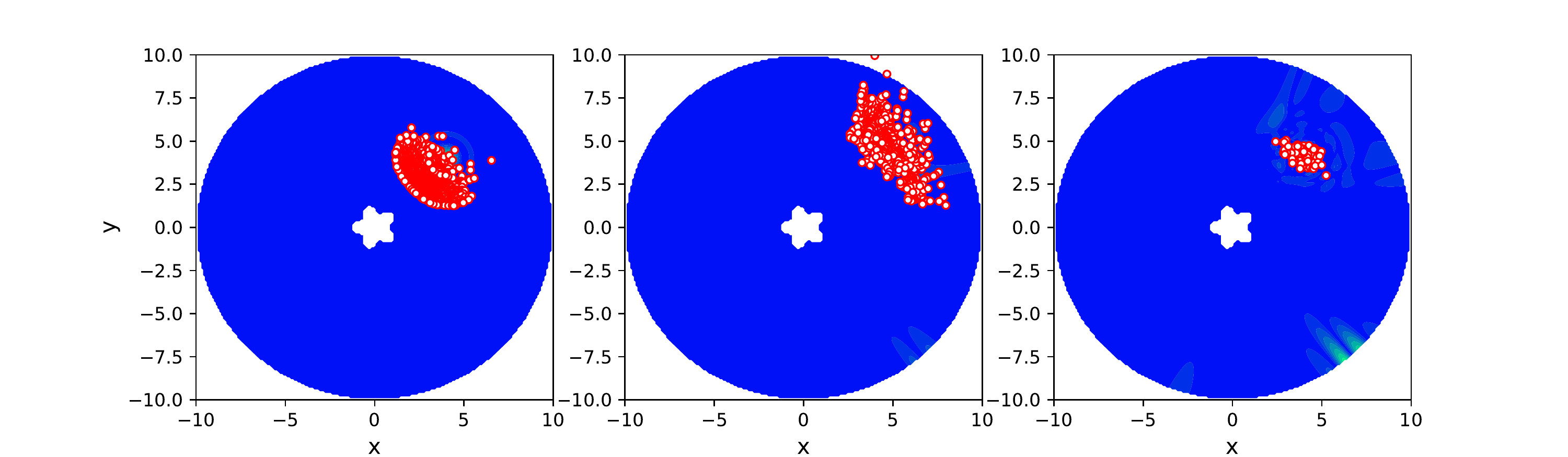}
    \caption{The distribution of collocation points obtained by SAIS for the first three updates.}
    \label{unbounded_2d_samples}
    \end{figure}

The challenge of selecting a training data  for unbounded domain arises from the fact that the true solution is unknown. To verify the effectiveness of our method, we starting with $5000$ collocation points randomly chosen in $[-2,2]^2\setminus \Omega$ which is far away from the peak. The estimated failure probability becomes smaller than the tolerance after 4 updates. The distribution of newly added collocation points obtained by \textit{SAIS}  for the first three iterations is shown in Fig.\ref{unbounded_2d_samples}. It is seen that the update points produced by \textit{SAIS} are concentrated near the peak $(4,4)$.

 \begin{figure}[t]
    \begin{center}
    \begin{overpic}[width=0.32\textwidth, clip=true,tics=10]{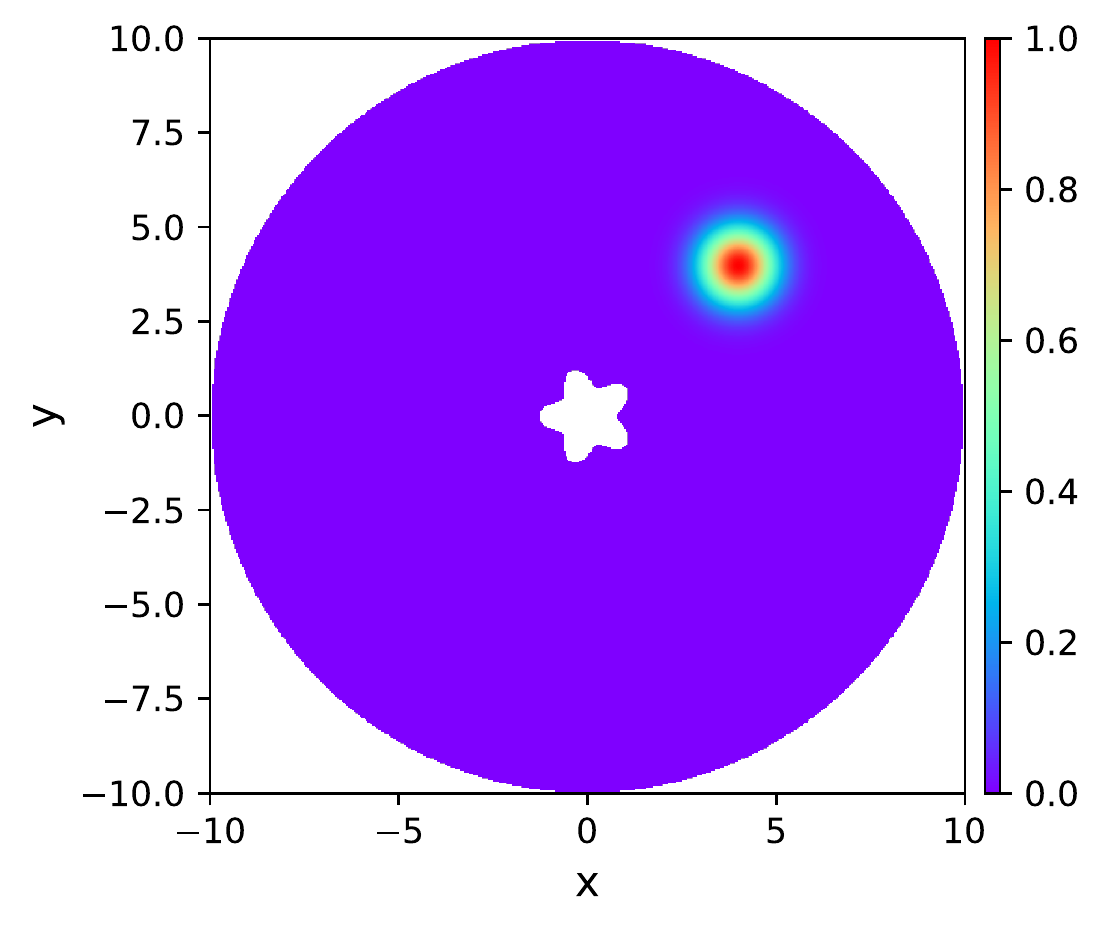}
    % \put (30,30) {\scriptsize {\bf Samples distribution by  ARA}}
    \end{overpic}
    % \hspace{-0.5cm}
    \begin{overpic}[width=0.32\textwidth, clip=true,tics=10]{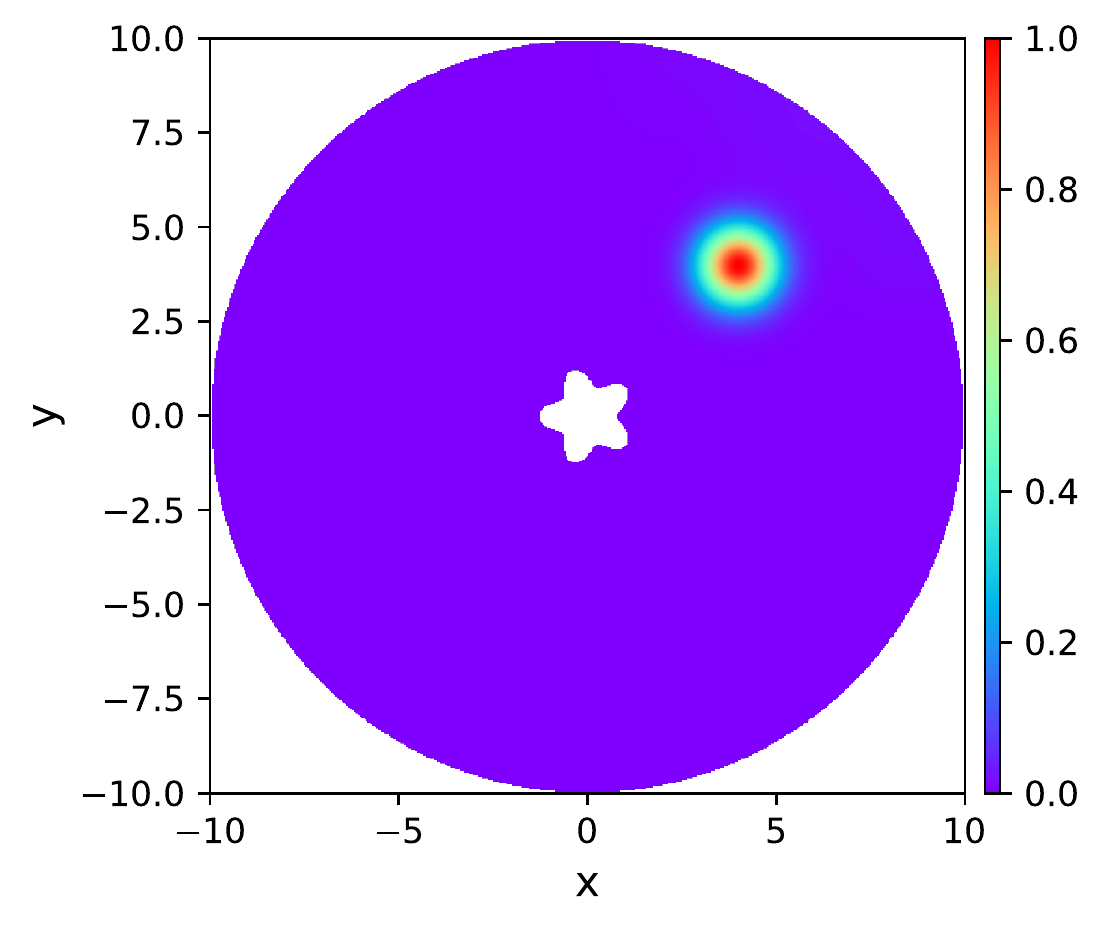}
    % \put (30,30) {\scriptsize {\bf Samples distribution by  ARA}}
    \end{overpic}
    \begin{overpic}[width=0.34\textwidth, clip=true,tics=10]{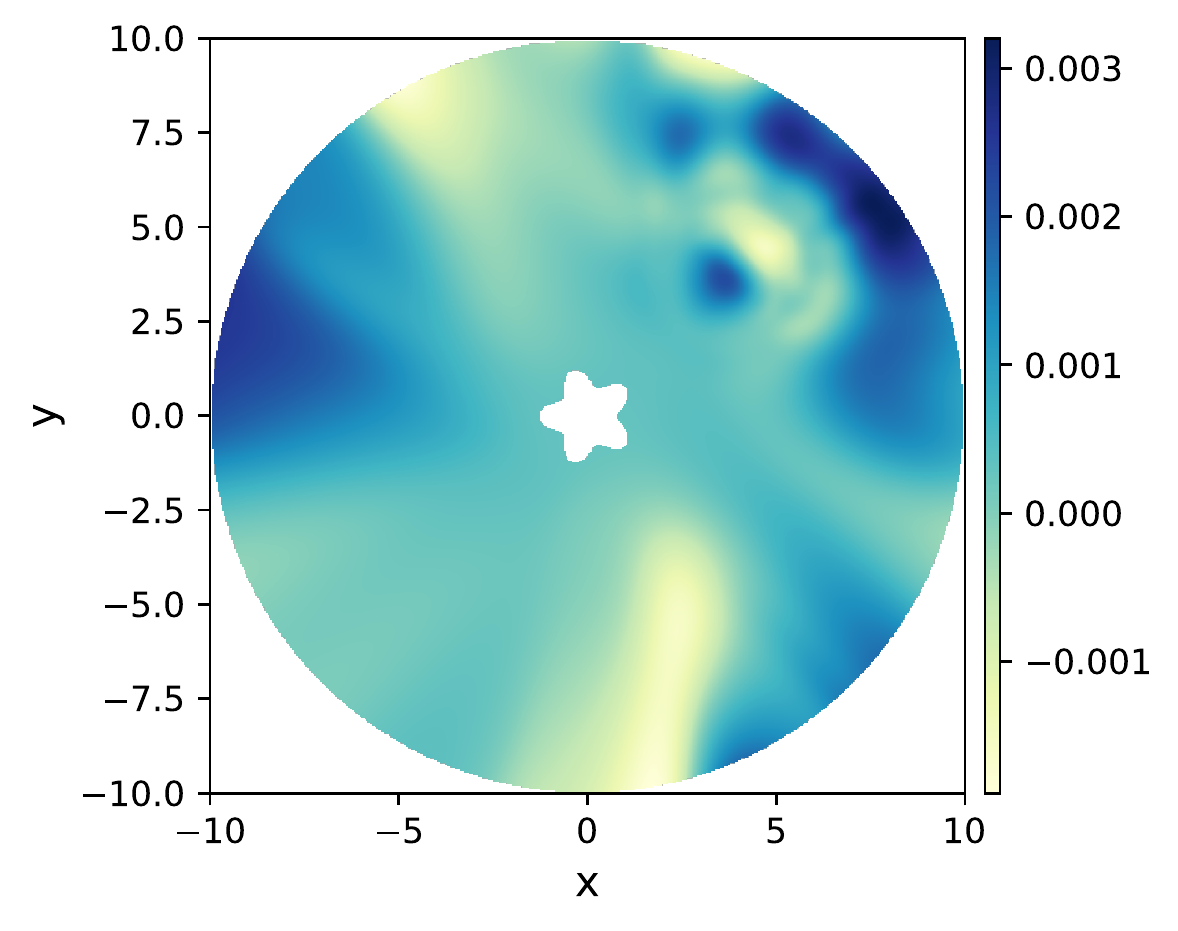}
        % \put (30,30) {\scriptsize {\bf Samples distribution by  ARA}}
        \end{overpic}
    \end{center}
%    \begin{center}
%        \begin{overpic}[width=1\textwidth, clip=true,tics=10]{figures/solution_x_unbounded_2d.pdf}
%        \end{overpic}
%    \end{center}
    \vspace{-0.2cm}
    \caption{The  exact solution, numerical solution and absolute error from left to right.}
    \label{unbounded_2d_solution}
    \end{figure}
    \begin{figure}[t]
        \centering
        \includegraphics[width = 0.5\textwidth]{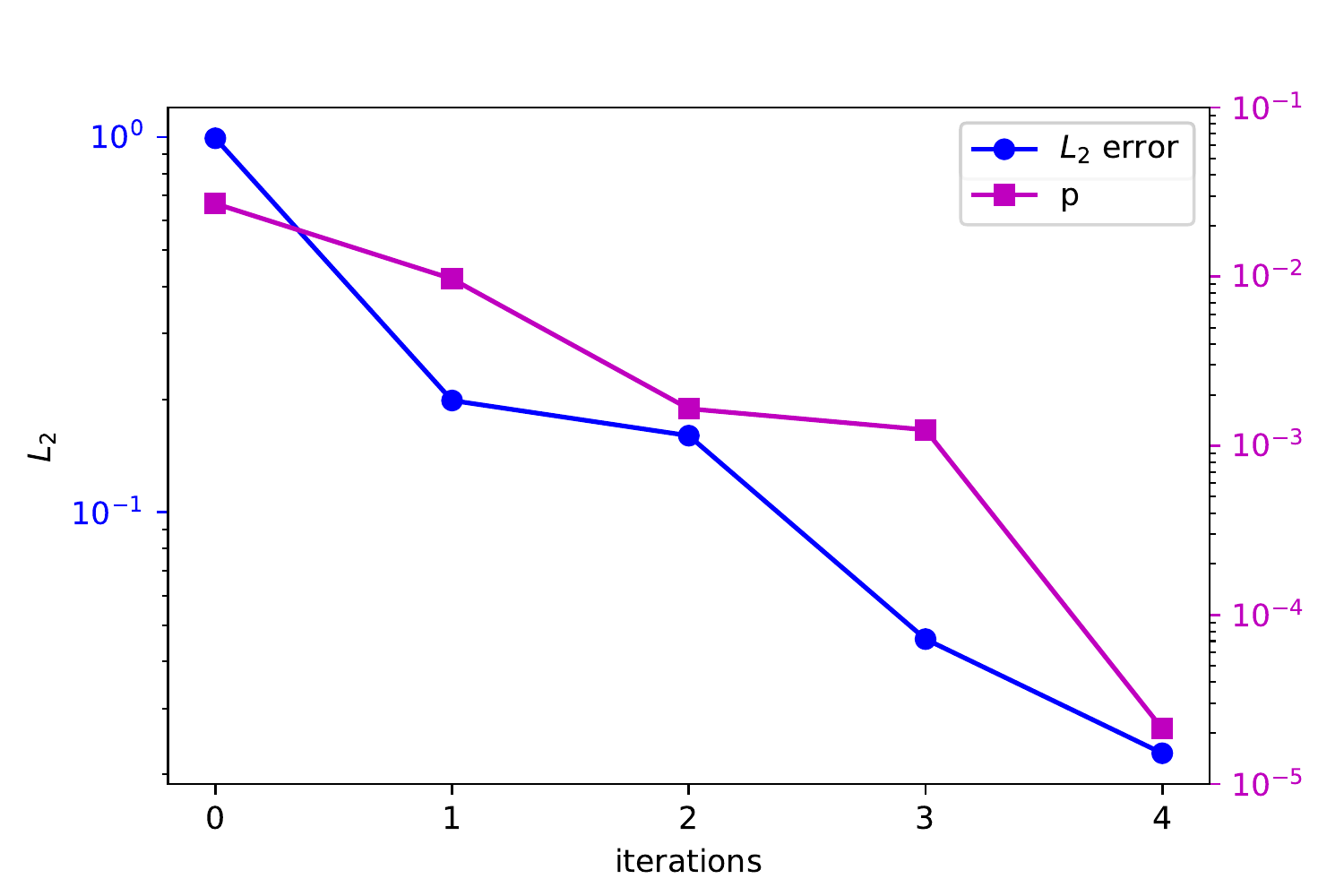}
        \caption{The relative $L_{2}$ error and the  estimated failure probability.}
        \label{unbounded_2d_error}
        \end{figure}

Fig. \ref{unbounded_2d_solution} shows the predicted solution over a disk with a radius of 10. It is seen that the numerical solution match well with the exact solution.
We also provide the  $L_{2}$ errors  in Fig. \ref{unbounded_2d_error} to demonstrate the efficiency of our adaptive framework. The figure shows that, after four updates, our technique can produce a  prediction error of only $2.28\times 10^{-2}$. Additionally, the prediction error and predicted failure probability have a similar pattern, making it useful for creating stop criteria.

\subsection{Time-dependent problem in an unbounded  domain}
We finally consider the following time-dependent equation
\begin{equation}
    \label{ubnounded_time_dependent_problem}
    \begin{split}
        &u_t(x, t) =  u_{xx}(x, t) + f(x, t), \quad (x, t)\in \mathbb{R} \times[0,1] \\
        &u(x, 0) = u_0(x), \quad x \in \mathbb{R}.
    \end{split}
\end{equation}
We choose the following reference solution
\begin{equation}
    \label{unbounded_time_dipendent_true_solution}
    u(x, t) = \frac{e^{- \frac{\left(x - 10\right)^{2}}{4 t + 4}}}{\sqrt{t + 1}},
\end{equation}

\begin{figure}[t]
    \centering
    \includegraphics[width = 0.55\textwidth]{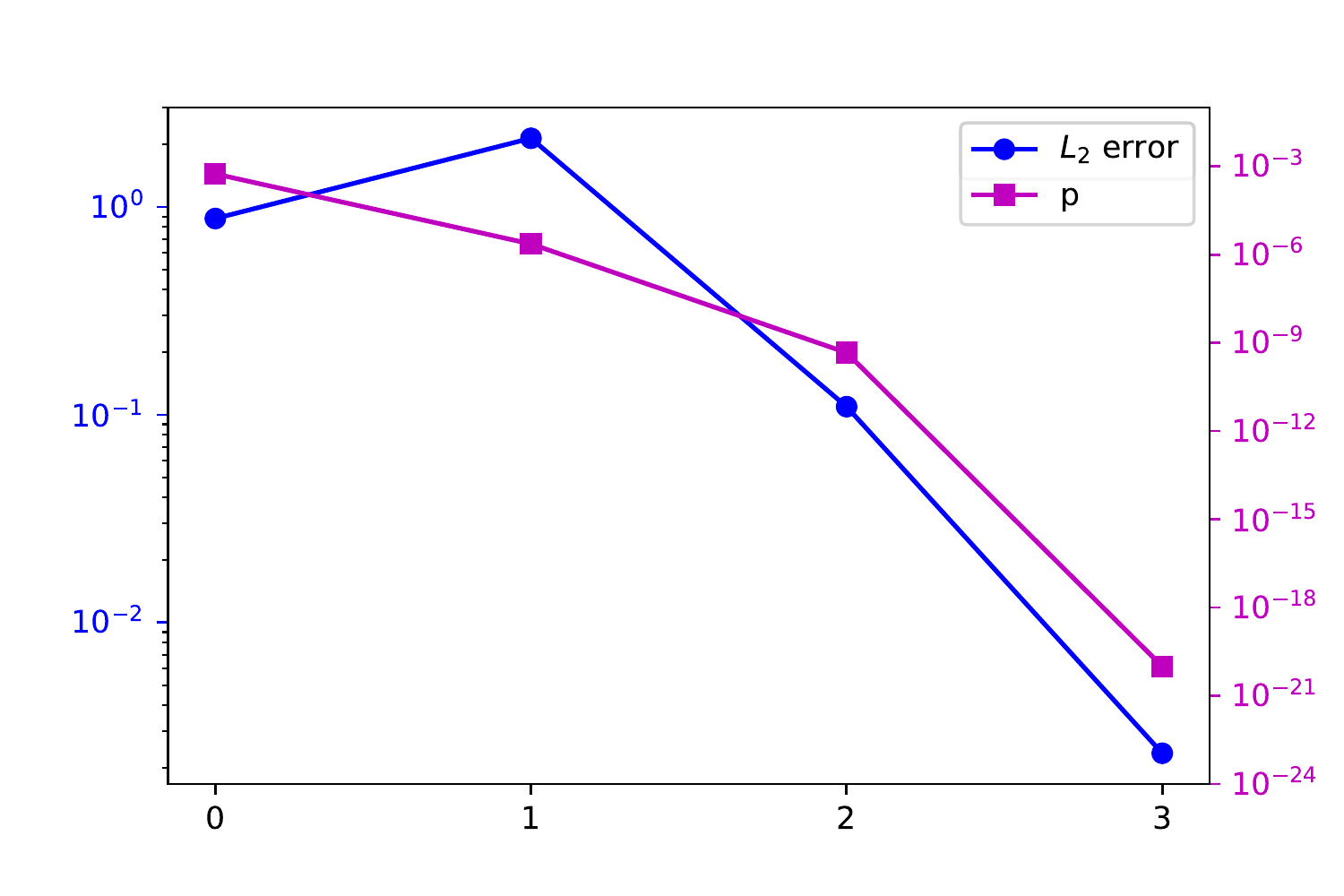}
    \vspace{-0.2cm}
    \caption{The relative $L_{2}$ error and the  estimated failure probability.}
    \label{unbounded_error}
    \end{figure}

In this example, we set the tolerances $\epsilon_{r} = 0.01$ and $\epsilon_{p} = 10^{-10}$, respectively. The hyper-parameters used in the \textit{SAIS} algorithm are set to $N_{1} = 600, N_{2} = 2000, p_{0}= 0.05$. The prior distribution of $(x, t)$ is a two dimensional Gaussian distribution with mean vector $\mu_{1} = [0, 0]$ and covariance matrix $\Sigma_{1} = 3I_{2}$.  Moreover, the initial 3000 training points are generated in $[-6, 0]\times[0,1]$, which is far away from the peak.   The training can be terminated after three updates to the training dataset because the estimated failure probability is smaller than the tolerance (see Fig.\ref{unbounded_error}).

\begin{figure}[t]
    \begin{center}
        \begin{overpic}[width=0.242\textwidth, clip=true,tics=10]{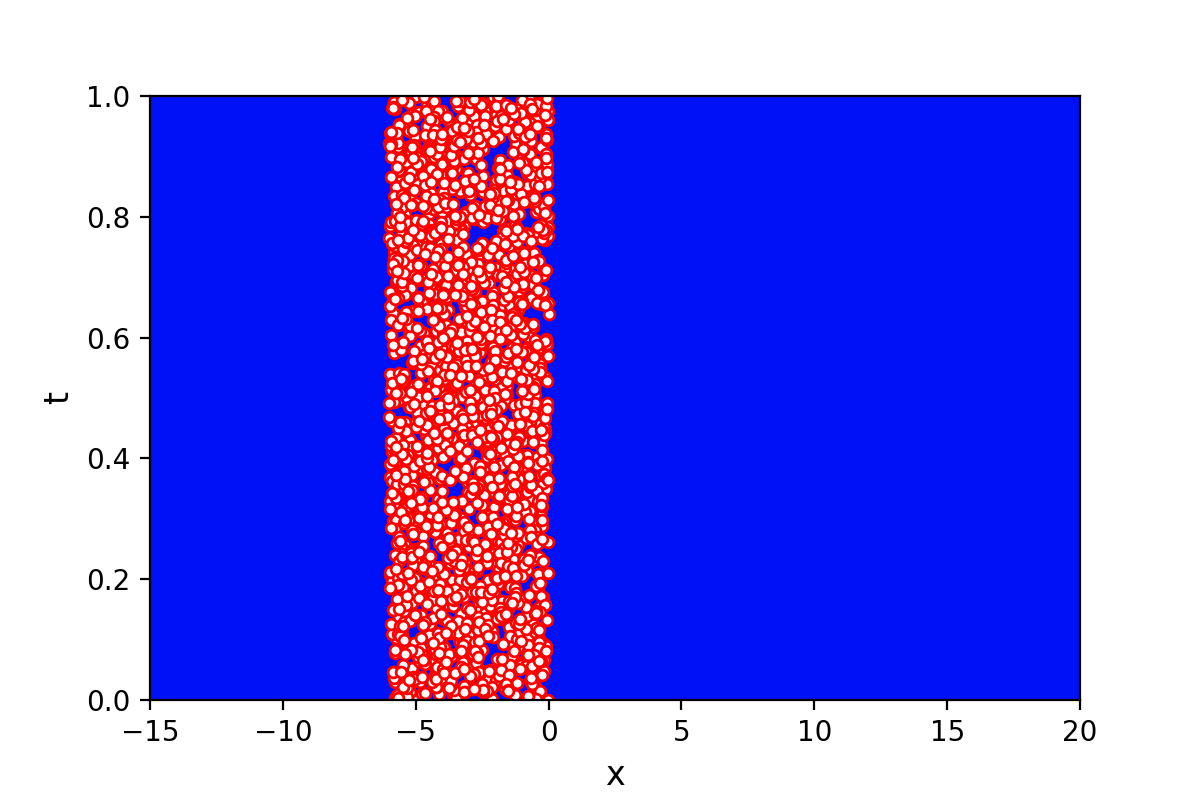}
        \put (25,65) {\scriptsize {\bf Initial  points}}
        \end{overpic}
     \begin{overpic}[width=0.242\textwidth, clip=true,tics=10]{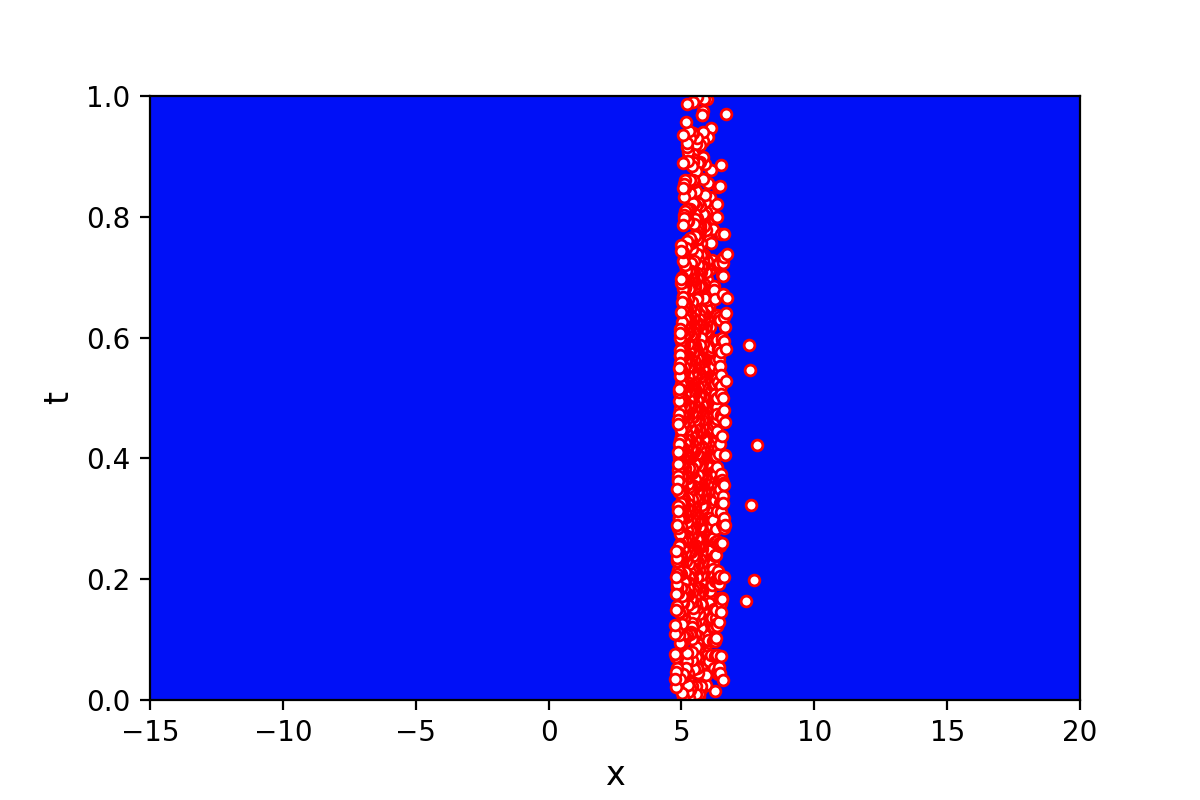}
             \put (35,65) {\scriptsize {\bf 1-updated}}
        \end{overpic}
             \begin{overpic}[width=0.242\textwidth, clip=true,tics=10]{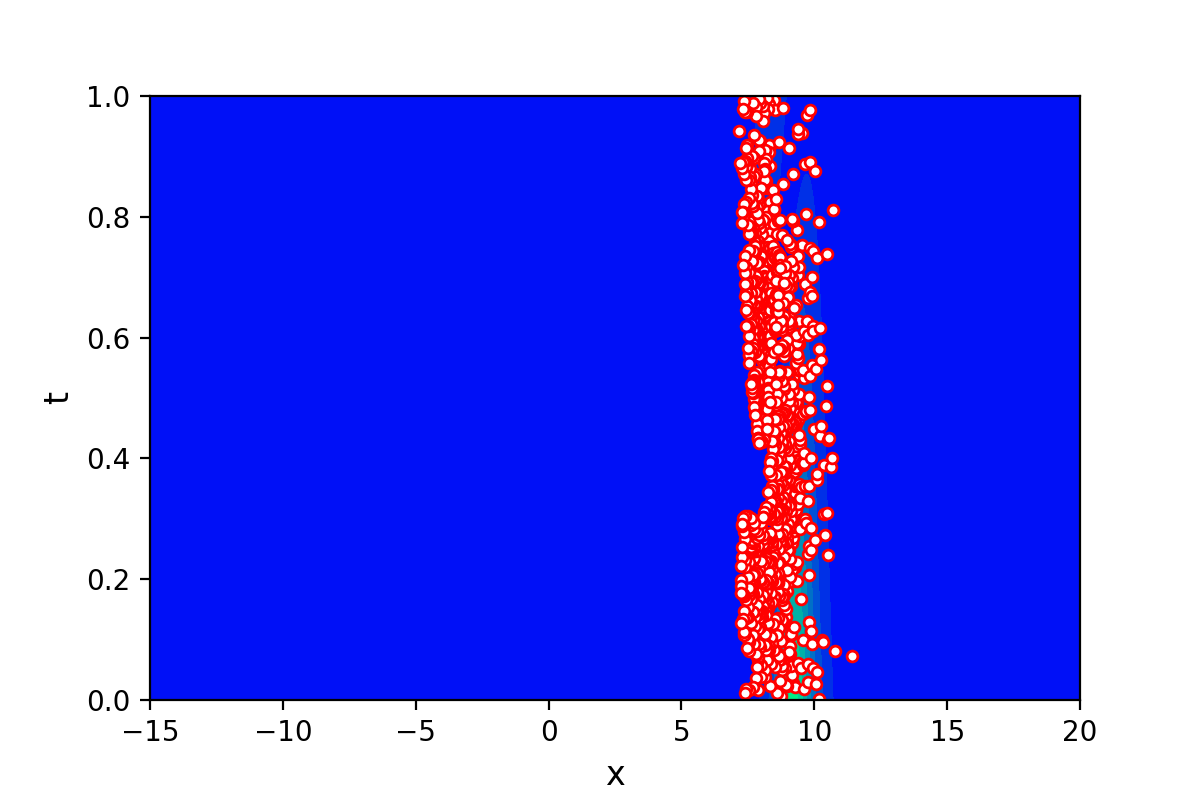}
             \put (35,65) {\scriptsize {\bf 2-updated}}
        \end{overpic}
   \begin{overpic}[width=0.242\textwidth, clip=true,tics=10]{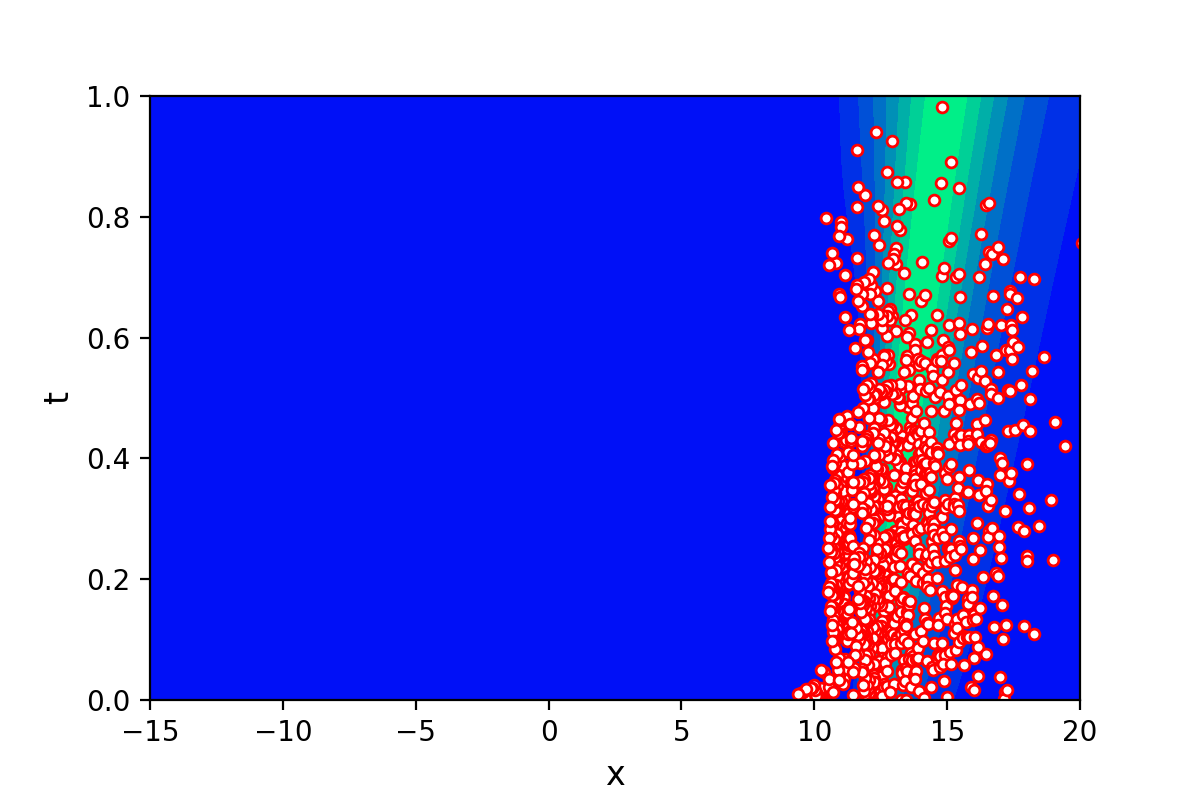}
           \put (35,65) {\scriptsize {\bf 3-updated}}
        \end{overpic}
\end{center}
    \begin{center}
     \begin{overpic}[width=0.242\textwidth, clip=true,tics=10]{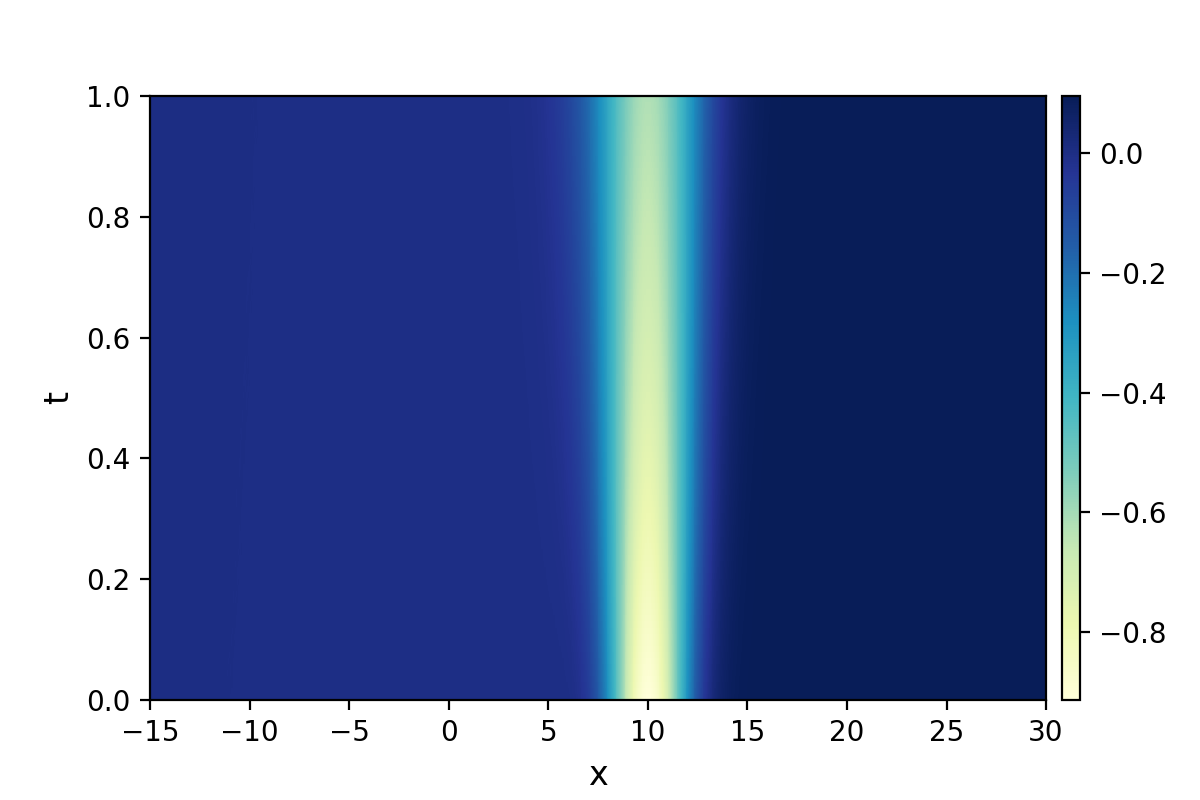}
        \end{overpic}
              \begin{overpic}[width=0.242\textwidth, clip=true,tics=10]{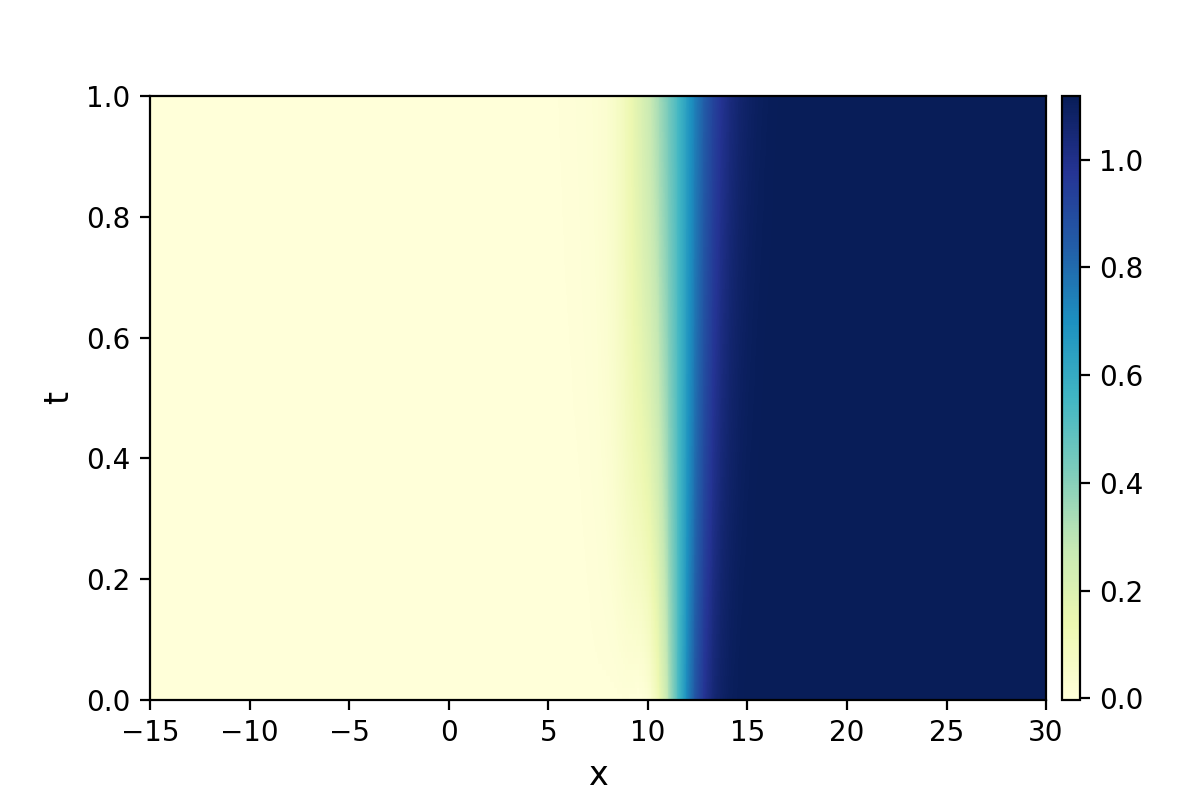}
        \end{overpic}
                      \begin{overpic}[width=0.242\textwidth, clip=true,tics=10]{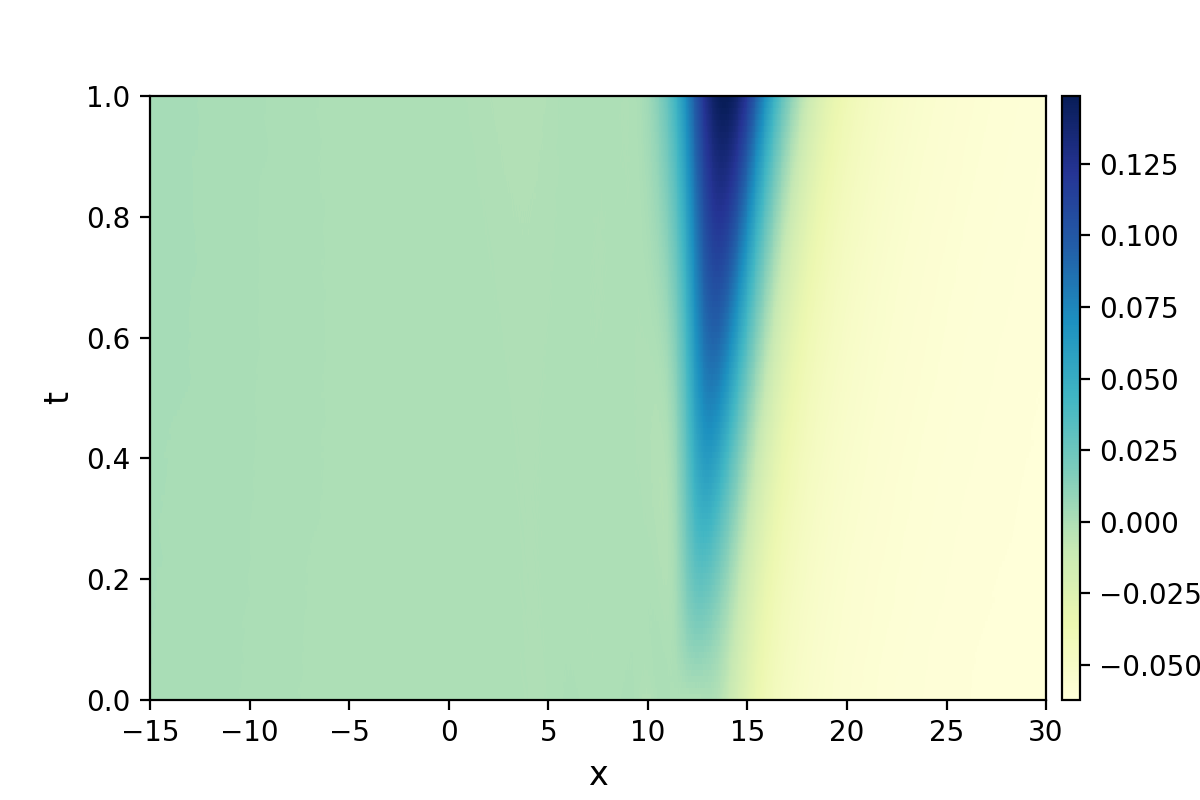}
        \end{overpic}
                \begin{overpic}[width=0.242\textwidth, clip=true,tics=10]{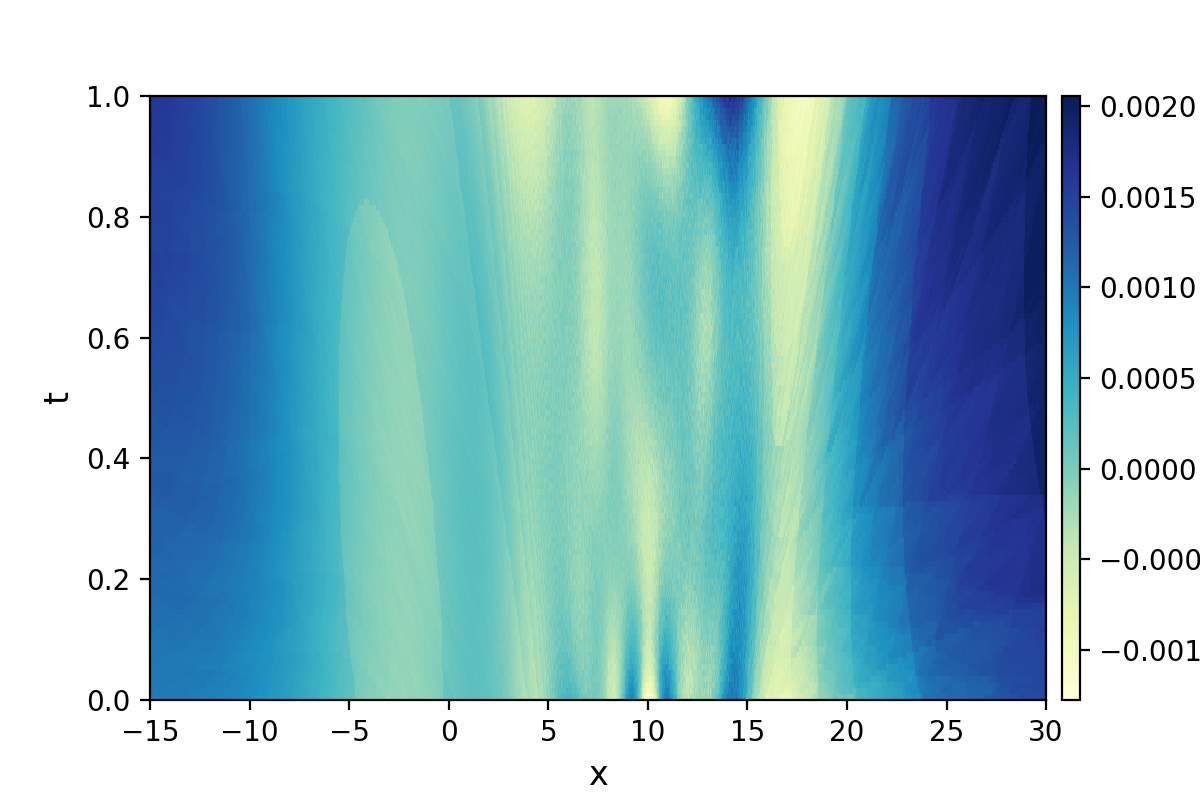}
        \end{overpic}
    \end{center}
    \begin{center}
           \begin{overpic}[width=0.242\textwidth, clip=true,tics=10]{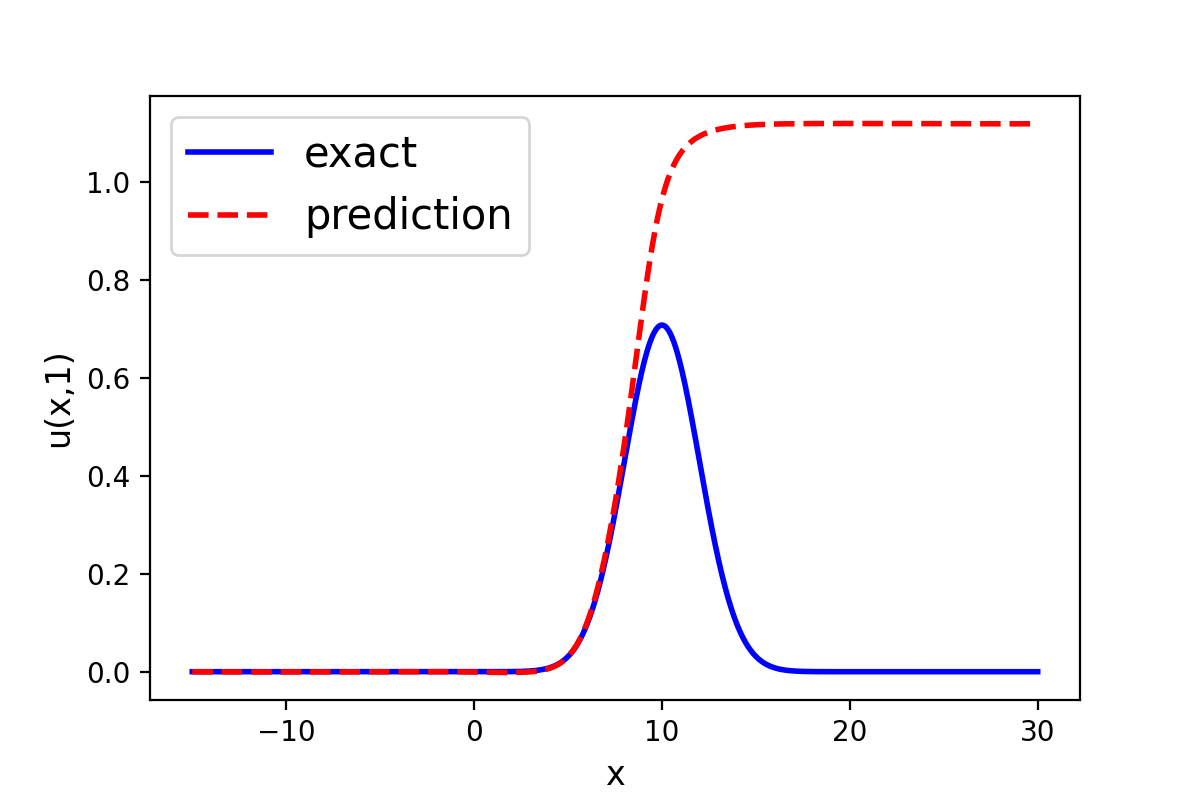}
        \end{overpic}
 \begin{overpic}[width=0.242\textwidth, clip=true,tics=10]{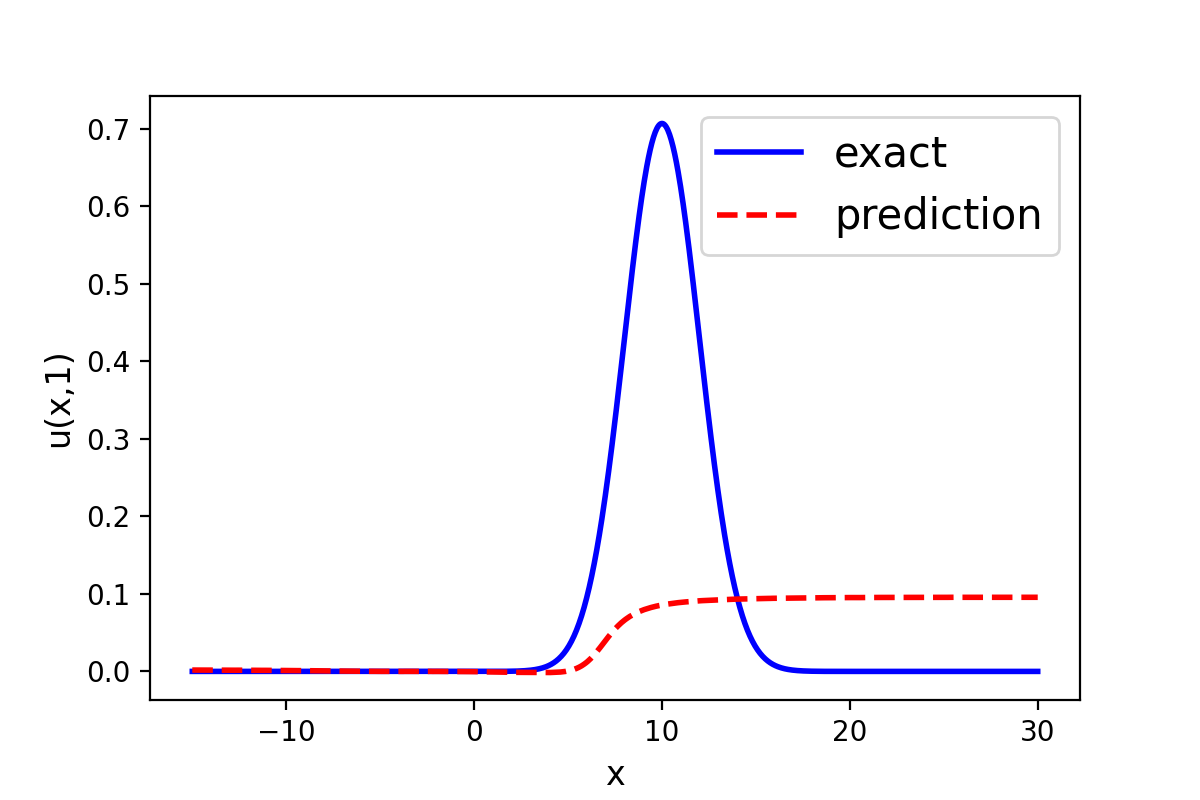}
        \end{overpic}
                \begin{overpic}[width=0.242\textwidth, clip=true,tics=10]{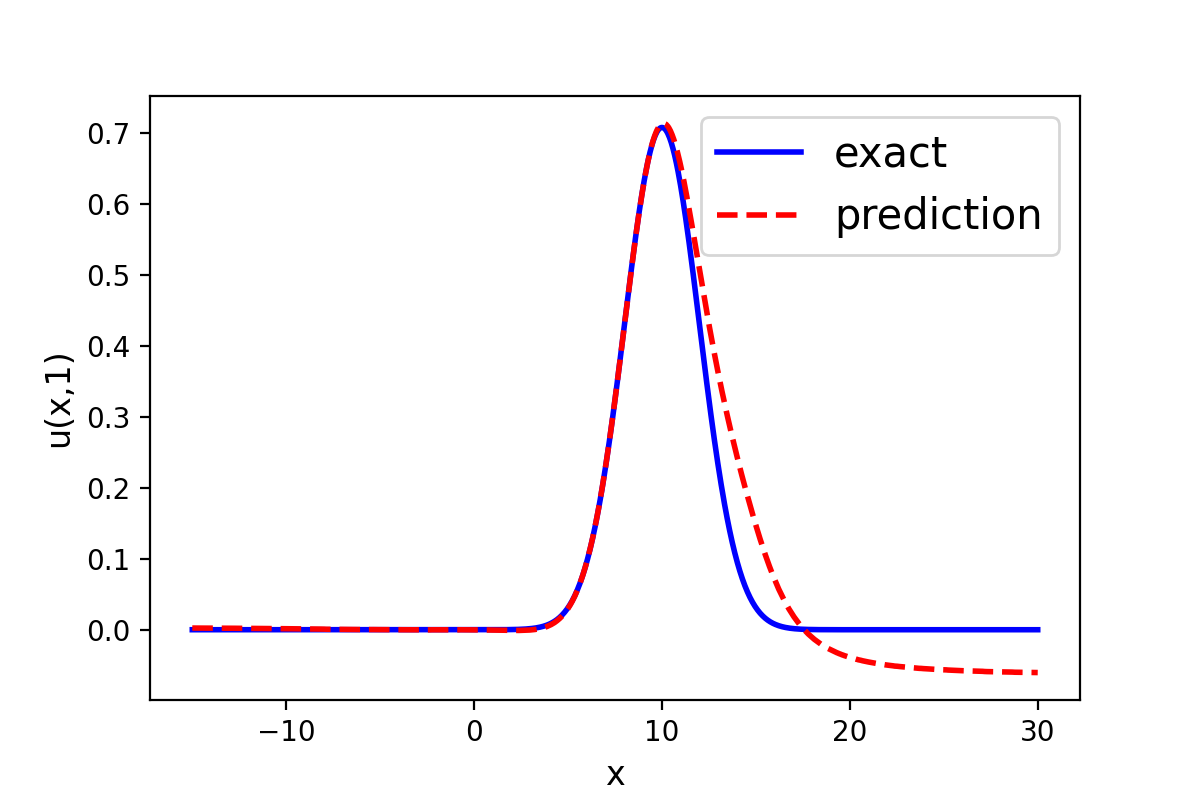}
        \end{overpic}
        \begin{overpic}[width=0.242\textwidth, clip=true,tics=10]{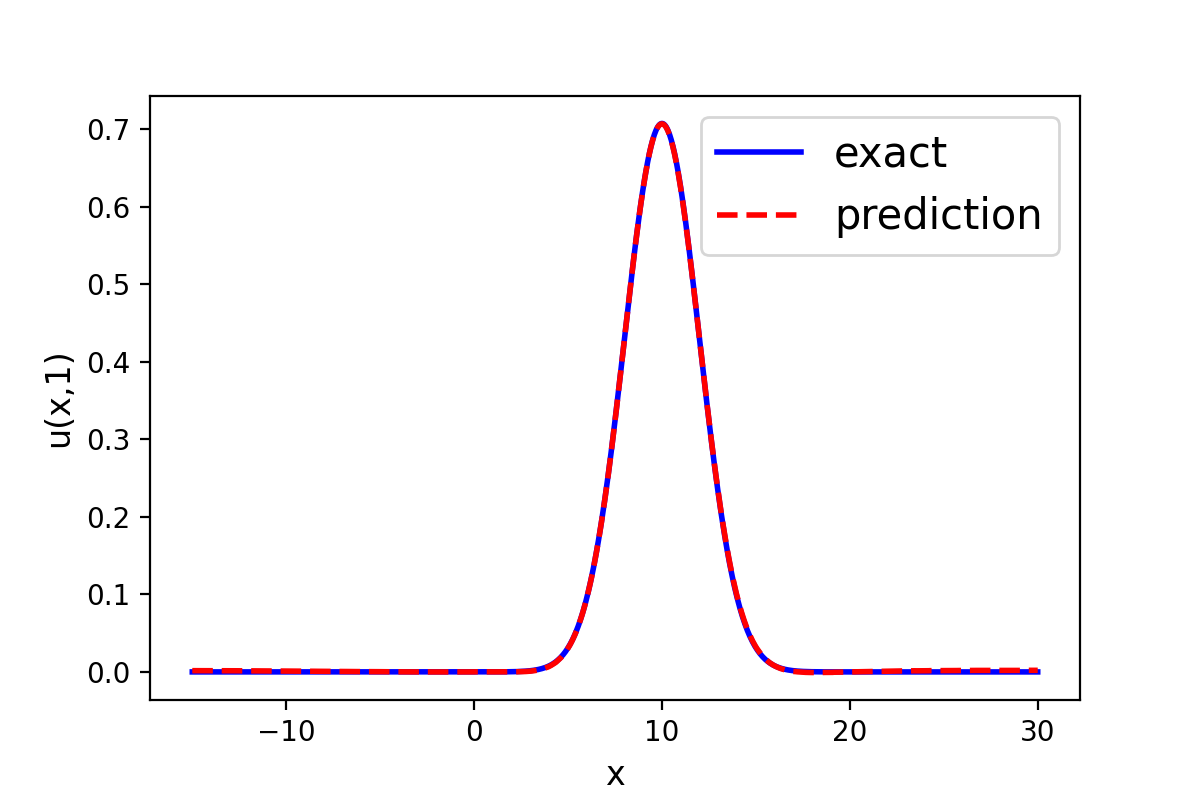}
        \end{overpic}
    \end{center}
    \caption{From top to bottom: the distribution of (updated) samples, the corresponding absolute error and the predicted values at $t = 1$. }
    \label{unbounded_samples}
\end{figure}

The distributions of the  collocation points are shown in Fig.\ref{unbounded_samples}.  We can clearly see that the collocation points produced using the  \textit{SAIS} algorithm automatically move toward the region where the residual error is the greatest.  This phenomenon demonstrates the effectiveness of our FI-PINNs framework in dealing with unbounded time-dependent problems.  Moreover, the absolute error showed in Fig.\ref{unbounded_samples}  gradually decreases during the training process. Finally,  we can obtain an absolute error less than $2\times 10^{-3}$ over the whole test domain. To be more clear, we plot the predicted values at $t = 1$. It is obvious that the initial predicted solution differs a lot from the exact solution. While as the training continues, the predicted solution gradually converges to the true solution.  After three updates, the relative prediction error shown in Fig.\ref{unbounded_error} is smaller than $3\times 10^{-3}$.  The corresponding failure probability shares the same trend. This phenomenon indicates that our FI-PINNs framework and the SAIS sampling strategy are suitable for this kinds of time-dependent problems as well.

\section{Concluding remarks}
\label{sec:con}

In this work, we present failure-informed PINNs that combine the PINNs with adaptive sampling procedures. The key idea in our FI-PINNs is to define a failure probability as an error indicator, based on the residual. In particular, we propose to use truncated Gaussian as a simple model to estimate the failure probability and to generate new training points. Rigorous error bound of FI-PINNs is also presented. We illustrate the performance FI-PINNs via several PDE problems that include PDEs with singular solutions,  PDEs in unbounded domains, and time-dependent PDEs. It is shown that for all the test problems, FI-PINNs can effectively capture the solution structure. Further studies include applications of FI-PINN to more complex real problems.  We remark again here that our truncated Gaussian model can also be changed into more complex models such as Gaussian mixture model. We also mention again that DAS \cite{tang2021deep} is in fact a more general model for this purpose, for which the main issue is the computational complexity. An optimal way for selecting the density model by balancing the effectiveness and the efficiency will be part of our future studies.

\section*{Acknowledgments}
LY's work was supported by the NSF of China (No.12171085). TZ's work was supported by the National Key R$\&$D Program of China (No. 2020YFA0712000), the NSF of China (under grant numbers 11822111, 11688101 and 11731006),  the Strategic Priority Research Program of Chinese Academy of Sciences (No.  XDA25000404), and youth innovation promotion association (CAS).

\appendix

\section{Further discusses}
In this section, we will show two straightforward extensions of the FI-PINNs framework. We shall show in Section \ref{mGau} that more general density models such as Gaussian mixture model can be used in our framework.  We shall also discuss in Section \ref{ACeq} how the causality-based weighted residual can be used as a performance function.

\subsection{SAIS using Gaussian mixture model} \label{mGau}

The key idea here is to utilize a mixed Gaussian model as the intermediate proposal distribution $h_k$ in our SAIS algorithm.  Note that the density function for the Gaussian mixture distribution usually takes the following form:
$$
\sum^{M}_{m=1} \pi_m \mathcal{N}(\mu_{m}, \Sigma_{m})
$$
where $\mathcal{N}(\mu_{m}, \Sigma_{m})$ is a Gaussian probability density function with mean value $\mu_{m}$ and covariance matrix $\Sigma_{m}$. For every $m$, $\pi_m\in(0,1)$ satisfy $\sum^{M}_{m=1}\pi_m=1$.
Similar to truncated Gaussian model, we also employ  an adaptive method that starts with an initial proposal and incrementally updates the intermediate proposal $h_k$ using samples $\mathcal{D}_k:=\{\widetilde{\mb{x}}_{i}^{k}\}_{i=1}^{N_{1}}$.   How to extract $\mathcal{D}_k$ from the intermediate proposal is completely explained in Section \ref{sec:sais}.  Once we have the candidate set $\mathcal{D}_k$, we can estimate the parameters for Gaussian mixture model $(\pi_m, \mu_{m}, \Sigma_{m})$ using the expectation-maximization (EM) algorithm \cite{bishop2006pattern}.

\begin{figure}[t]
\centering
\includegraphics[width = 0.5\linewidth]{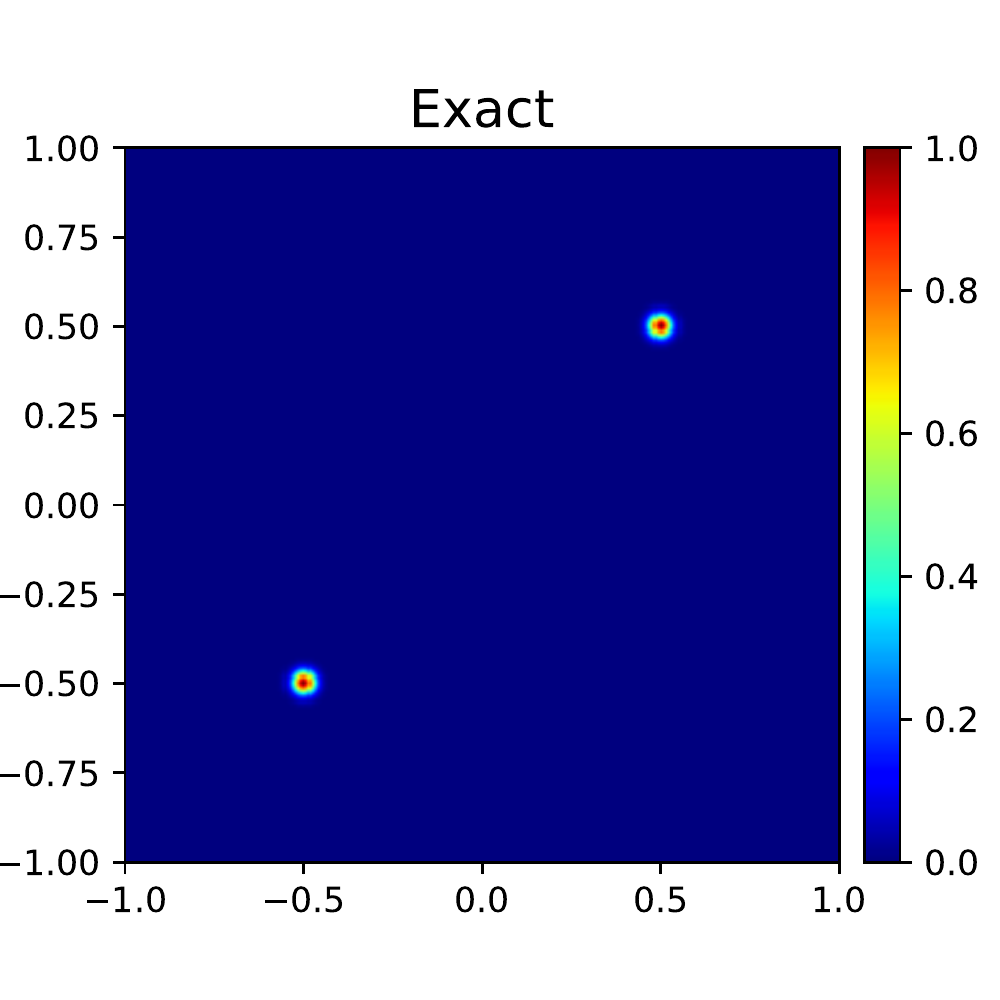}
\caption{Exact solution for the two-dimensional peak test problem.}
\label{two_peak_exact_solution}
\end{figure}

In order to demonstrate the effectiveness of the SAIS algorithm with Gaussian mixture model, we consider Example \ref{sec:52} with two peaks. We define
$$
u(x,y)=e^{-1000[(x-0.5)^2+(y-0.5)^2]}+e^{-1000[(x+0.5)^2+(y+0.5)^2]},
$$
as the exact solution, which is shown in Fig.\ref{two_peak_exact_solution}.

\begin{figure}[t]
    \begin{center}
             \begin{overpic}[width=0.9\textwidth,trim= 55 10 5 2, clip=true,tics=10]{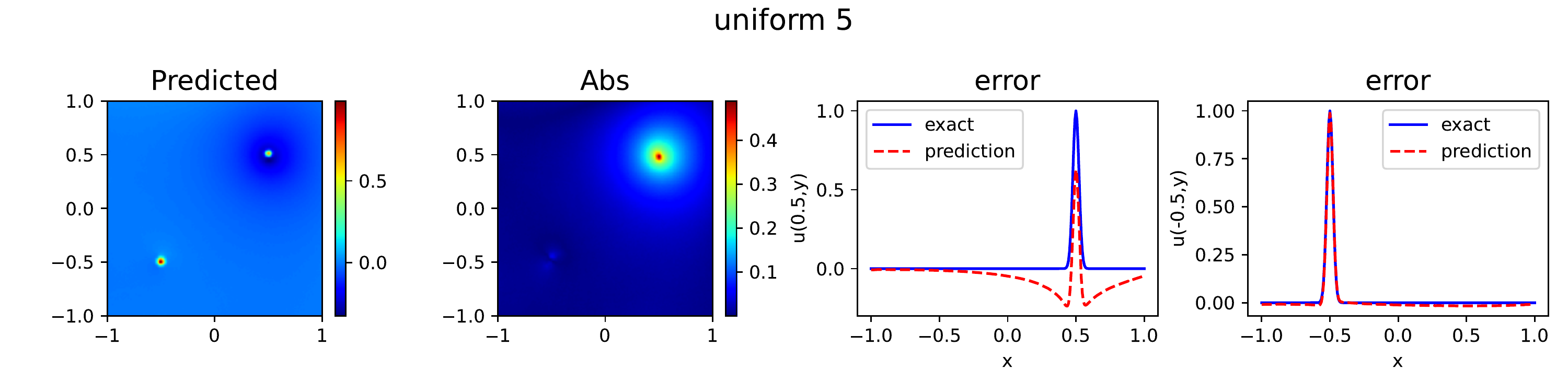}
         %    \put (30,22) {\scriptsize {\bf Numerical results  obtained by  \textit{Uniform}}}
            \end{overpic}
        \end{center}
        \vspace{0.3cm}
      \begin{center}
        \begin{overpic}[width=0.9\textwidth,trim= 55 0 5 2, clip=true,tics=10]{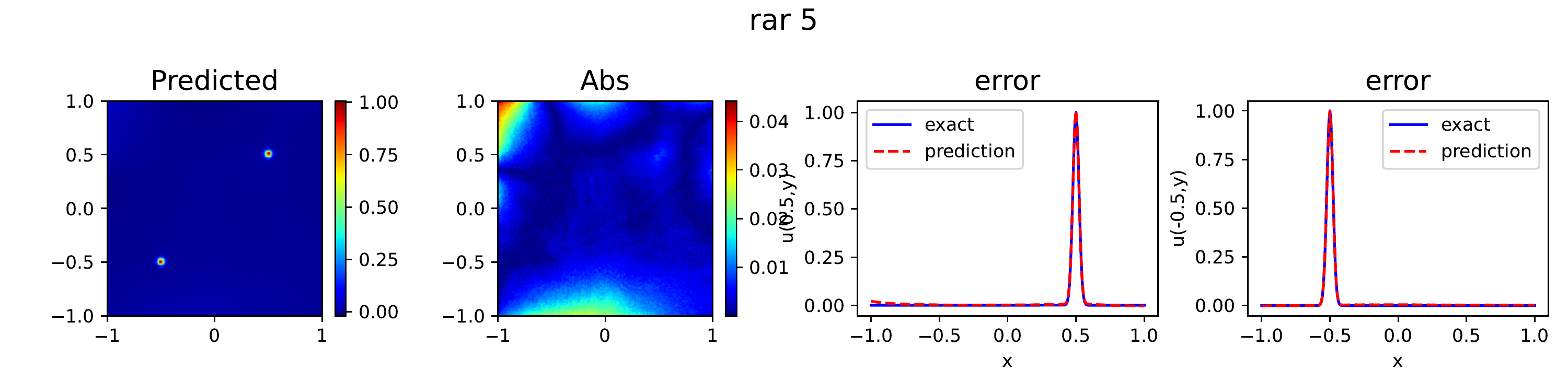}
    %      \caption{Solution errors obtained by ARA}
            % \put (30,26) {\scriptsize {\bf Numerical results  obtained by \textit{RAR}}}
      \end{overpic}
    \end{center}
    \vspace{0.3cm}
        \begin{center}
      \begin{overpic}[width=0.9\textwidth,trim= 55 0 5 2, clip=true,tics=10]{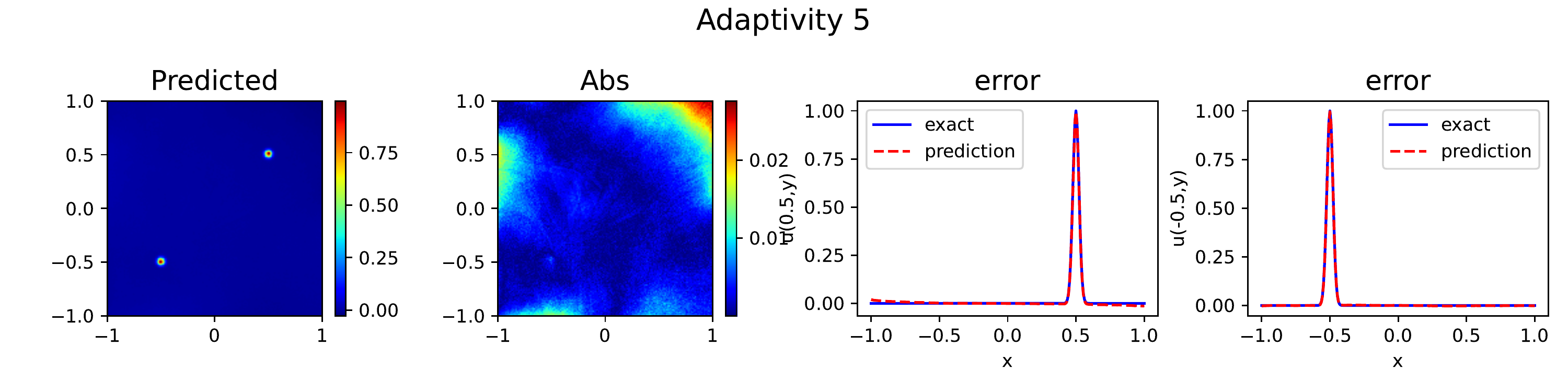}
         %  \put (30,26) {\scriptsize {\bf Numerical results  obtained  by \textit{SAIS}}}
      \end{overpic}
      \end{center}
           \vspace{-.2cm}
    \caption{The predicted solution, absolute error  and the predicted curves at $x = 0.5, -0.5$ obtained by the three different sampling strategies.}
    \label{two_peak_solution}
      \end{figure}

\begin{figure}[t]
    \begin{center}
        \begin{overpic}[width=0.242\textwidth, clip=true,tics=10]{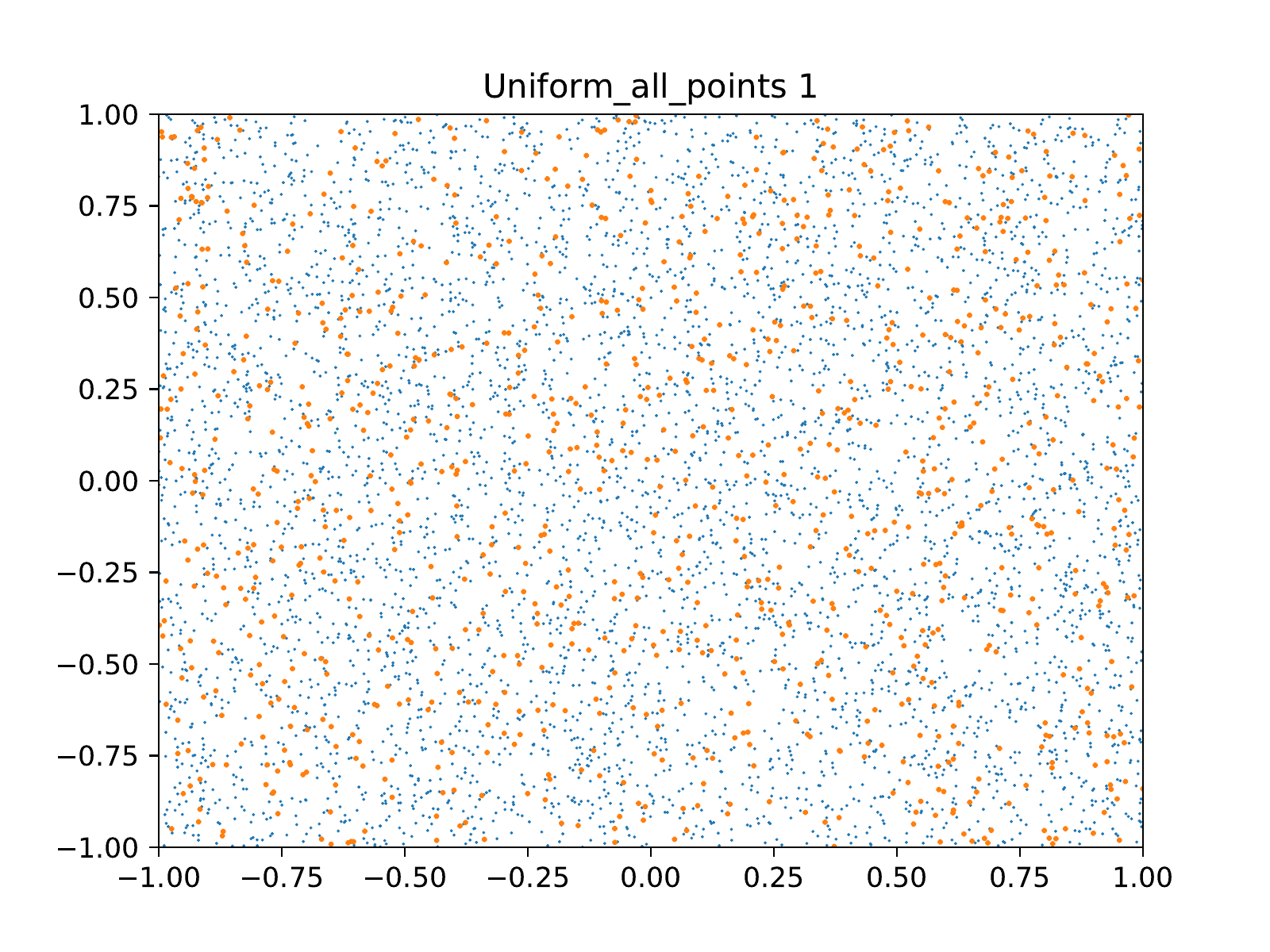}
        \put (25,72) {\scriptsize {\bf 1-updated}}
        \end{overpic}
     \begin{overpic}[width=0.242\textwidth, clip=true,tics=10]{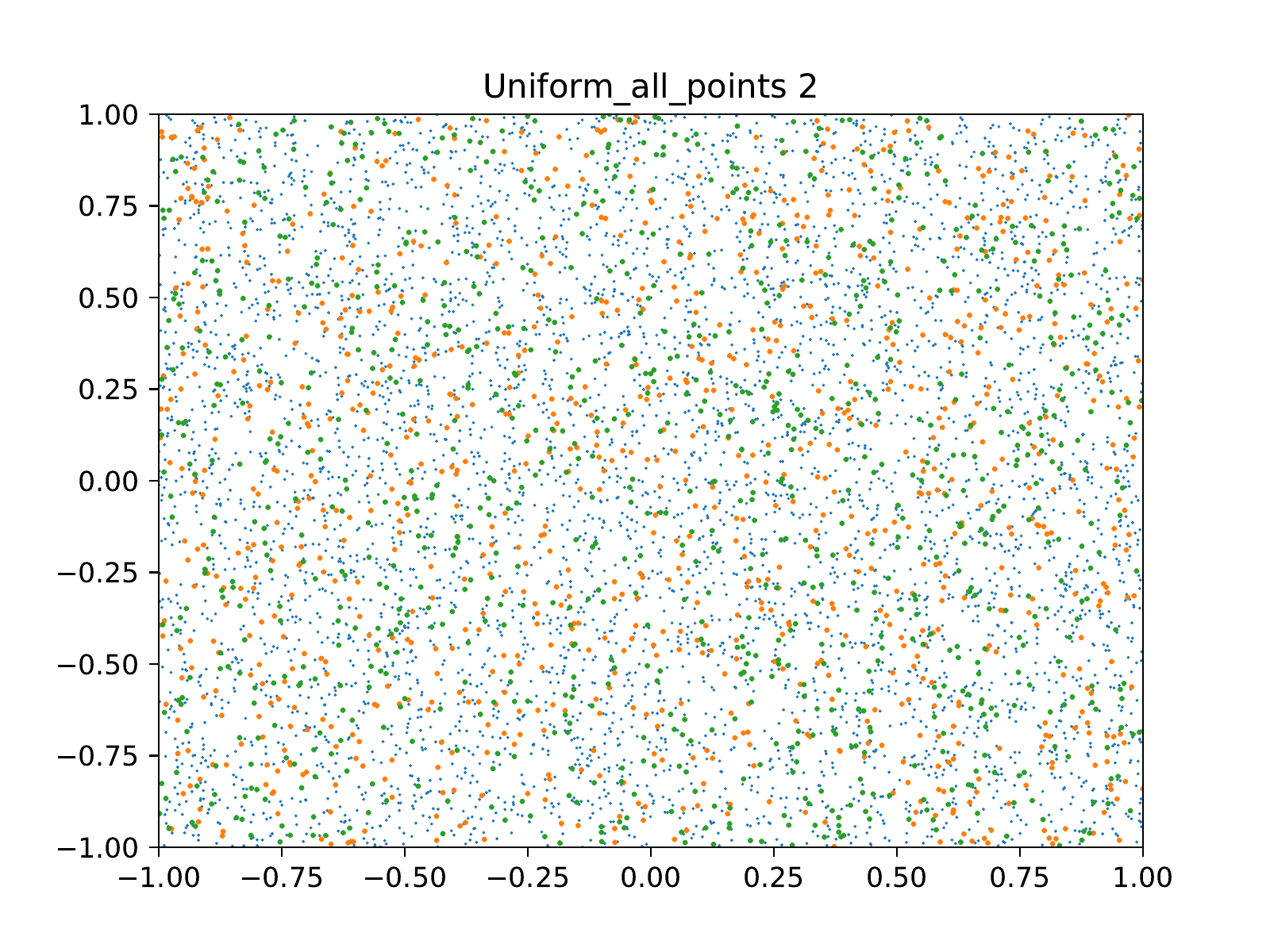}
             \put (35,72) {\scriptsize {\bf 2-updated}}
        \end{overpic}
             \begin{overpic}[width=0.242\textwidth, clip=true,tics=10]{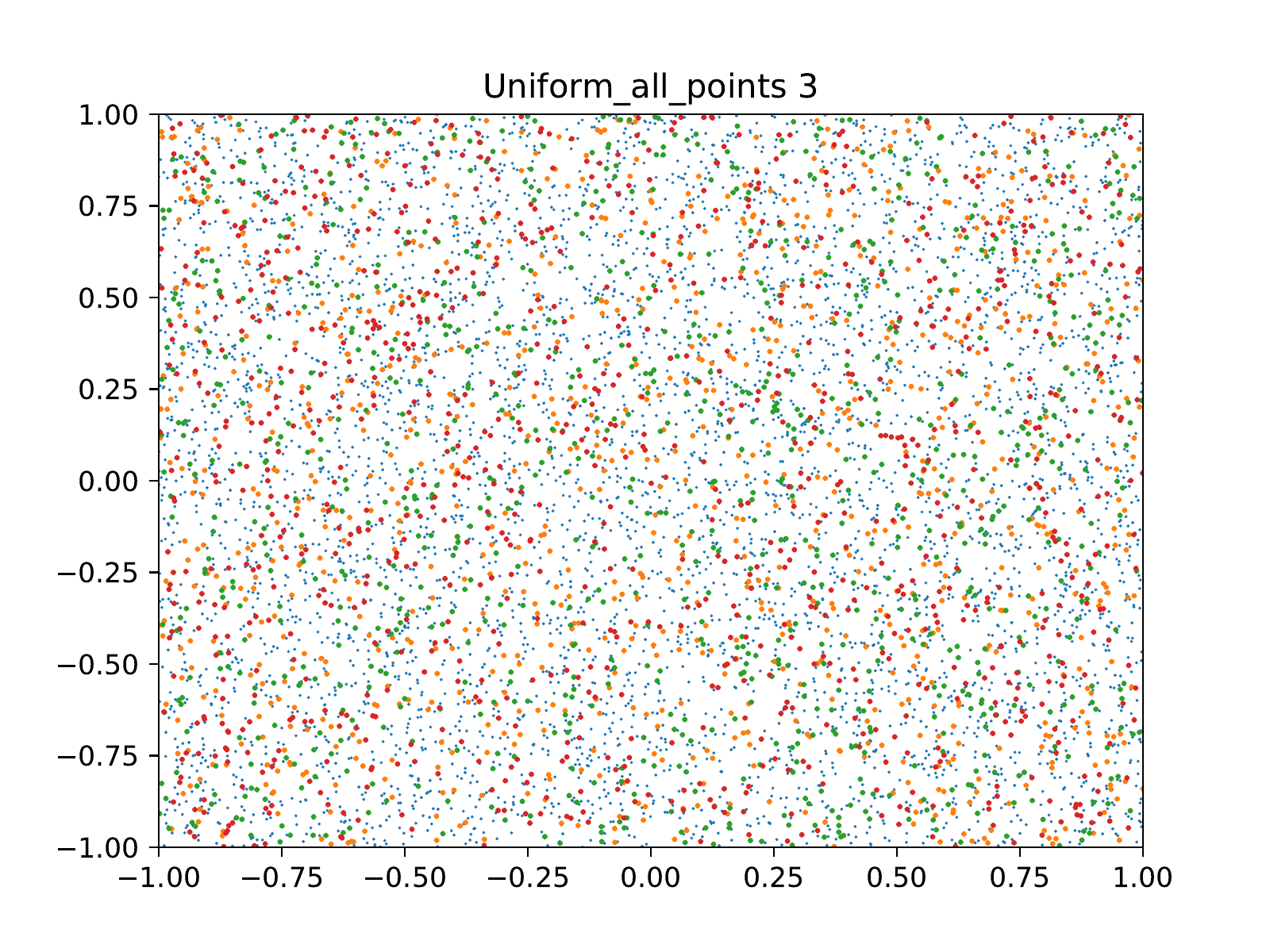}
             \put (35,72) {\scriptsize {\bf 3-updated}}
        \end{overpic}
   \begin{overpic}[width=0.242\textwidth, clip=true,tics=10]{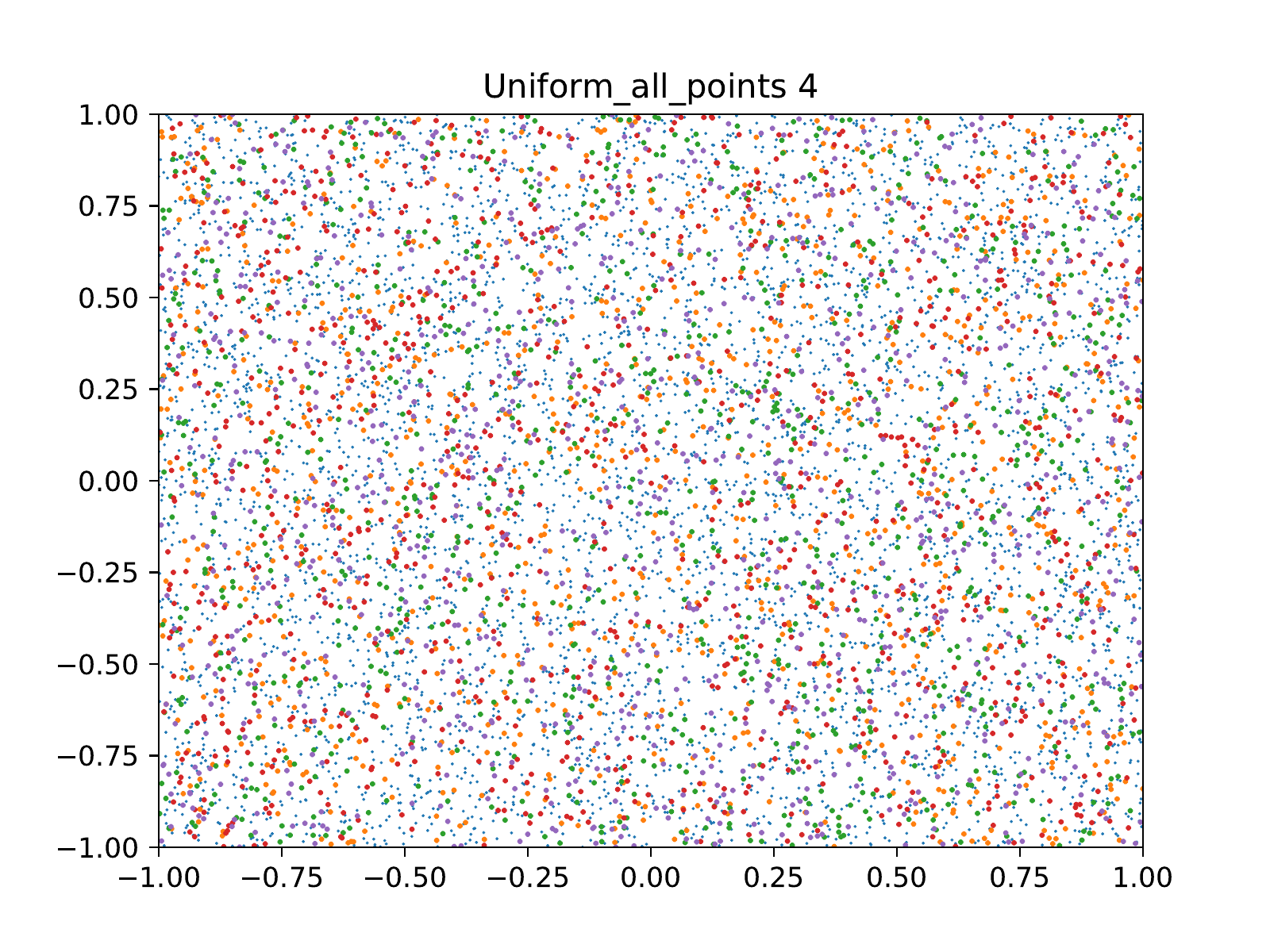}
           \put (35,72) {\scriptsize {\bf 4-updated}}
        \end{overpic}
\end{center}
    \begin{center}
     \begin{overpic}[width=0.242\textwidth, clip=true,tics=10]{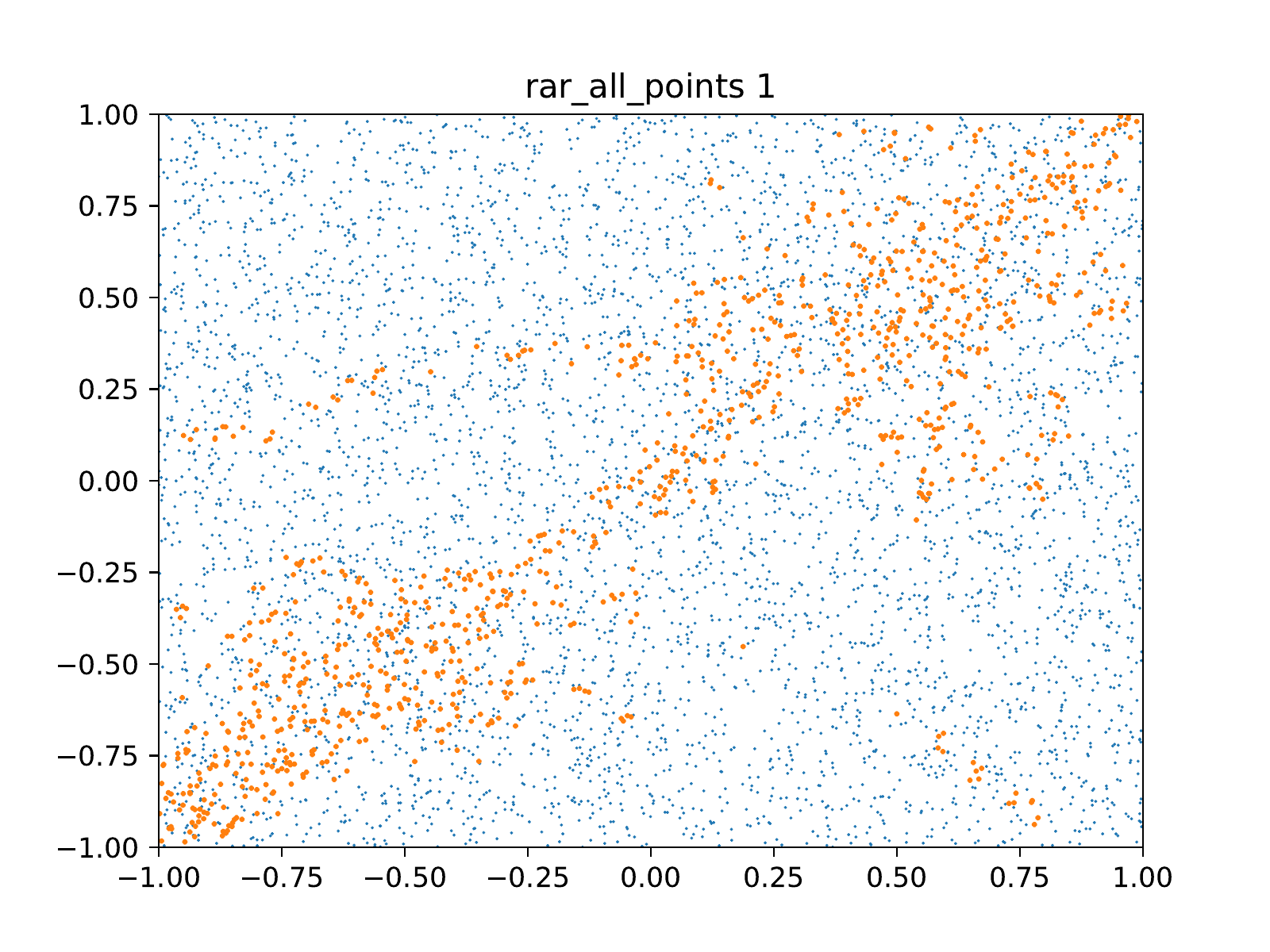}
        \end{overpic}
              \begin{overpic}[width=0.242\textwidth, clip=true,tics=10]{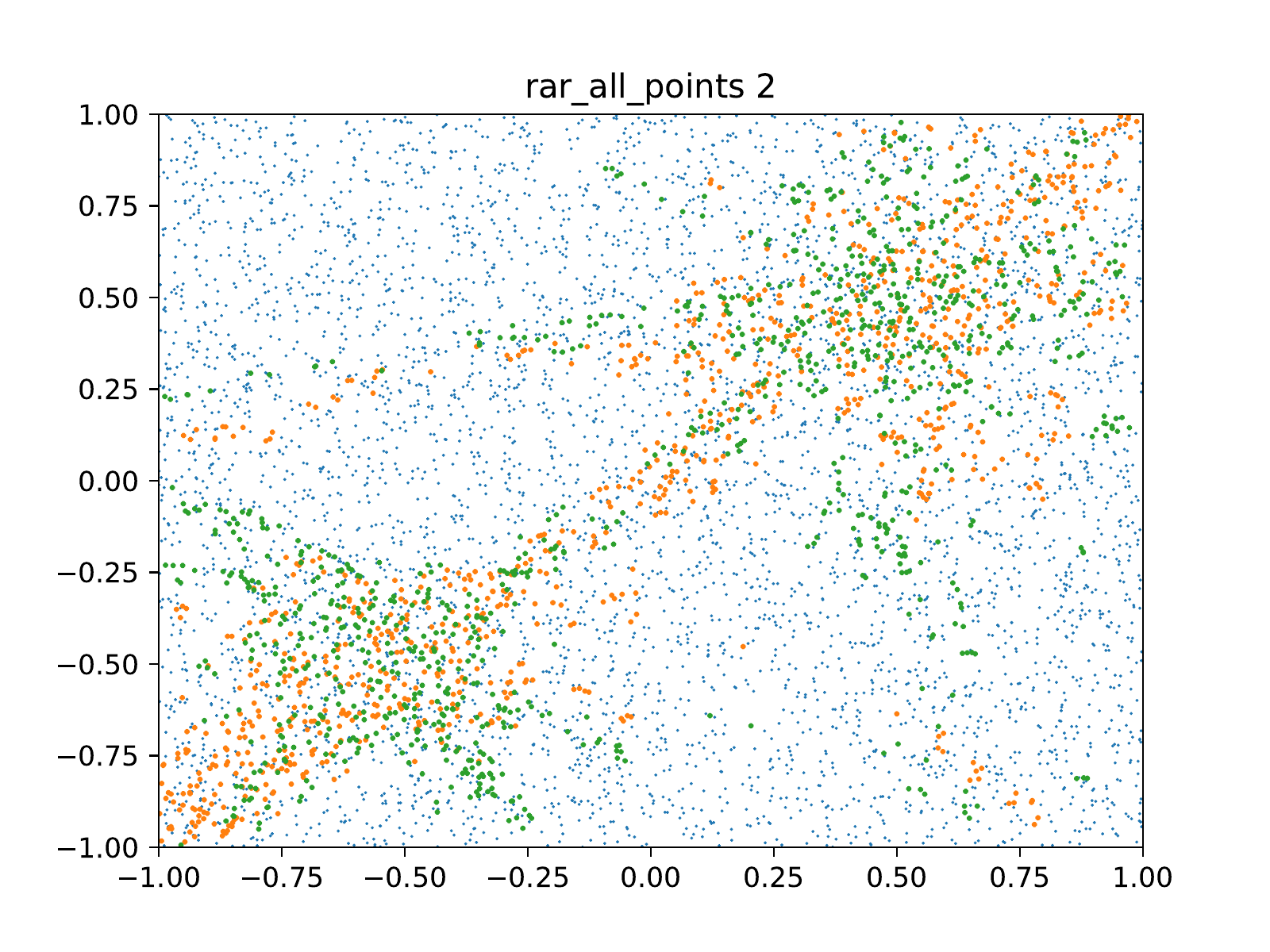}
        \end{overpic}
                      \begin{overpic}[width=0.242\textwidth, clip=true,tics=10]{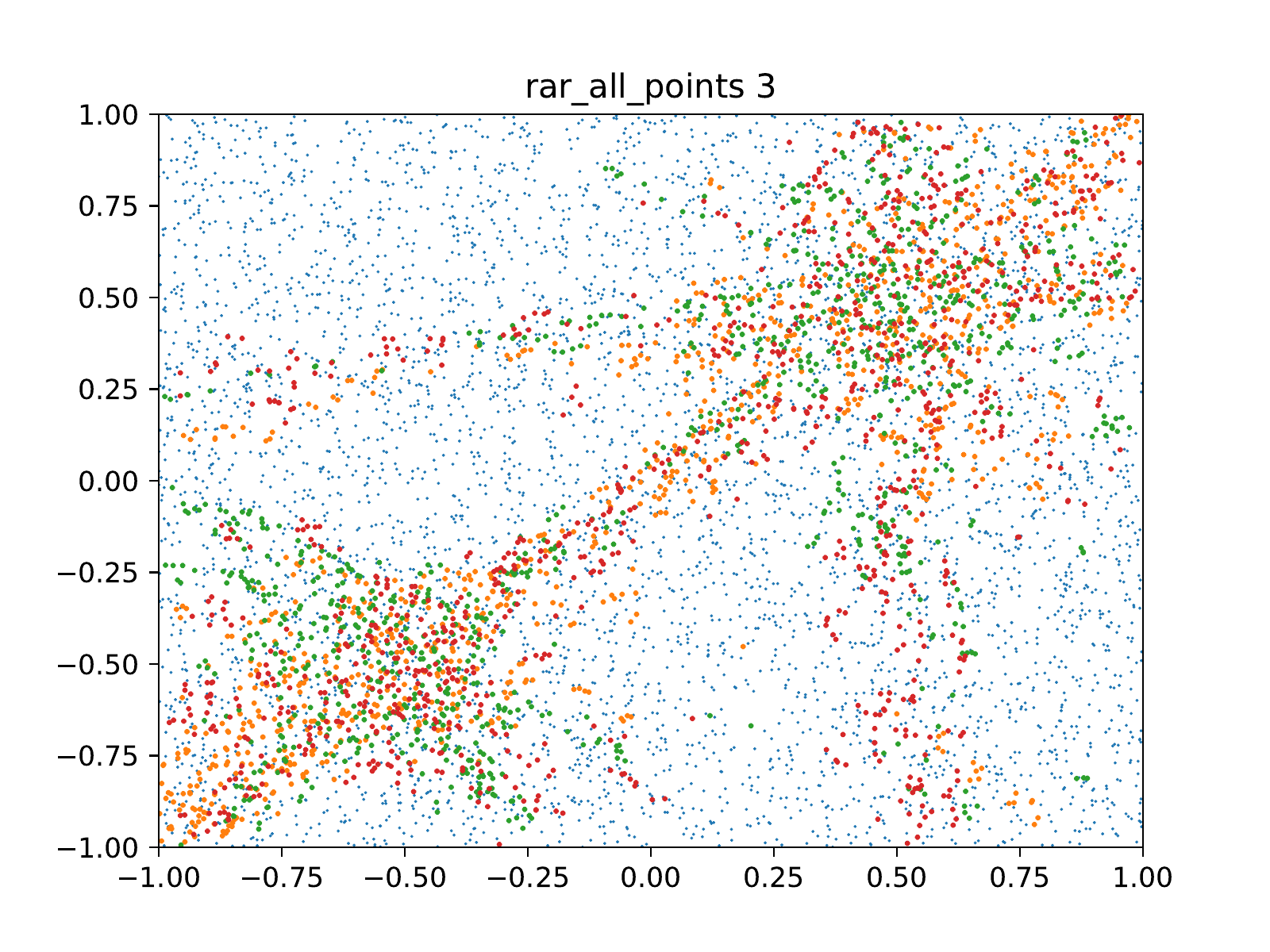}
        \end{overpic}
                \begin{overpic}[width=0.242\textwidth, clip=true,tics=10]{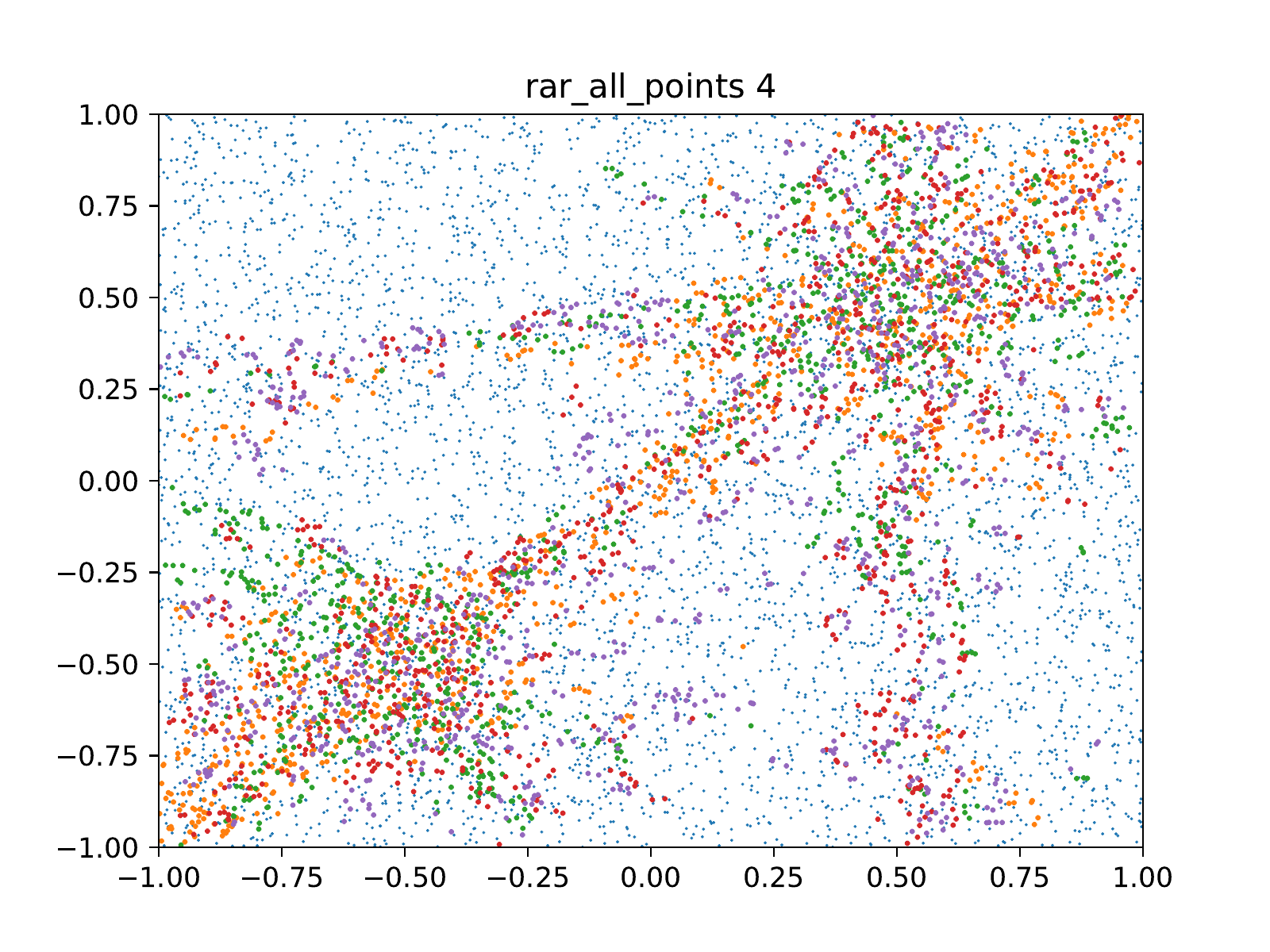}
        \end{overpic}
    \end{center}
    \begin{center}
           \begin{overpic}[width=0.242\textwidth, clip=true,tics=10]{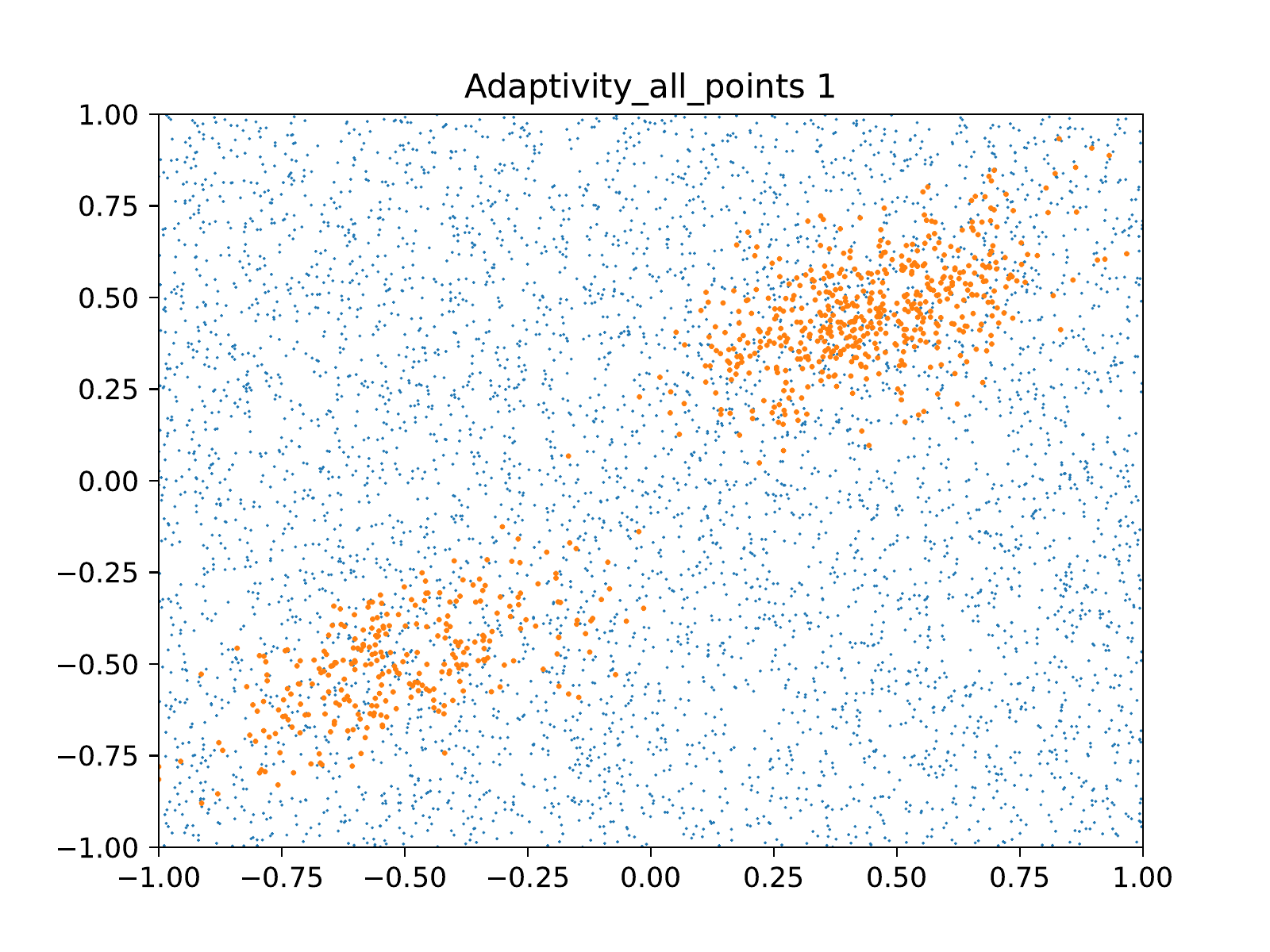}
        \end{overpic}
 \begin{overpic}[width=0.242\textwidth, clip=true,tics=10]{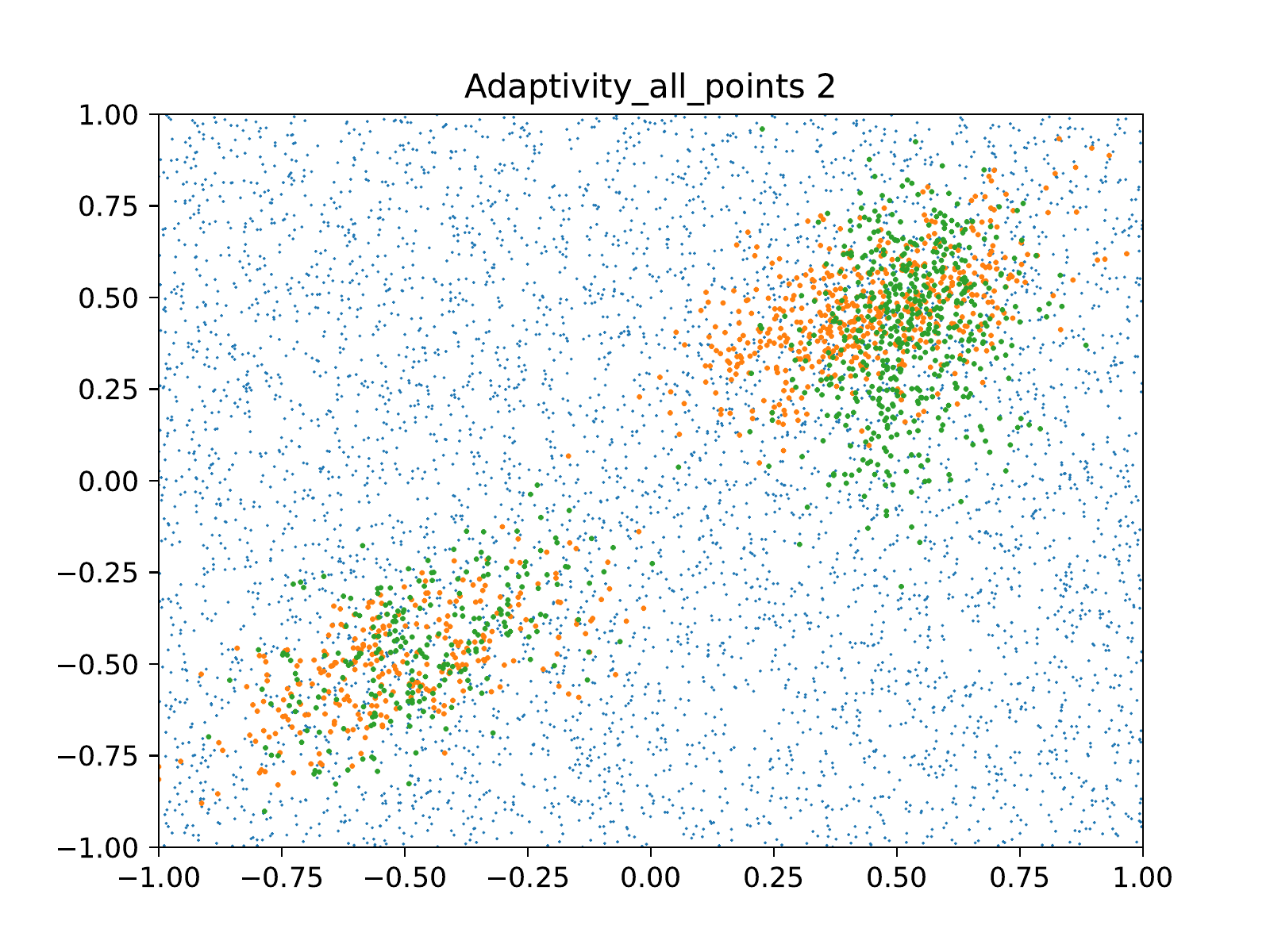}
        \end{overpic}
                \begin{overpic}[width=0.242\textwidth, clip=true,tics=10]{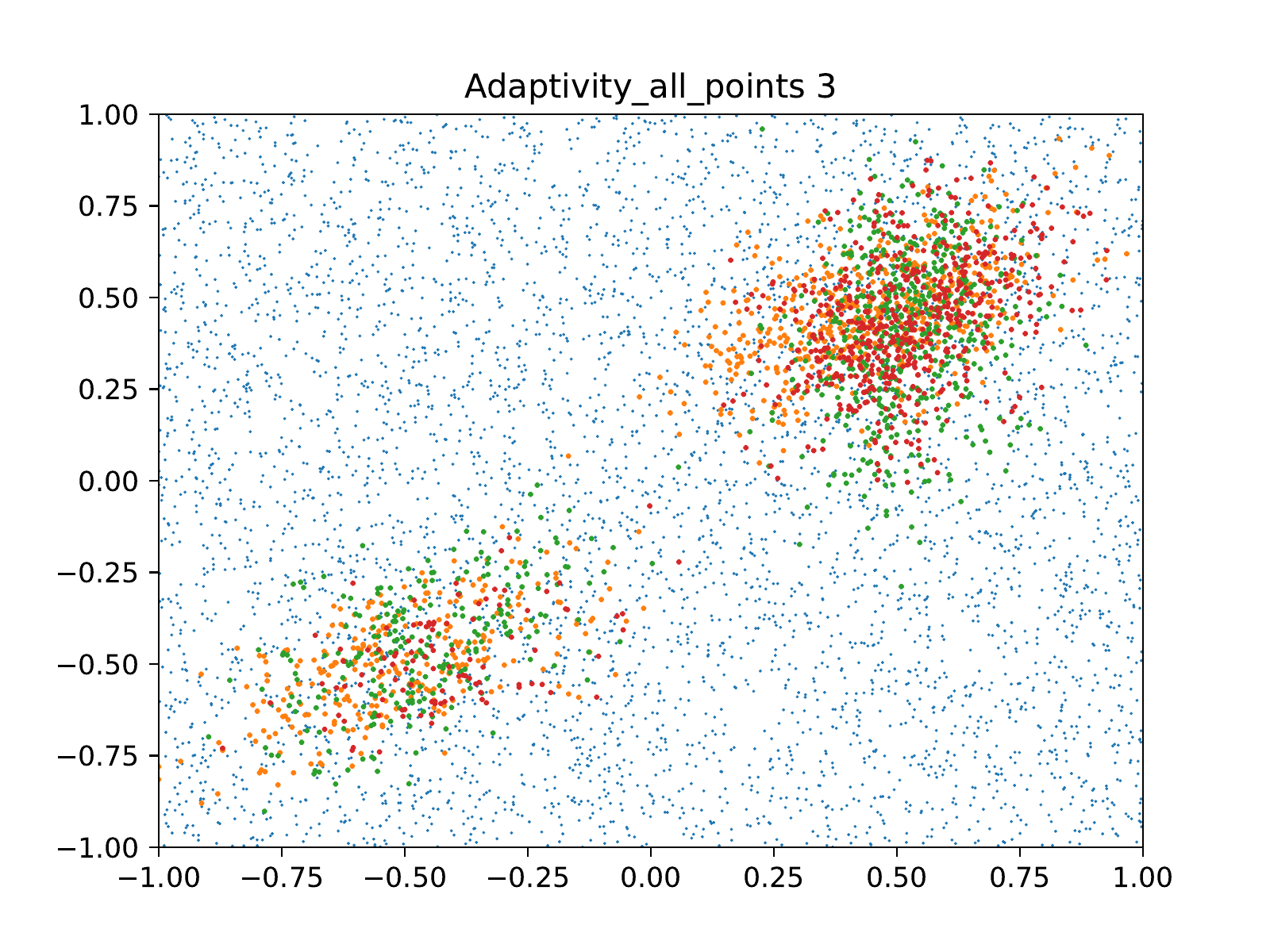}
        \end{overpic}
        \begin{overpic}[width=0.242\textwidth, clip=true,tics=10]{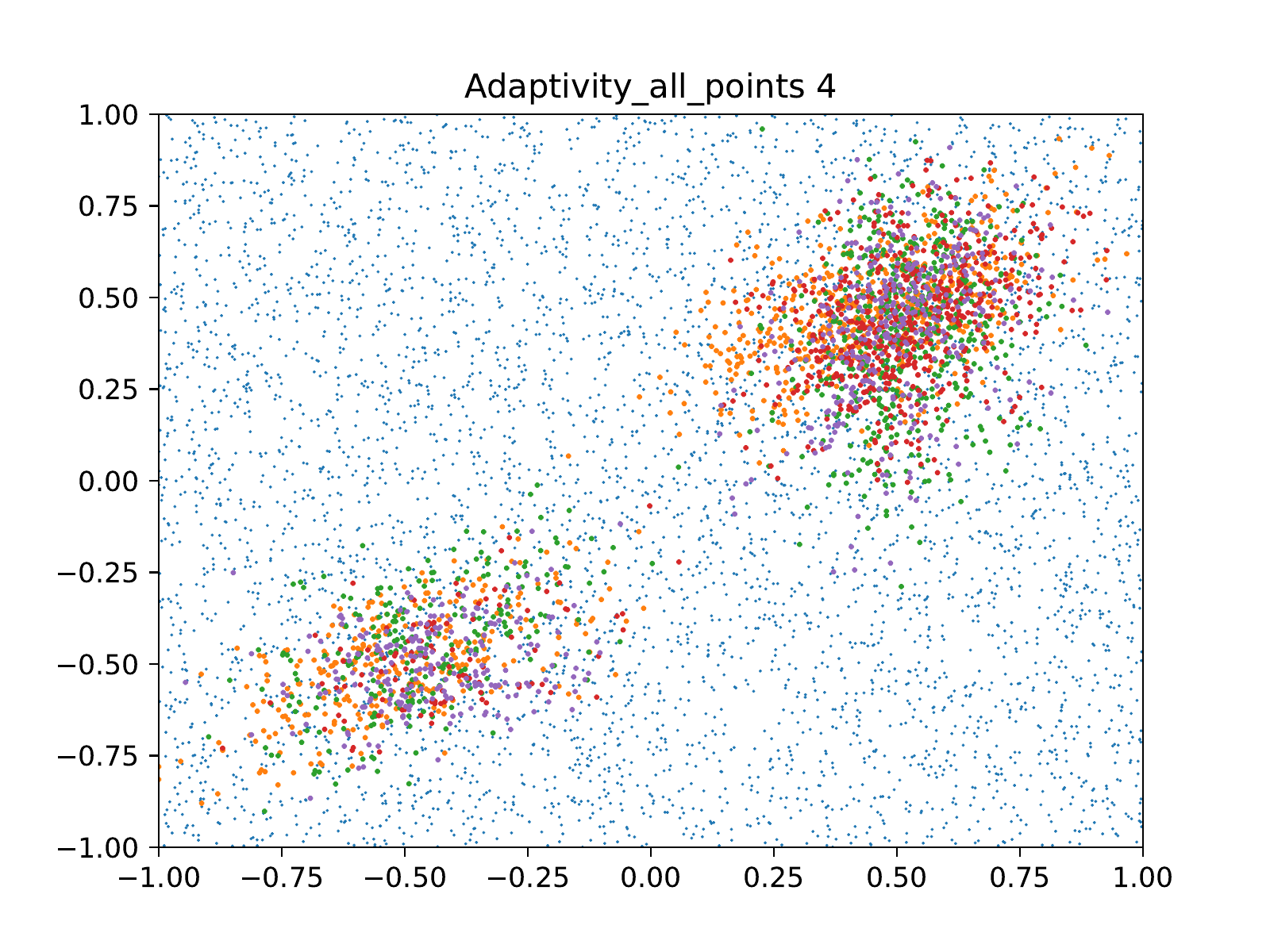}
        \end{overpic}
    \end{center}
    \caption{The first four updates of the  samples distribution  obtained using the three sampling strategies: Uniform, RAR and SAIS (from top to bottom). }
    \label{twopeak_samples}
\end{figure}

The corresponding numerical results are display in Figs.  \ref{two_peak_solution} and \ref{twopeak_samples}.  We observe very similar behavior for the predicted results as in the earlier experiments: the samples generated by the \textit{SAIS} strategy are concentrated more in the areas where the residual error is high, and the predicted values produced by the \textit{SAIS} strategy closely match the exact solution. This further demonstrates the adaptability of the FI-PINNs framework.

\subsection{Causal extension of FI-PINNs} \label{ACeq}

In this section, we present a causal extension of FI-PINNs to deal with time-dependent PDEs.   In comparison to FI-PINNs, causal FI-PINNs have two modifications: (1)  the {\it limit-state function} is defined using a causal formulation of the performance function, and (2) the network is trained using a causal weight residual loss. To this end, we divide the whole domain into $N_{t} $ subdomains in the time direction. Then we can define the weighted residual loss:
    \begin{equation}\label{wloss}
      \mathcal{L}_{c}(\boldsymbol{\theta})=\frac{1}{N_{t}} \sum_{i=1}^{N_{t}} w_{i} \mathcal{L}_{c}\left(t_{i}, \boldsymbol{\theta}\right),
    \end{equation}
    where $\mathcal{L}_{c}\left(t_{i}, \boldsymbol{\theta}\right)$ is the residual loss in the $i$th subdomain and the weights
    \begin{equation}
        w_{i}=\exp \left(-\epsilon \sum_{k=1}^{i-1} \mathcal{L}_{c}\left(t_{k}, \boldsymbol{\theta}\right)\right), \text { for } i=2,3, \ldots, N_{t},
    \end{equation}
where $\epsilon$ is a causality parameter.  Now  we can modify the performance function  of FI-PINNs as  $\mathcal{Q}(x) = w_{i}|r(x, t_{i})|,$ where $r(x, t_{i})$ is the residual function in the $i$th subdomain.  Using this performance function, we can define our FI-PINNs framework.  We train the network using the weighted residual loss of Eq. (\ref{wloss}), just like the causal PINNs\cite{wang2022respecting} do. The proposed method is now refereed to as causal FI-PINNs.

In order to  test the effectiveness of the causal FI-PINNs, we consider the following one-dimensional Allen-Cahn equation \cite{wang2022respecting}
    \begin{equation}
    \begin{array}{l}
        u_{t}-0.0001 u_{x x}+5 u^{3}-5 u=0, \quad t \in[0,1], x \in[-1,1], \\
        u(0, x)=x^{2} \cos (\pi x), \\
        u(t,-1)=u(t, 1), \\
        u_{x}(t,-1)=u_{x}(t, 1).
    \end{array}
    \end{equation}

We use the same network architectures as in \cite{wang2022respecting} for the NN model, which consists of  a  fully connected network with 6 hidden layers, each with 128 neurons, and  \texttt{tanh} as activation function. For ease of use, we also impose periodic boundary conditions as hard constraints.  The optimizer is the \texttt{Adam} with 20000  iterations.  For simplicity, we
set the  causality parameter $\epsilon =100$.  The initial number of collocation points is 1000. The failure probability tolerance $\epsilon_{p}$ is set to 0.001.
In the figures, we will use ``\textit{Uniform}" to denote the conventional causal PINNs, ``\textit{SAIS}" to denote the causal FI-PINNs algorithm.

\begin{figure}[t]
        \begin{center}
            \begin{overpic}[width = 0.9\textwidth]{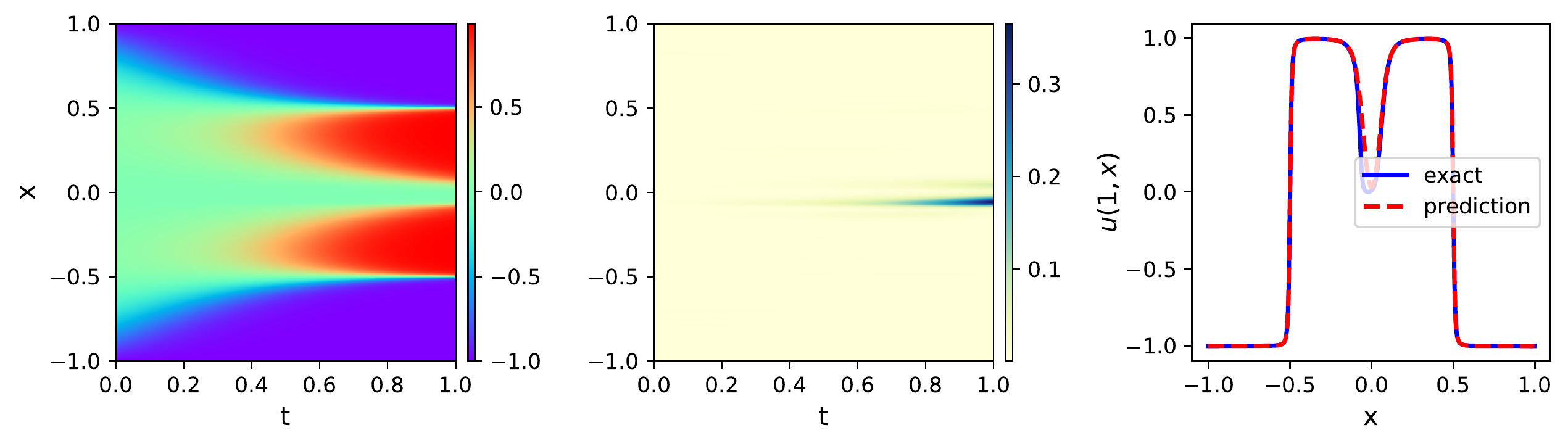}
                \put(35,29){\small Numerical results by \textit{Uniform}}
            \end{overpic}
        \end{center}
            \vspace{0.3cm}
        \begin{center}
            \begin{overpic}[width = 0.9\textwidth]{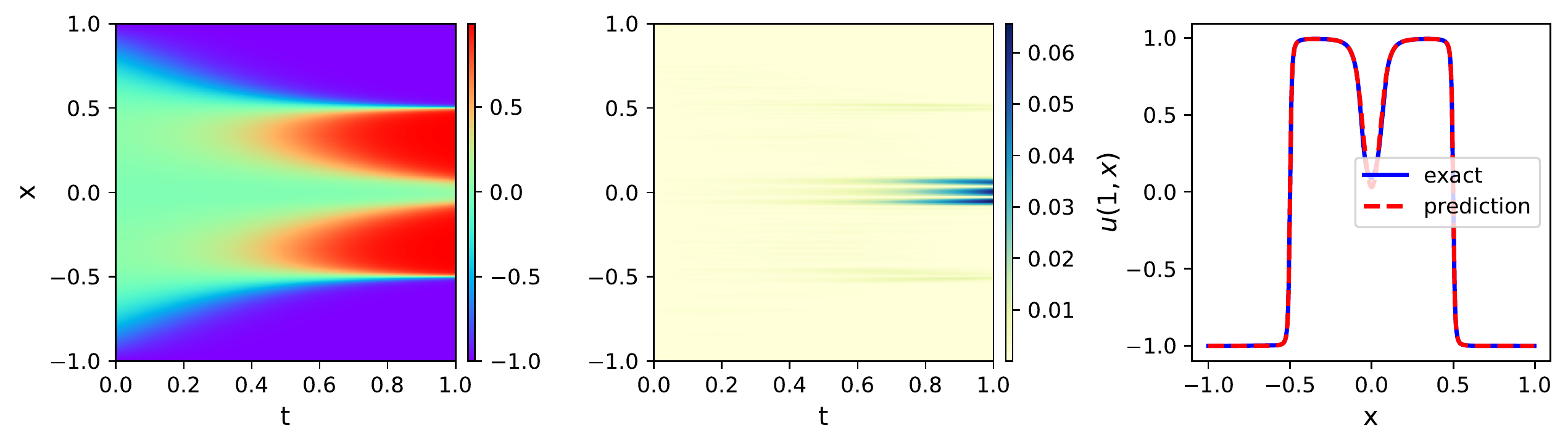}
            \put(37,29) {\small Numerical results by \textit{SAIS}}
            \end{overpic}
        \end{center}
        % \vspace{-0.3cm}
        \caption{Numerical results obtained by different methods.}
        \label{Ac_results}
    \end{figure}

\begin{figure}[htbp]
        \begin{center}
            \begin{overpic}[width = 0.7\textwidth]{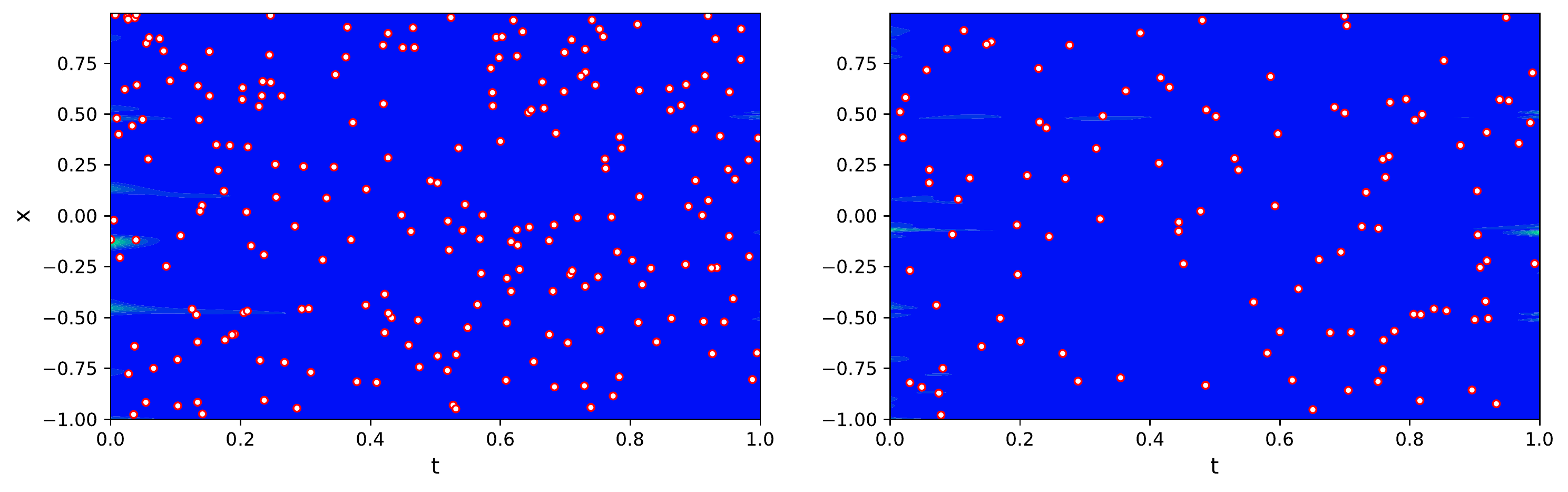}
                \put(35,33){\small Samples by \textit{Uniform}}
            \end{overpic}
        \end{center}
            \vspace{0.5cm}
        \begin{center}
            \begin{overpic}[width = 0.7\textwidth]{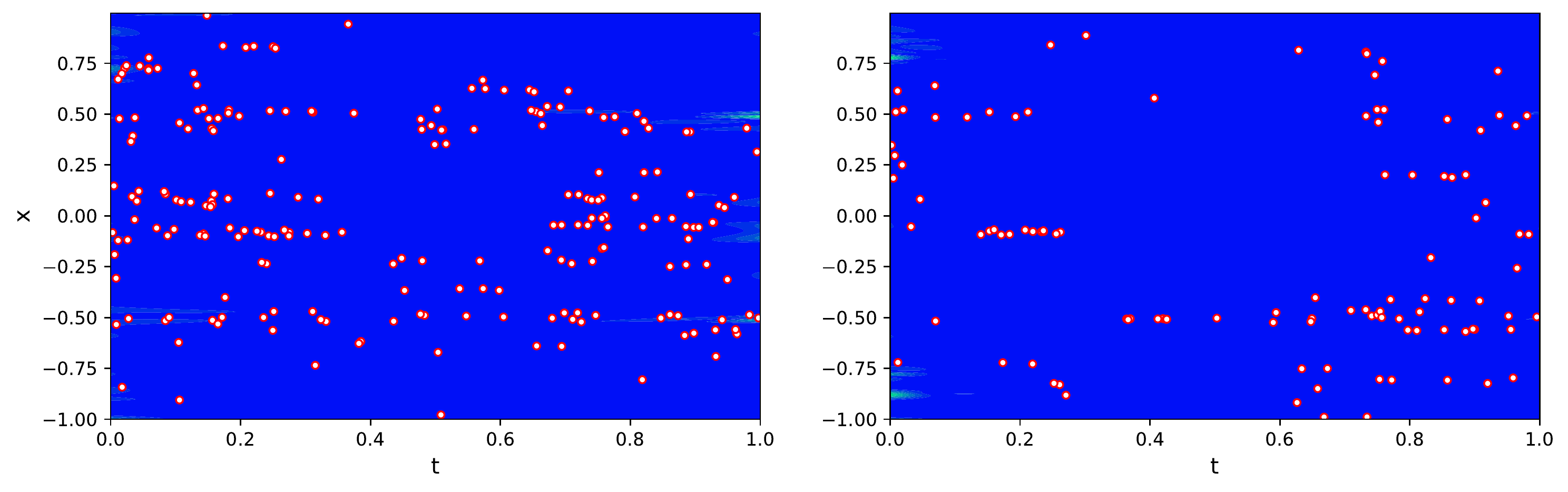}
            \put(40,33) {\small Samples by \textit{SAIS}}
            \end{overpic}
        \end{center}
        % \vspace{-0.3cm}
        \caption{The $3th, 5th$ sample distribution (from left to right) obtained by different methods.}
        \label{Ac_samples}
    \end{figure}

\begin{figure}[htbp]
        \begin{center}
            \begin{overpic}[width = 0.35\textwidth]
                {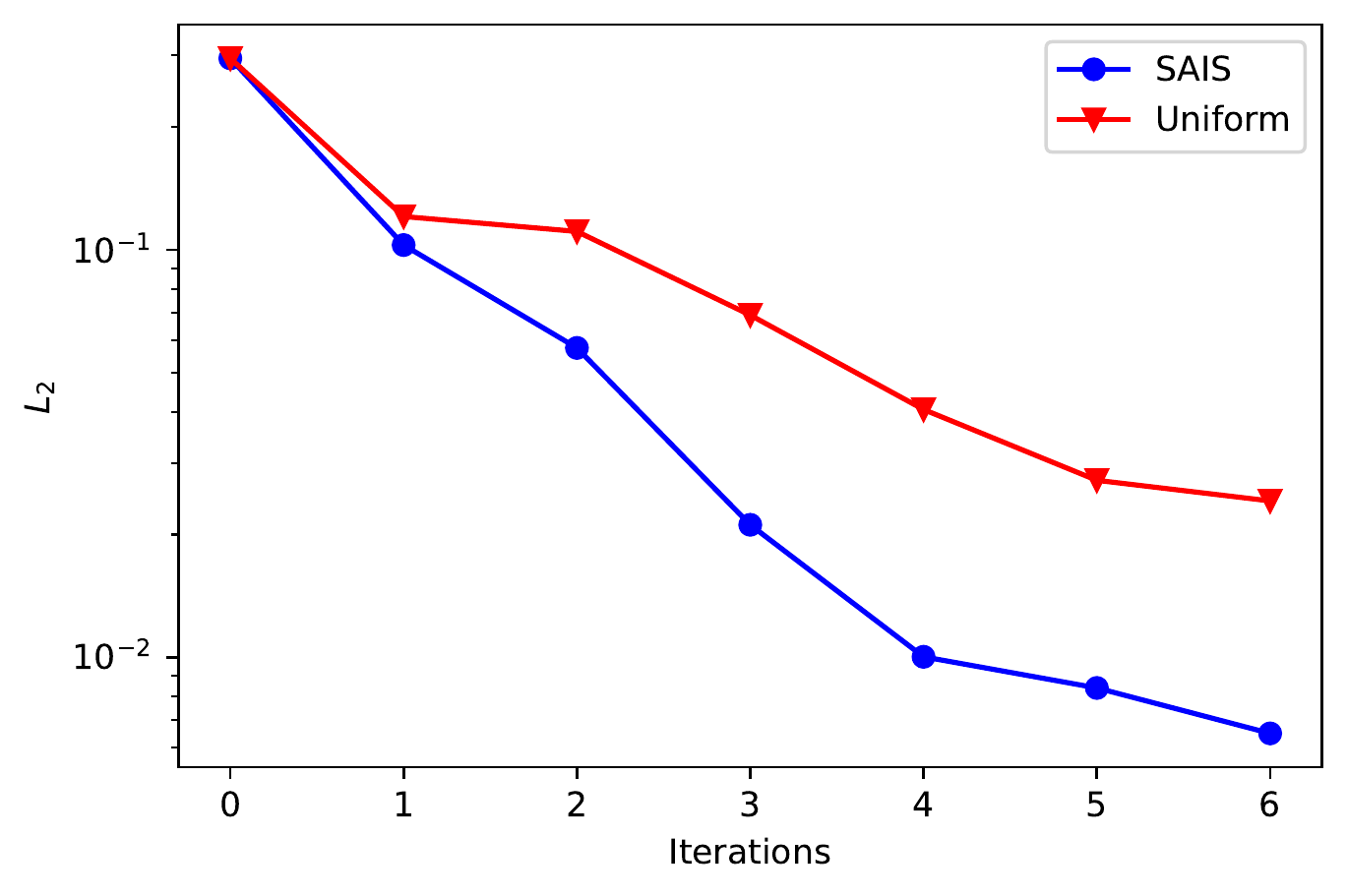}
            \end{overpic}
            \begin{overpic}[width = 0.365\textwidth]
                {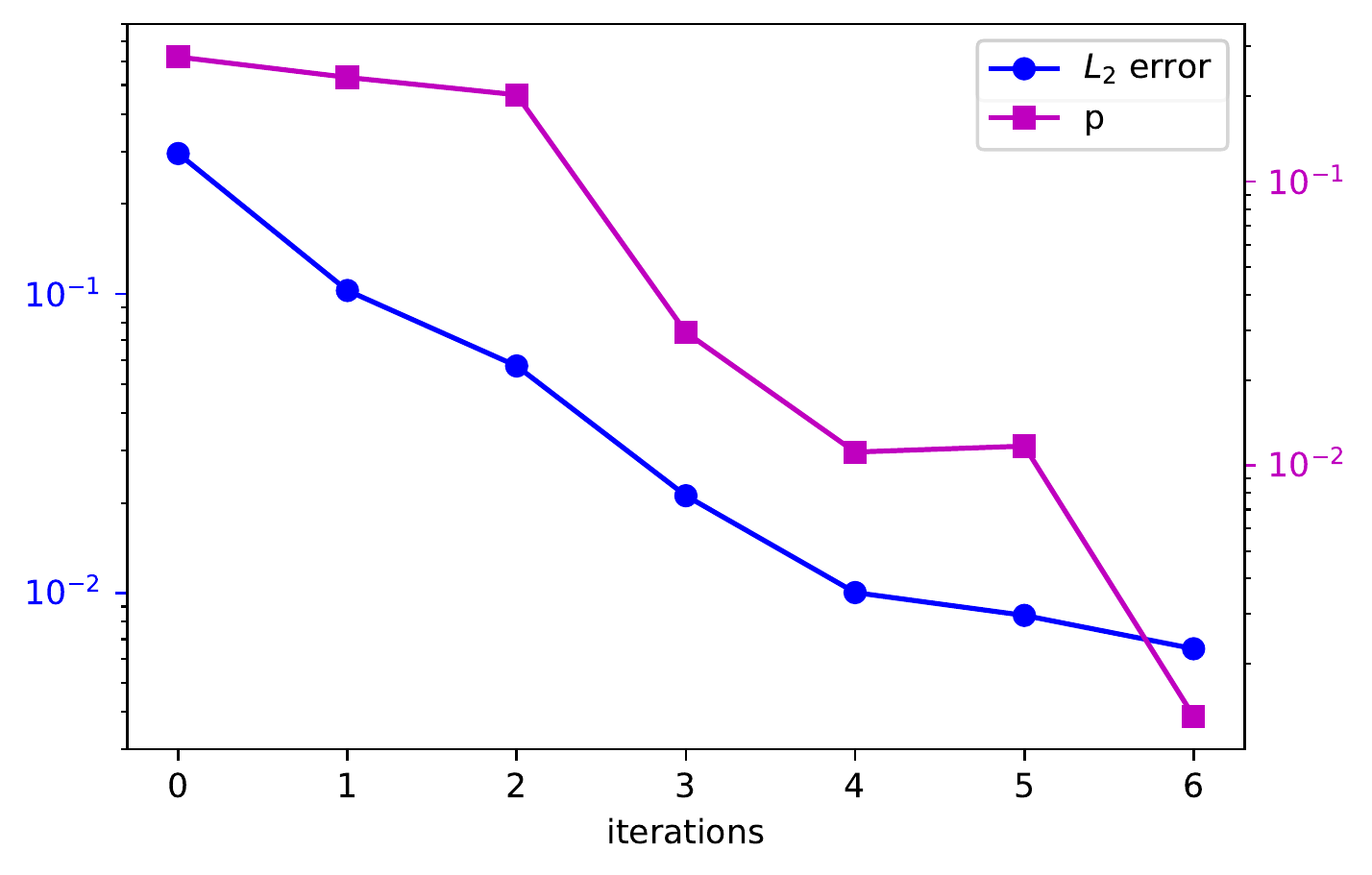}
            \end{overpic}
        \end{center}
        \caption{Relative $L_{2}$ error (left) and the corresponding failure probability (right).}
        \label{Ac_error}
    \end{figure}

The corresponding numerical results are display in Figs. \ref{Ac_results}- \ref{Ac_error}.  We also compare the predicted error obtained by the conventional FI-PINNs  using  various sampling techniques, as shown in Table \ref{Ac_compared_error}.  Our causal FI-PINNs strategy, which requires less than 4000 training points, can achieve  a  prediction error of only $6.5\times 10^{-3}$ after 6 updates.  The size of training dataset is much smaller than the one used in casual PINNs\cite{wang2022respecting}, which achieves the same order of accuracy and save computational cost at the same time.  As can be seen, our casual FI-PINNs produce significantly lower predicted error than the conventional casual PINNs when given the same number of training dataset. The relative error was found to be $2.42\times 10^{-2}$ for causal PINNs strategy in this case.   Uniform training points were used in the original casual PINNs sampling strategy, which resulted in worse predicted error. It should be noted that, in this example, FI-PINNs does not improve numerical results as the {\it principle of causality} for solving time-dependent PDEs\cite{wang2022respecting}.

    \begin{table}
      \caption{The final predicted error achieved by different methods.}
    \label{Ac_compared_error}
    \renewcommand\arraystretch{1.3}
    \centering
    \begin{tabular}{ccc}
        \toprule
      Method& Relative $L_{2}$ error\\
      \hline
      PINNs& $6.79\times 10^{-1}$\\
      FI-PINNs & $5.03\times 10^{-1}$ \\
      casual PINNs & $2.42\times 10^{-2}$ \\
      casual FI-PINNs & $\mathbf{6.50\times 10^{-3}}$\\
      \bottomrule
    \end{tabular}
    \end{table}

%\bibliographystyle{siam}
%\bibliography{refer}

\end{document}